\documentclass[11pt]{article}

\usepackage{latexsym,amsmath,amsfonts,amssymb,amstext,
euscript,array,enumerate}
\usepackage{color}
\usepackage{comment}

\usepackage[colorlinks=true,linkcolor=blue,urlcolor=blue]{hyperref}

\usepackage{amsthm}
\usepackage{bm}
\allowdisplaybreaks

\numberwithin{equation}{section}

\theoremstyle{plain}
\newtheorem{theorem}{Theorem}[section]
\newtheorem{lemma}[theorem]{Lemma}
\newtheorem{proposition}[theorem]{Proposition}

\newtheorem{definition}[theorem]{Definition} 

\numberwithin{equation}{section}

\theoremstyle{Remark}
\newtheorem{remark}[theorem]{Remark}

\def\esssup_#1{\underset{#1}{\mathrm{ess\,sup\, }}}
\def\essinf_#1{\underset{#1}{\mathrm{ess\,inf\, }}}

\def\qed{{\hfill\hbox{\enspace${ \square}$}} \smallskip}
\def\sqr#1#2{{\vcenter{\vbox{\hrule height .#2pt \hbox{\vrule
 width .#2pt height#1pt \kern#1pt \vrule
width .#2pt} \hrule height .#2pt}}}}
\def\square{\mathchoice\sqr54\sqr54\sqr{4.1}3\sqr{3.5}3}

\def\ds{\begin{displaystyle}}
\def\eds{\end{displaystyle}}

\def\<{\langle }
\def\>{\rangle }

\def \N{\mathbb{N}}
\def \R{\mathbb{R}}

\def \E{\mathbb{E}}
\def \F{\mathbb{F}}

\def \P{\mathbb{P}}

\def\beqs{\begin{eqnarray*}}
\def\enqs{\end{eqnarray*}}
\def\beq{\begin{eqnarray}}
\def\enq{\end{eqnarray}}

\newcommand{\spernutxa}[1]{\mathbb{E}^{t,x,a}_{\nu} \left[ #1 \right]}                               

\newcommand{\spertxa}[1]{\mathbb{E}^{t,x,a} \left[ #1 \right]}                               

\DeclareMathAlphabet{\mathonebb}{U}{bbold}{m}{n}                           %
\newcommand{\one}{\ensuremath{\mathonebb{1}}}                               

\allowdisplaybreaks

\title{ Optimal control of infinite-dimensional \\Piecewise Deterministic Markov Processes: a BSDE approach. Application to the control of an excitable cell membrane.}

\author{Elena BANDINI\thanks{Dipartimento di Matematica e Applicazioni, Universit\`a di Milano-Bicocca, Via R. Cozzi 55, 20125 Milano, Italy; e-mail: \texttt{elena.bandini@unimib.it}}
\and
Mich\`ele THIEULLEN \thanks{Laboratoire de Probabilit\'es, Statistique et Mod\'elisation (LPSM, UMR 8001), Sorbonne Universit\'e - Campus Pierre et Marie Curie, Boite 158, 4 Place Jussieu, 75252 Paris Cedex 05, France; e-mail: \texttt{michele.thieullen@upmc.fr}}
}

\begin{document}

\date{}
\maketitle

\begin{abstract}
In this paper we consider the optimal control of Hilbert space-valued infinite-dimensional Piecewise Deterministic Markov Processes (PDMP) 
and we prove that the corresponding value function can be represented via a Feynman-Kac type formula
through the solution of a constrained Backward Stochastic Differential Equation. 
A fundamental step consists in showing 
that the corresponding integro-differential Hamilton-Jacobi-Bellman equation
 has a unique viscosity solution, by proving a suitable comparison theorem. We apply our results
  to the control of a PDMP Hodgkin-Huxley model with spatial component, previously 
  studied 
   in \cite{ReTreThieu}, \cite{Ren} and inspired by optogenetics.

\end{abstract}

\noindent{\small\textbf{Keywords:} infinite-dimensional PDMPs, 
 constrained backward stochastic differential equations, 
integro-differential Hamilton-Jacobi-Bellman equation, 
viscosity solutions in infinite dimensions,
 spatio-temporal Hodgkin-Huxley models.
}

\medskip

\noindent{\small\textbf{MSC 2010:} 93E20, 60H10, 60J25.}

 \section{Introduction}\label{Sec_introduction}

In this  paper we consider  optimal control problems for  Hilbert space-valued infinite-dimensional Piecewise Deterministic Markov Processes, 
and we prove that the corresponding value function can be represented through a Feynman-Kac formula by means of the solution of a constrained Backward Stochastic Differential Equation (BSDE). As an intermediate step, we also show that the corresponding   Hamilton-Jacobi-Bellmann (HJB)   has a unique viscosity solution by providing a  comparison theorem for suitable Integro Partial Differential Equations (IPDE). We apply our theoretical results to the control of a PDMP Hodgkin-Huxley model with spatial component,  previously considered in \cite{ReTreThieu}, \cite{Ren} and inspired by optogenetics. 

The Feynman-Kac type representation for the value function is obtained by  implementing the randomization procedure introduced in \cite{KhPh} for jump-diffusions, later extended in \cite{BandiniFuhrman} and \cite{BandiniPDMPsNoBordo} respectively to the case of  finite-dimensional pure jump Markov processes and of finite-dimensional PDMPs. 
The control randomization method is particularly useful to probabilistically represent 
the value function associated to stochastic control problems,  where the laws of the family of controlled processes are not dominated by a common measure. Roughly speaking, the randomization principle consists in enlarging the state space by an additional independent piecewise constant component corresponding to the control, and in subsequently  generating a family of dominated laws and an auxiliary control problem, where the cost is optimized with respect to the intensity of the extended  pure jump component. The value function of this latter 
(randomized) control problem can be represented by means of the solution of a constrained BSDE,  namely a backward equation driven by a random measure with a sign constraint on its martingale part. In order to be able to relate this backward equation to the HJB equation associated to the primal problem, one has 
to show that the randomized value function does not depend on the additional component, and that it provides a solution to the above-mentioned  HJB equation. 
Afterwards, the Feynman-Kac representation formula for the original value function comes  from the uniqueness of the viscosity solution to the corresponding HJB equation. 
 We refer the reader to the introduction of \cite{KhPh} for an extended  exposition of the issues involved. 
 Note that the randomization procedure is a very  general methodology which  applies even if the laws of the controlled processes are  dominated.
The Feynman-Kac representation formula can be used to design algorithms based on the numerical approximation of the solution to
the corresponding constrained BSDE, and therefore to get probabilistic numerical approximations for the value
function of the addressed optimal control problem, see e.g. \cite{KLP}.

In our infinite-dimensional setting,  we prove existence and uniqueness (in a suitable sense) of the solution of such a constrained BSDE and its independence with respect to the additional component. We also provide a randomized dynamic principle which enables us to establish that the value function of the randomized problem  is a viscosity solution of
 the  
HJB-IPDE 
 on the Hilbert space.
Viscosity solutions for partial differential equations in infinite dimension  with unbounded linear terms have been first studied in  \cite{CrandallLionsV}, where  the notions of 
$B$-upper/lower-semicontinuity are introduced, and subsequently considered by many other authors, see e.g. \cite{gozswi15} for a modern and detailed exposition on this topic.  Recently the papers \cite{SwiechZabczykEXISTENCE} and \cite{SwiechZabczykUNIQUENESS} have addressed respectively existence and uniqueness for an 
HJB-IPDE 
 resulting from the control of an Hilbert space-valued SDE driven by a L\'evy process. 
 Notice that in our framework 
  we
  do not  ask that our PDMP is a strong solution to some SDE.
   Our approach is instead based on  the study of the local characteristics of the PDMP in the spirit of the theory developed in  \cite{Da}.
We prove a comparison theorem which implies  the uniqueness of the viscosity solution of our 
 HJB-IPDE. 
 The appropriate definition of viscosity solution, on which the  comparison theorem
   relies,  is derived suitably extending  the one provided in    \cite{SwiechZabczykUNIQUENESS}.

Our theoretical results are applied to the control of a PDMP Hodgkin-Huxley model with spatial component. Hilbert space-valued PDMP models describing the spatio-temporal evolution of a neuron with a finite number of ion channels (or more general excitable membranes) have been rigorously settled in \cite{Riedler2011}. In particular it was proved in \cite{WainribRiedler} that such PDMP converge to the spatio-temporal Hodgkin-Huxley model proposed in \cite{HH} when the number of channels goes to infinity, see also \cite{Austin}.   Inspired by optogenetics, optimal control  of general infinite-dimensional PDMP has been previously considered in \cite{ReTreThieu}, \cite{Ren}. In particular the results in \cite{ReTreThieu} were applied to a tracking problem for a Hilbert space-valued Hodgkin-Huxley type PDMP.  In that paper, as in the present one, piecewise open loop controls (see e.g. \cite{Ver}) were considered, and the control acted on the three characteristics of the PDMP. However, the main tools were relaxed controls and the optimal control theory  of Markov Decision Processes, see \cite{BeSh}. 
Moreover, even if an 
HJB-IPDE 
  were written down in that paper, no study was conducted about existence or uniqueness of its solutions.   We also mention the more recent paper \cite{CalviaInf}, which exploits Markov Decision Processes in infinite dimension in the framework of stochastic filtering.

Many generalizations of the present work may be possible. 
For instance, it would be interesting to treat the general case with infinite-dimensional PDMPs  on a  state space with  boundary,  from which   additional instantaneous jumps  into the interior of the domain may occur (in the finite-dimensional case, this feature has been recently considered in \cite{BandiniPDMPsBordo}).
Moreover,  in our application section  we have considered the classical case of a Laplacian operator,  but other  operators could be addressed as well. Finally, a challenging future development 
would consists in applying our results 
to the infinite-dimensional PDMP that naturally arise in filtering problems.

The  paper is organized as follows. In Section \ref{Sec:PDP_Sec_control_problem} we  construct our infinite-dimensional controlled  
PDMP and we define the related optimal control problem.
In particular, inspired by \cite{J}, we provide a canonical construction  of the PDMP state process in infinite dimension, by suitably extending the finite-dimensional construction
 implemented  in \cite{BandiniPDMPsNoBordo},  \cite{BandiniFuhrman}.
  We then set the associated control problem, 
  and we establish in Theorem \ref{Sec:PDP_Thm_existence} 
   that the corresponding value function is a viscosity solution of the   
HJB equation \eqref{Sec:PDP_HJB}-\eqref{Sec:PDP_HJB_T}.
In Section \ref{Sec:PDP_Section_dual_control} we describe the control randomization method in our setting,  and we introduce  the randomized optimal control problem. Then in Section \ref{Sec:PDP_Sec_ConstrainedBSDE} we  define and  study the related constrained BSDE, and we address
the Feynman-Kac representation. 
  As described above,  the first step of the randomization approach consists in proving that the randomized value function does not depend on the additional component,  and that satisfies a suitable randomized dynamic programming principle, see respectively  Proposition \ref{Sec:Prop_Feynman_Kac_HJB} and Theorem \ref{Sec:PDP_THm_Feynman_Kac_HJB}. 
Then in  Theorem \ref{Sec:first_main result} we show  that also the randomized value function is a viscosity solution to the HJB equation.
The last step towards the Feynman-Kac representation  
consists in   the comparison Theorem \ref{Sec:PDP_Thm_uniqueness}, which provides uniqueness of the viscosity solutions to our HJB-IPDE equation. Section \ref{example} is devoted to the application of our results 
  to the control of a spatio-temporal Hodgkin-Huxley type model. Finally,  Sections \ref{Sec_mainproofs_Sec2} and \ref{Sec_mainproofs_Sec4} are devoted to the proofs of the results provided respectively in Sections \ref{Sec:PDP_Sec_control_problem} and \ref{Sec:PDP_Sec_ConstrainedBSDE}.

  \section{Optimal  control of infinite-dimensional PDMPs}\label{Sec:PDP_Sec_control_problem}
 In the present section we are going to formulate
 an  optimal control problem for infinite-dimensional piecewise deterministic Markov processes, and to discuss its solvability.
The PDMP state space $E$ is a real separable Hilbert space, equipped with the norm $||\cdot||$ and the inner product $\langle \cdot, \cdot \rangle$, with  corresponding Borel  $\sigma$-field $\mathcal E$. 
  In addition,  we introduce  a Polish space $A$, endowed with its Borel $\sigma$-field $\mathcal{A}$, called the space of control actions.
 The other data of the problem
 consist in four functions $f$,   $b$,  $\lambda$ on $E \times A$, $g$ on $E$, in a  probability transition kernel $Q$ from $(E \times A, \mathcal E \otimes \mathcal A)$ to $(E, \mathcal E)$, and in an operator $L$, satisfying the following conditions.

  \medskip
 
  \noindent \textbf{(HL)}

  \medskip 
  
\noindent (i) $L$ is  a   linear, densely defined, maximal monotone operator in $E$. Moreover, there exists an operator $B$ on $E$  bounded, linear, positive (i.e., $\langle Bx,\,x \rangle > 0$ for every $x \in E$, $x \neq 0$) and self-adjoint, such that $L^\ast B$ is bounded on $E$, and,  for some $c_0\geq 0$,   
 \begin{equation}\label{LstarB_property}
 \langle (L^\ast B + c_0 B)\,x, \,x\rangle \geq 0 \quad \forall x \in E.
 \end{equation}
 We define the space $E_{-1}$ to be the completion of $E$ under the norm 
 $
 ||x||_{-1} = ||B^{1/2} x||.
 $
 $E_{-1}$ is an Hilbert space equipped with the inner product  
 $
 \langle x, x\rangle_{-1} = \langle B^{1/2} x, B^{1/2} x \rangle.
 $
Moreover,  
 \begin{equation}\label{-1norm}
 	||x||_{-1} \leq ||B^{1/2}||\,||x||, \quad x \in E.
 \end{equation}

\noindent (ii) $-L$ generates a strongly continuous semigroup $(S(u))_{u \geq 0}$ such that, for any $u>0$, $S(u)$ is a contraction on $E$ with respect to $||\cdot||_{-1}$. 
\begin{remark}
$-L$ is the generator of a strongly continuous semigroup of contractions $(S(u))_{u \geq 0}$ on $E$, see e.g. Theorem B.45 in  \cite{gozswi15}.
\end{remark}
 \begin{definition}\label{D:Bcontinuity}
 	We say that a function $u: W  \rightarrow \R$ is $B$-upper-semicontinuous (resp., $B$-lower-semicontinuous) on $W \subset [0,\,T] \times E$ if, whenever $t_n \rightarrow t$, $x_n \rightharpoonup  x$, $B\, x_n \rightarrow B\,x$, $(t,x) \in W$, then $\limsup_{n \rightarrow \infty}u(t_n,x_n) \leq u(t,x)$ (resp. $\liminf_{n \rightarrow \infty}u(t_n,x_n) \geq u(t,x)$). The function $u$ is $B$-continuous on $W$ if it is $B$-upper-semicontinuous and $B$-lower-semicontinuous on $W$.
 \end{definition}

\noindent  In the assumptions below $C$ is a generic constant which may vary from line to line.

 \medskip

 \noindent \textbf{(H$\textup{b$\lambda$Q}$)}

  \medskip

\noindent 
(i) $b: E \times A \mapsto E$, $\lambda: E \times A \mapsto R_+$ are bounded    continuous functions such that    
 	\begin{equation*}
 	\left\{
 	\begin{array}{ll}
 	||b(x,a)-b(x',a)|| \leqslant C \,||x-x'||_{-1},\quad x,\,x' \in E, a \in A\\
 |\lambda(x,a)-\lambda(x',a)| \leqslant C \,||x-x'||_{-1}, \quad x,\,x' \in E, a \in A.
 	\end{array}
 	\right.  	
 	\end{equation*}
 \noindent	
 (ii)
 	$Q$ maps 
 	$E \times A$ 
 	into the set of probability measures on $(E, \mathcal E)$, and is a
 	continuous
 	stochastic  kernel
 	(see e.g. Proposition 7.30 in \cite{BeSh}).
 	Moreover, 
for any real function $\varphi$ continuous on $(\varepsilon,\,T-\varepsilon)\times E$ for any $\varepsilon >0$ and   bounded,  and for every $R >0$, 
 we have, for all  $s, s' \in (\varepsilon,\,T-\varepsilon)$, 
\begin{align}
 &\left|\int_{E} \varphi(s,y) Q(z,a,dy)-\int_{E} \varphi(s,y) Q(z',a,dy)\right| \leq C\,  \omega(||z-z'||_{-1}),\quad z, z' \in E, a \in A\label {hp_Q},\\
 & \left|\int_{E} [\varphi(s,y) -\varphi(s',y)] \, Q(z,a,dy)\right| \leq C \sigma_R(|s-s'|), \quad z \in E: ||z|| \leq R, \,\,a \in A.\label {hp2_Q}
\end{align}
where $\omega$ is a  modulus of continuity, and  $\sigma_R( \cdot)$ is a  modulus of continuity depending on $R$.
 \medskip

 \noindent \textbf{(H$\textup{fg}$)}\quad
 $f: E \times A \mapsto \R_+$,  $g: E  \mapsto \R_+$ are 
   continuous and bounded functions, such that 
\begin{align*}
|f(x,a)-f(x',a)| + |g(x)-g(x')|&\leqslant C \,\omega(||x-x'||_{-1}),\quad    a \in A,
\end{align*}
 for all $x,\,x' \in E$, where $\omega$ is a modulus of continuity.

\subsection{The 
 optimal control problem}

We construct the controlled  process $X$ in a canonical way. 
We start by fixing $(t,x) \in [0,\,T] \times E$,
and we set $\Omega = [0,\,T] \times E \times \Omega'$, where $\omega = (t, x, \omega')$, $\Omega'^t$ being  the set of sequences  $\omega' = (t_n, e_n)_{n \geq 1}$ contained in $((0,\infty)\times E \cup \{(\infty, \Delta)\})$, where $\Delta \notin E$ is an isolated point adjoined to $E$, such that $t_n \leq t_{n+1}$, and $t_n < t_{n+1}$ if $t_n < \infty$. 
On the sample space $\Omega$ we define the canonical functions $T_n: \Omega \rightarrow (t,\,\infty]$, $E_n : \Omega \rightarrow E \cup \{\Delta \}$ as follows: $T_0(\omega) = t$,  $E_0(\omega) = x$, and for $n \geq 1$, $T_n(\omega)= t_n$, $E_n(\omega) = e_n$, and $T_\infty(\omega) = \lim_{n \rightarrow \infty}t_n$. We also introduce the associated integer-valued counting measure on $(0,\,\infty)\times E$ given by $p(ds\,dy)= \sum_{n \in \N} \delta_{(T_n, E_n)} (ds,dy)$.

The class of admissible control law $\mathcal A^t_{ad}$ is  the set of all Borel-measurable maps $\alpha: [t,\,\infty) \times E \rightarrow A$
of  the form:
 \begin{equation}\label{Sec:PDP_open_loop_controls}
 \alpha_s=\alpha_0(s-t,x)\,\one_{[t,\,T_1)}(s)+ \sum_{n=1}^{\infty}\alpha_n(s-T_n,E_n)\,\one_{[T_n,\,T_{n+1})}(s),\quad s \in [t,\,T],
 \end{equation}
 where $(\alpha_n)_n$, $\alpha_n: \R_+ \times E \rightarrow A$, is a sequence of measurable functions,  see for instance \cite{Da},
 \cite{CoDu}, \cite{Almudevar}.
 In other words, at each jump time $T_n$, we choose an open loop control $\alpha_n$ depending on the initial
 condition 
 $E_n$ and on the time elapsed up to  $T_n$, to be used until the next jump time.
 We define the controlled  process $X: \Omega \times [t,\,\infty) \rightarrow E \cup \{\Delta\}$ setting
   \begin{equation}\label{Sec:PDP_controlledX}
   X_s=
   \left\{
   \begin{array}{ll}
   \phi^{\alpha_0}(s-t,  x) \quad &\textup{if}\,\,s \in [t,\, T_{1}),\\
   \phi^{\alpha_n}(s- T_n, E_n)\quad &\textup{if}\,\, s \in [T_n,\, T_{n+ 1}),\,\, n \in \N \setminus \{0\},
   \end{array}	
   \right.
   \end{equation}
 where $\phi^\beta(s, x)$
 is the unique mild solution  to the parabolic  partial differential equation
\begin{align}\label{phialpha}
 \dot x(s) = - L x(s) + b(x(s), \beta(s)),\quad x(0)=x \in E,
\end{align}
 with $\beta$ an $\mathcal A_{ad}^0$-measurable function, namely 
 \begin{align}\label{flow}
 	\phi^\beta(s,x) = S(s)x + \int_0^s S(s-r) b(\phi^\beta(r,x), \beta(r))dr. 
 \end{align}
One can easily prove the  following result, see e.g. Lemma 3.5 in \cite{ReTreThieu}.
 \begin{proposition}\label{P:controlled_flow}
 Let \textbf{\textup{(HL)}} and  \textup{\textbf{(H$\textup{b$\lambda$Q}$)}}  hold. 
Then, for every 
$R>0$, $t \in [0,\,T]$, $t < s'<s$,  
 $\alpha\in \mathcal A_{ad}^t$,
 there exists a constant $C$, only depending on $T$, 
 such that
\begin{align}
&||\phi^{\alpha}(s-t,x)-\phi^{\alpha}(s-t,x')||\leq C \, \omega(||x-x'||),\quad  x,x' \in E,\label{contrflowestimate_x}\\ 
&||\phi^{\alpha}(s-t,x)-\phi^{\alpha}(s'-t,x)||\leq C\,\sigma_R(s-s'),\quad x \in E: ||x||\leq R, \label{contrflowestimate_s}\\ 
&||\phi^{\alpha}(s-t,x)|| \leq C (1 + ||x||), \quad x \in E, \label{contrflowestimate_bound}\\
&||\phi^{\alpha}(s-t,x)-\phi^{\alpha}(s-t,x')||_{-1}\leq C\,  \omega(||x-x'||_{-1})\quad  x,x' \in E,\label{contrflowestimate-1}\\
&||\phi^{\alpha}(s-t,x)-\phi^{\alpha}(s'-t,x)||_{-1}\leq C\,\sigma_R(s-s'),\quad x \in E: ||x||\leq R \label{contrflowestimate-BIS}.
\end{align}
where  $\omega$ is a modulus of continuity, and $\sigma_R$ is a modulus of continuity depending on $R$.  
 \end{proposition}

Set $\mathcal F_0 = \mathcal B([0,\,T])\otimes \mathcal E \otimes \{\emptyset, \Omega'\}$
 and, for all $s \geq t$,  $\mathcal{G}_s^t= \sigma(p((t,r] \times B): r \in (t,s], B \in \mathcal E)$. For all $s\geq t$, let $\mathcal F_s^t$ be the $\sigma$-algebra generated by $\mathcal F_0$ and $\mathcal G_s^t$.
In the following all the concepts of measurability for stochastic processes will refer to the right-continuous, natural  filtration $\mathbb F^t = (\mathcal F_s)_{s \geq t}$.  By the symbol $\mathcal P^t$ we will denote the $\sigma$ algebra of $\mathbb F^t$-predictable subsets of $[t,\,\infty) \times \Omega$.

For every  initial time and starting  point $(t,x)\in  [0,\,T] \times E$ and for each $\alpha \in \mathcal A^t_{ad}$, by Theorem 3.6 in \cite{J} there exists a unique probability measure on $(\Omega, \mathcal F^t_{\infty})$,  denoted by $\P^{t,x}_{\alpha}$, such that its restriction to $\mathcal F^t_t$ is  $\delta_{x}$, and
  the  $\mathbb F^t$-compensator 
 under $\P^{t,x}_{\alpha}$ of the measure $p(ds\,dy)$ is
 $$
 \tilde{p}^{\alpha}(ds\,dy)=
 \sum_{n=1}^{\infty}\one_{[T_n,\,T_{n+1})}(s)\,\lambda(X_{s},\alpha_n(s-T_n,  E_n)
 )\,Q(X_s,\alpha_n(s-T_n,  E_n)
 , dy)\,ds.
 $$ 
 We will denote by $\E^{t,x}_{\alpha}$ the expectation under $\P^{t,x}_{\alpha}$.
   The following  proposition can be obtained by suitably extending the analogous finite-dimensional result, 
  see  Theorem 1.2 in  \cite{ReTreThieu}.
 \endproof
\begin{proposition}\label{T_ItoFormula1}
 Assume that Hypotheses \textup{\textbf{(HL)}} and \textup{\textbf{(H$\textup{b$\lambda$Q}$)}} hold.
For any $(t,x)\in [0,\,T] \times E$ and $\alpha \in \mathcal A_{ad}^t$, let $s \mapsto \phi^{\alpha}(s,x)$ be the unique mild solution to \eqref{phialpha} with $\beta = \alpha$, and $X$ be the process  in \eqref{Sec:PDP_controlledX} with law $\P^{t,x}_{\alpha}$.
Then $X$ is an homogeneous strong Markov process.

Moreover, let   $\mathcal D$ 
be  the set of all measurable functions $\psi:\R_+ \times  E \rightarrow \R$ which are absolutely continuous on $\R_+$ as maps $s \mapsto \psi(s,\phi^{\alpha}(s-t,x))$, for all $x\in E$, and 
such that the map $(x,s,\omega) \mapsto \psi(s,y)- \psi(s,X_{s-})$ is a valid integrand for the random measure $Q$, and set 
\begin{align*}
\bar {\mathcal D}:=\{&\psi \in \mathcal D,\,\psi\in C^{1}(\R\times E):\\
& \,D \psi(s,x)\in E \,\,\textup{if}\,\, x \in E, D \psi(s,x),  \frac{\partial \psi}{\partial s}(s,x) \,\,\textup{bounded if}\,\,x\,\,\textup{bounded}\}.
\end{align*}
 where  $D \psi$ is  the unique element of $E$ such that 
$\frac{d \psi}{d x}[s,x](y) = \langle y,D \psi(s,x)\rangle$, $y \in E$,
where $\frac{d \psi}{d x}[s,x]$ denotes the Fr\'echet-derivative of $\psi$ w.r.t. $x \in E$ evaluated at $(s,x) \in [0,\,T] \times E$.
Let $t < \bar T < T$, $\hat \tau$ be a stopping time such that $\hat \tau  \in [t,\,\bar T]$,     let $\tau_R$ be the exit time of $X$ from $\{y:\,||y||\leq R\}$,  $R>0$, and set $\tau= \hat \tau \wedge \tau_R$.
  Then,  for every $\psi \in \bar {\mathcal D}$, 
 \begin{align}\label{itoformula}
&\E^{t,x}_{\alpha}\left[\psi(\tau,X_\tau)\right] = \psi(t,x)+ \E^{t,x}_{\alpha}\left[\int_t^\tau \left(\frac{\partial \psi}{\partial t}(r,X_r) + \langle b(X_r, \alpha_r),\,D\psi(r,X_r)\rangle\right)dr\right]\\
	&- \E^{t,x}_{\alpha}\left[\int_t^\tau  \langle L\,X_r, \,D\psi(r,X_r)\rangle \,dr\right]+\E^{t,x}_{\alpha}\left[\int_t^\tau\int_E (\psi(r,y)-\psi(r,X_r))\,\lambda(X_r, \alpha_r)\,Q(X_r, \alpha_r,\,dy) \,dr\right]\notag.
\end{align} 
 \end{proposition}

\noindent At this point,  we define for any $(t,x) \in [0,\,T] \times E$ and $\alpha \in \mathcal A^t_{ad}$, the  functional cost
 \begin{equation}\label{Sec:PDP_functional_cost}
 J(t,x,\alpha) =
 \E^{t,x}_{\alpha}\left[\int_{t}^{T}  f(X_{s},\alpha_s)\,ds + g(X_T)\right]
 \end{equation}
 and the value function of the   control problem
 \begin{equation}\label{Sec:PDP_value_function}
 V(t,x) = \inf_{\alpha \in \mathcal{A}^t_{ad}}J(t,x,\alpha),
 \end{equation}

 \begin{proposition}\label{P:dynprogpr}
 	Assume that Hypotheses \textup{\textbf{(HL)}}, \textup{\textbf{(H$\textup{b$\lambda$Q}$)}}
 	and  \textup{\textbf{(H$\textup{fg}$)}}  hold. Then the value function $V$  
 	in \eqref{Sec:PDP_value_function}  is  bounded and uniformly continuous in the $|\cdot| \times ||\cdot||_{-1}$ norm. 
Moreover,  $V$  
  satisfies the following
	dynamic programming principle (DPP):
\begin{align}\label{DynProgPr}
	V(t,x)= \inf_{\alpha \in \mathcal{A}^t_{ad}}\E^{t,x}_{\alpha}\left[\int_{t}^{T_1 \wedge T}  f(X_{s},\alpha_s)\,ds + V(T_1 \wedge T,X_{T_1 \wedge T})\right]\quad  t \in [0,\,T], \,x \in E.
	\end{align}
	 \end{proposition}
\proof
See Section \ref{Sec_proof_contV}.
\endproof 
One can prove that formula \eqref{DynProgPr} also holds with $h \wedge T \wedge T_1$,  for any deterministic time $h >t$,  in place of $ T \wedge T_1$. More generally,  previous result can be extended as  follows.
\begin{proposition}\label{P:DPPNEW} 
Under the same hypotheses of Proposition \ref{P:dynprogpr},  the (DPP) \eqref{DynProgPr} can be  extended to the form 
	\begin{align}\label{DPP2}
	V(t,x)= \inf_{\alpha \in \mathcal{A}^t_{ad}}\E^{t,x}_{\alpha}\left[\int_{t}^{\theta}  f(X_{s},\alpha_s)\,ds + V(\theta,X_{\theta})\right]\quad  t \in [0,\,T], \,x \in E, 
	\end{align}
with 	
$$
\theta := \tau \wedge T_1 \wedge T,  \quad \tau :=\inf \left\{s \geqslant t: 
(s,X_{s}) \notin B\left((t,x); 
\rho\right)
\right\},  
$$
where $B((t,x); \rho):=\{(s,y) \in (t,\,T) \times  E:\,\,||y-x||< \rho,  |s-t|<\rho\}$, $(t,x) \in [0,\,T] \times E$, $\rho >0$.	
	\end{proposition}
\proof
See Section \ref{proofP:DPPNEW}.
\endproof

 \subsection{The related HJB equation}
 Let us now consider the 
HJB-IPDE  
 associated to the optimal control problem: this is the following parabolic nonlinear equation on $[0,\,T]\times E$:
 	\begin{align}\label{Sec:PDP_HJB}
 	&\frac{\partial v}{\partial t} (t,x)
 	-\langle L\,x,\,D v(t,x)\rangle  +\inf_{a \in A}\{\mathcal L^a v(t,x)  + f(x,a) \}=0,\\
 	&v(T, x) = g(x),\label{Sec:PDP_HJB_T}
 	\end{align}
 	where 
 	$\mathcal L^a$ is the time-homogeneous operator depending on $a \in A$
 	defined  as
 	\begin{equation}\label{Sec:PDP_ext_gen_HJB}
 	\mathcal L^a \psi(t,x) := \langle b(x,a), D \psi(t,x)\rangle  +\lambda(x, a)\int_{E} (\psi(t,y)-\psi(t,x))\,   Q(x,a, dy). 
 	\end{equation}

  \begin{remark}\label{Sec:PDP_Rem_form_HJB}
 	The HJB equation 
 	\eqref{Sec:PDP_HJB}-\eqref{Sec:PDP_HJB_T} 
 	can be rewritten as
 \begin{eqnarray}
 H^{v}(x,v,Dv) = 0\label{Sec:PDP_HJB2}\\
 v(T, x) = g(x), \label{Sec:PDP_HJB_T2}
 \end{eqnarray}
 where
 \begin{align*}
 H^\psi(z, v, p) &= \frac{\partial v}{\partial t}
-\langle L\,z,\,p\rangle  + \inf_{a \in A}\left\{b(z,a) \cdot  p +\int_{E} (\psi(y)-\psi(z))\, \lambda(z, a)\,  Q(z,a, dy) +f(z,a) \right\}.
 \end{align*}
 \end{remark}

   \begin{definition}\label{D:testfunc}
 	We  say that a function $\psi$ is a test function if 
 	$
 	\psi(t,x) = \varphi(t,x) + \delta(t,x)\, h(||x||)
 	$,
 	where 
 	\begin{itemize}
 		\item [(i)]  $\psi, \frac{\partial \varphi}{\partial t}$, $D \varphi$, 
 		$L^\ast D\varphi$, $\frac{\partial \delta}{\partial t}$, $D \delta$,
 		 $L^\ast D\delta$ are uniformly continuous on $(\varepsilon,\,T-\varepsilon) \times E$ for every $\varepsilon >0$, $\delta \geq 0$  is $B$-continuous  and bounded, $\varphi$ is $B$-lower semicontinuous and bounded.
 		\item [(ii)] $h \in C^2(\R)$ with $h', h''$ uniformly continuous, $h$ is even and bounded, 
 		$h'(r) \geq 0$ for $r \in (0,\,+ \infty)$.
 	\end{itemize}
 \end{definition}

\begin{definition}\label{Sec:PDP_Def_viscosity_sol_HJB} Viscosity solution to \eqref{Sec:PDP_HJB}-\eqref{Sec:PDP_HJB_T}.
 	\begin{itemize}
 		\item[(i)] A bounded $B$-upper-semicontinuous  function $u:(0,\,T)\times E\rightarrow \R$  is  a \emph{viscosity subsolution} of \eqref{Sec:PDP_HJB}
 		if, 
 		whenever $u- \psi$ has a global maximum at a point $(t,x)$ for a test function $\psi$, then 
 		\begin{align*}
 		&
 		\frac{\partial \psi}{\partial t} (t,x)
 		-\langle x,\,L^\ast \, D \varphi(t,x) + h(||x||)\,L^\ast D\delta(t,x)\rangle \\
 		& + \inf_{a \in A}\left\{\langle b(x,a),  D \psi(t,x) \rangle  +\int_{E} (\psi(t,y)-\psi(t,x))\, \lambda(x, a)\,  Q(x,a, dy) +f(x,a) \right\}  \geq  \,\,0.
 		\end{align*}
 	 	\item[(ii)] A bounded $B$-lower-semicontinuous  function $w:(0,\,T) \times E\rightarrow \R$  is  a \emph{viscosity supersolution} of \eqref{Sec:PDP_HJB}
 	 	if,
 		 whenever $w+ \psi$ has a global minimum at a point $(t,x)$ for a test function $\psi$, then 
 		\begin{align*}
 		& 
 		-\frac{\partial \psi}{\partial t}(t,x)
 		 +\langle x,\,L^\ast \, D \varphi(t,x) + h(||x||)\,L^\ast D\delta(t,x)\rangle \\
 		& + \inf_{a \in A}\left\{\langle b(x,a),  -D \psi(t,x) \rangle  -\int_{E} (\psi(t,y)-\psi(t,x))\, \lambda(x, a)\,  Q(x,a, dy) +f(x,a) \right\}  \leq  \,\,0.
 		\end{align*}
 		\item[(iii)] A \emph{viscosity solution} of \eqref{Sec:PDP_HJB}-\eqref{Sec:PDP_HJB_T} is a function which is both a viscosity subsolution and a viscosity supersolution.
 	\end{itemize}
 \end{definition}

The following lemma will play a fundamental role in the following.

\begin{lemma}\label{Gpsi} 

Let $\psi(s,y) = \varphi(s,y) + \delta(s,y) \,h(||x||)$ be a test function of the type introduced in Definition \ref{D:testfunc}. For any $a\in A$, define
\begin{equation}\label{defGapsi}
 G_a^\psi\left(s,z\right) :=-\frac{\partial \psi}{\partial s} (s,z)+\langle z,\,L^\ast \, D \varphi(s,z) + h(||z||)\,L^\ast D\delta(s,z)\rangle+f(z,a) - \mathcal L^a \psi(s,z)
 \end{equation}
 where $\mathcal L^a$ is defined in (\ref{Sec:PDP_ext_gen_HJB}).
Then, for any $t \in (\varepsilon,\,T-\varepsilon)$, $\varepsilon >0$, $x\in E$,  and  any measurable function $\alpha_0 : \R_+ \times E \rightarrow A$, the map 
\begin{align*}
r  \mapsto  G_a^\psi\left(r,\phi^{\alpha_0}(r-t,x)\right)  
 \end{align*}
is continuous on $[t,\,T-\varepsilon)$, $\varepsilon >0$, uniformly in $a$ and $\alpha_0$.
In particular, for any $t \in (\varepsilon,\,T-\varepsilon)$, $\varepsilon >0$, $x\in E$,  and  any measurable function $\alpha_0 : \R_+ \times E \rightarrow A$, the map 
\begin{align*}
r  \mapsto \mathcal G^{\alpha_0}(r) :=
\inf_{a \in A} G_a^\psi\left(r,\phi^{\alpha_0}(r-t,x)\right)  
 \end{align*}
is continuous on $[t,\,T-\varepsilon)$, $\varepsilon >0$, uniformly in $\alpha_0$, and uniformly on  $B_R(x):=\{x \in E: ||x||\leq R\}$, $R >0$.
\end{lemma}
\proof
See Section \ref{Proof Lemma Gpsi}.
\endproof

\begin{theorem}\label{Sec:PDP_Thm_existence}
	Let \textbf{\textup{(HL)}}, \textup{\textbf{(H$\textup{b$\lambda$Q}$)}} and  \textup{\textbf{(H$\textup{fg}$)}} hold. Then the value function  $V$ 
 	 provides
 	a
 	viscosity solution to \eqref{Sec:PDP_HJB}-\eqref{Sec:PDP_HJB_T}.
\end{theorem}
\proof
See Section \ref{Sec_visc_prop_V}.
\endproof
 
 \section{Control randomization}\label{Sec:PDP_Section_dual_control}
 In this section we start to  implement the control randomization method.
As a first step, for an initial time $t \geq 0$ and a starting point $x \in E$, we construct an (uncontrolled) PDMP $(X,I)$ with values in $E \times A$ by specifying its local characteristics, see \eqref{Sec:PDP_h_XI}-\eqref{Sec:PDP_lambda_XI}-\eqref{Sec:PDP_Q_XI} below. Next we formulate an auxiliary optimal control problem where,
 roughly speaking, we optimize a  functional cost by modifying the intensity of the process $I$ over a suitable family of probability measures.

 \subsection{Construction of randomized state systems 
 }\label{Sec:PDP_Section_control_rand}
 Let $E$ still denote a real separable Hilbert space  Borel $\sigma$-field $\mathcal E$, and  $A$ be a Polish space with corresponding Borel $\sigma$-field $\mathcal A$. Let moreover $b$, $\lambda$ and  $Q$ be respectively two real  functions on $E \times A$ and  a   probability transition  from $(E \times A, \mathcal E \otimes \mathcal A)$, satisfying \textbf{(H$\textup{b$\lambda$Q}$)}
 as before.
 We denote by $\phi(s, x,a)$
 the unique mild  solution to the parabolic partial differential equation
 \begin{equation}\label{PDE}
 \dot x(s) = - L\,x(s) +  b(x(s), a), \quad x(0)=x \in E, \,\,a \in A.
\end{equation}
 In particular, $\phi(s,x,a)$ corresponds to the function $\phi^\beta(s,x)$ introduced in Section \ref{Sec:PDP_Sec_control_problem} when $\beta(s)\equiv a$ and, for every  
$x,x' \in E$, $0 < s'<s<T$,  $a\in A$, satisfies 
\begin{align}\label{flowestimateRand}
||\phi(s,x,a)-\phi(s',x',a)||_{-1}\leq C \omega(||x-x'||_{-1} + (s-s'))
\end{align}
with $C$ a constant only depending on $T$, and $\omega$  some modulus of continuity
 by Proposition \ref{P:controlled_flow}. This fact will be of great use in the sequel.

 Let us now introduce another finite measure $\lambda_0$ on $(A, \mathcal A)$ satisfying the following assumption:
 
 \medskip
 
 \noindent \noindent \textbf{(H\textup{$\lambda_0$})}\quad $\lambda_0$ is a finite measure on $(A, \mathcal A)$ with full topological support.
 
 \medskip
 
 \noindent The existence of such a measure is guaranteed by the fact that  $A$ is a separable space  with   metrizable topology. We define
 \begin{eqnarray}
 \tilde \phi(t,x,a)&:=& (\phi(t,x,a), \quad a),
 \label{Sec:PDP_h_XI}\\ 
 \tilde{\lambda}(x,a) &:=& \lambda(x,a)+ \lambda_0(A),\label{Sec:PDP_lambda_XI}\\  \tilde{Q}(x,a,dy\,db) &:=& \frac{\lambda(x,a)\,Q(x,a,dy)\,\delta_{a}(db) + \lambda_0(db)\,\delta_{x}(dy)}{\tilde{\lambda}(x,a)}.\label{Sec:PDP_Q_XI}  \end{eqnarray}
 We wish to construct a PDMP $(X,I)$ 
 with enlarged state space $E \times A$ and  local characteristics $(\tilde \phi,\tilde \lambda, \tilde Q)$.
 Firstly, we need to introduce  a suitable sample space to describe the jump mechanism of the  process $(X,I)$ on $E \times A$.
 Accordingly, we  set  $\Omega'$ as the set of sequences $\omega'=(t_n, e_n, a_n)_{n \geq 1}$ contained in $((t,\,\infty) \times E \times A) \cup \{ (\infty, \Delta, \Delta')\}$, where $\Delta \notin E$ (resp. $\Delta' \notin A$) is adjoined to $E$ (resp. to $A$) as an isolated point. In the sample space $\Omega = [0,T] \times E \times A \times \Omega'$ we define the random variables $T_n : \Omega \rightarrow (t,\,\infty]$, $E_n : \Omega \rightarrow E \cup \{\Delta\}$, $A_n : \Omega \rightarrow A \cup \{\Delta'\}$, as follows: writing $\omega =(t,x,a,\omega')$ in the form $\omega = (t,x,a,t_1, e_1, a_1, t_2, e_2, a_2,...)$, we set for $n \geq 1$,
 \begin{align*}
 & T_n (\omega )= t_n, \qquad T_\infty  (\omega )= \lim_{n\to\infty} t_n,
 \qquad T_0(\omega)=t,\\
 & E_n(\omega)=e_n,\qquad A_n(\omega)=a_n,
 \qquad E_t(\omega)=x, \qquad A_t(\omega)=a.
 \end{align*}
 We  define the process $(X,I)$ on $(E \times A) \cup \{\Delta, \Delta'\}$ setting
 \begin{align}\label{Sec:PDP_XI_def}
 (X,I)_s&=
 \left\{
 \begin{array}{lll}
 (\phi(s-t, x,a),a)& \quad \textup{if}\,\, s < T_1,\\
(\phi(s-T_n, E_n,A_n),A_n)& \quad \textup{if}\,\, T_n \leq s < T_{n+1},\,\,\textup{for}\,\,n \in \N,\\
 (\Delta,\Delta')&  \quad \textup{if}\,\,s \geq T_{\infty}.
 \end{array}	
 \right.
 \end{align}
 In $\Omega$ we introduce for all $s\geq t$ the $\sigma$-algebras
 $\mathcal G^t_r=\sigma(N(s,G)\,:\, s\in (t,r], G\in\mathcal E \otimes \mathcal A)$ generated
 by the counting processes  $N(s,G)=
 \sum_{n \in \N}\one_{T_n\le s}\one_{(E_n,A_n)\in G}$, and the
 $\sigma$-algebra $\mathcal F^t_s$ generated by $\mathcal F_0$ and $\mathcal G^t_s$, where
 $\mathcal F_0=\mathcal B([0,\,T]) \otimes \mathcal E\otimes \mathcal A \otimes \{\emptyset,\Omega'\}$.
 We still denote by $\F^t = (\mathcal F^t_s)_{s\geq t}$ and $\mathcal P^t$ the corresponding filtration and predictable $\sigma$-algebra.
 The random measure $p$ is now defined on $(t,\,\infty) \times E \times A$ as
 \begin{equation}\label{Sec:PDP_p_dual}
 p(ds\,dy\,db)= \sum_{n\in \N} \delta_{(T_n,E_n, A_n)}(ds\,dy\,db).
 \end{equation}
 Given any starting point $(t,x,a) \in E \times A$,
 by  Theorem 3.6 in \cite{J}, there exists a unique
 probability measure    on $(\Omega, \mathcal F^t_{\infty})$, denoted by $\P^{t,x,a}$, such that its restriction to  $\mathcal F_0$ is $\delta_{(x,a)}$  and the $\mathbb F^t$-compensator  of the measure $p(ds\,dy\,db)$ under $\P^{t,x,a}$ is the random measure
 \begin{equation}\label{tildep}
 \tilde p(ds\,dy\,db)= \sum_{n \in \N}\one_{[T_n,\,T_{n+1})}(s)\,\Lambda( \phi(s-T_n, E_n,A_n),A_n, dy\,db)\,ds,
 \end{equation}
 where
 $$
 \Lambda(x,a,dy\,db)= \lambda(x,a)\,Q(x,a,dy)\,\delta_a(db)+ \lambda_0(db)\,\delta_x(dy), \quad \forall (x,a)\in E \times A.
 $$
 We denote by  $q=p-\tilde{p}$ the compensated martingale measure associated to $p$.

 The sample path of a process $(X,I)$ with values in $E \times A$, starting from a fixed initial point $(x,a) \in  E \times A$ at time $t$, can be defined iteratively by means of its local characteristics $(\tilde \phi, \tilde \lambda, \tilde Q)$ in the following way. Set
 $F(t,x,a;s)=e^{-\int_t^s (\lambda (\phi(r-t,x,a),a) + \lambda_0(A))\,dr}$.
 For any $B \in \mathcal E,\,C \in \mathcal A$, we have
 \begin{align}
 &\P^{t,x,a}(T_1 > s)=F(t,x,a;s),\quad s \geq t,\label{Sec:PDP_Abis}\\
 &\P^{t,x,a}(X_{T_1} \in B, I_{T_1}\in C|\,T_1)= \tilde Q(x, B \times C),\label{Sec:PDP_Bbis}
 \end{align}
 on $\{T_1 < \infty\}$, and, for every $n \geq 1$,  on $\{T_n < \infty\}$,
 \begin{align}
 &\P^{t,x,a}(T_{n+1} > s\,|\,\mathcal F_{T_n})=\exp\left(-\int_{T_n}^s (\lambda (\phi(r-T_n,X_{T_n}), I_{T_n})+ \lambda_0(A))\,dr\right),\quad s \geq T_n,\label{Sec:PDP_A_kbis}\\
 &\P^{t,x,a}(X_{T_{n+1}} \in B, I_{T_{n+1}} \in C 
 |\,\mathcal F_{T_n},\,T_{n+1})= \tilde Q(\phi(T_{n+1}-T_n,X_{T_n},I_{T_n}), I_{T_n}, B \times C). 
  \label{Sec:PDP_B_kbis}
 \end{align}
 We recall the following  result, that is a direct consequence of Theorem 4 in \cite{Riedler2011}.
\begin{proposition}\label{T_ItoFormula}
For any $(t,x,a)\in [0,\,T] \times E \times A$, let $\phi(t,x,a)$ be the unique mild solution to \eqref{PDE}, and $(X,I)$ be the process defined in \eqref{Sec:PDP_XI_def} with law $\P^{t,x,a}$. Then $(X,I)$ is an homogeneous strong Markov process.

Moreover,  denote by	 $\mathcal D$ 
 the set of all measurable functions $\varphi: E \times A\rightarrow \R$ which are absolutely continuous on $\R_+$ as maps $s \mapsto \varphi(\phi(s,x,a),a)$, for all $(x,a)\in E \times A$, 
and 	such that the map $(x,a,s,\omega) \mapsto \varphi(y,b)- \varphi(X_{s-}, I_{s-})$ is a valid integrand for the random measure \eqref{tildep}, and 
 set  
 $$
\bar {\mathcal D}:=\{\varphi \in D(\mathcal L),\,\varphi\in C^{1,0}(E \times A),\,D \varphi(x,a)\in E \,\,\textup{if}\,\, x \in E, D \varphi(x,a) \,\,\textup{bounded if}\,\,x\,\,\textup{bounded}\},
 $$ 
where  $D \varphi$ is  the unique element of $E$ such that 
$\frac{d \varphi}{d x}[x,a](y) = \langle y,D \varphi(x,a)\rangle$, $y \in E$,
where $\frac{d \varphi}{d x}[x,a]$ denotes the Fr\'echet-derivative of $\varphi$ w.r.t. $x \in E$ evaluated at $(x,a) \in E \times A$.
Then 
the extended generator of $(X,I)$ is given by
 \begin{align*}
 \mathcal{L}\varphi(x,a) :=& \langle -L x +  b(x,a),\,  D \varphi(x,a)\rangle + \int_{E}(\varphi(y,a) - \varphi(x,a))\, \lambda(x,a)\,Q(x,a,dy)\\
 &+\int_{A}(\varphi(x,b) - \varphi(x,a))\, \lambda_0(db), \quad  \textup{ for every } \varphi \in \bar {\mathcal D}.
 \end{align*}
 \end{proposition}

 \subsection{The randomized optimal control problem}\label{Sec:PDP_Section_dual_optimal_control}
 We now introduce a randomized   
 optimal control problem associated to the process $(X,I)$, and formulated in a weak form. For fixed $(t,x,a)$, we consider a family of probability measures $\{\P_{\nu}^{t,x,a},\,\nu \in \mathcal V\}$ in the space $(\Omega, \mathcal F_{\infty})$, 
 whose effect   is to change the stochastic intensity of the process $(X,I)$. 
 
 Let us proceed with precise definitions.
 We still  assume that {\bf (Hb$\lambda$Q)},
  {\bf (H$\lambda_0$)} and \textbf{(H$\textup{fg}$)} hold.
 We recall that $\mathbb{F}^t= (\mathcal{F}^t_s)_{s \geqslant t}$ is the augmentation of the natural filtration generated by $p$ in \eqref{Sec:PDP_p_dual}, and that $\mathcal{P}^t$
 denotes  the $\sigma$-field of $\mathbb{F}^t$-predictable subsets of $[t,\,\infty) \times \Omega$.
 We define
 $$
 \mathcal{V} = \{ \nu: \Omega \times [0,\,\infty) \times A \rightarrow (0,\,\infty)\,\, \mathcal{P}^0\otimes \mathcal A\text{-measurable and bounded}\}.
 $$
 For every $\nu \in \mathcal{V}$,
 we consider the predictable random measure
 \begin{align}\label{Sec:PDP_dual_comp}
 \tilde{p}^\nu(ds\,dy\,db)&:= \nu_s(b) \, \lambda_0 (db) \, \delta_{\{ X_{s-} \}}(dy) \, ds +\,
 \lambda(X_{s-},\,I_{s-})\,Q(X_{s-},\,I_{s-},\,dy)\, \delta_{\{ I_{s-} \}}(db)\, ds.
 \end{align}
 In particular, for any $t \in [0,\,T]$, by the Radon Nikodym theorem one can find two nonnegative functions $d_1$, $d_2$ defined on $\Omega \times [t,\,\infty) \times E \times A$, $\mathcal P \otimes \mathcal E \otimes \mathcal A$, such that
 \begin{eqnarray*}
 	\lambda_0(db) \, \delta_{\{ X_{s-} \}}(dy) \, ds &=& d_1(s,y,b)\, \tilde{p}(ds\,dy\,db)
 	\\
 	\lambda(X_{s-},\,I_{s-},\,dy) \, \delta_{\{ I_{s-} \}}(db) \, ds &=&
 	d_2(s,y,b)\, \tilde{p}(ds\,dy\,db),
 	\\
 	d_1(s,y,b) + d_2(s,y,b) &=& 1, \qquad \tilde{p}(ds\,dy\,db)\textup{-a.e}.
 \end{eqnarray*}
 and we have
 $
 d \tilde{p}^{\nu} = (\nu \, d_1  + d_2)\, d\tilde{p}$.
 For any  $t \in [0,\,T]$, $\nu \in \mathcal{V}$, consider then the Dol\'eans-Dade exponential local martingale
 $L^{t,\nu}$ defined 
 \begin{align}
 L_s^{t,\nu} &= \exp\bigg(\int_t^s\!\int_{E\times A}\log(\nu_r(b) \, d_1(r,y,b) \, + d_2(r,y,b))\,p(dr \, dy \,db)- \int_{t}^{s}\!\int_{A}(\nu_r(b) - 1)\lambda_0(db)\,dr\bigg)
 \nonumber\\
 &= e^{\int_{t}^{s}\int_{A}(1 - \nu_r(b) )\lambda_0(db)\,dr} \prod_{n \geqslant 1: t \leq  T_{n} \leq s}(\nu_{T_{n}}(A_n)\,d_1(T_{n},E_n,A_n) + d_2(T_{n},E_n,A_n)),\label{Sec:PDP_Lnu}
 \end{align}
 for $s\geq t$. 
 When $(L^{t,\nu}_s)_{s \geq t}$ is a true martingale on $[t,\,T]$,
 we can define  a probability measure $\P^{t,x,a}_{\nu}$ equivalent to $\P^{t,x,a}$ on $(\Omega,\,\mathcal{F}^t_T)$ by
 \begin{equation}\label{Sec:PDP_PnuT}
 \P^{t,x,a}_{\nu}(d\omega)=L_T^\nu(\omega)
 \,\P^{t,x,a}(d\omega).
 \end{equation}
 By the Girsanov theorem for point processes (see Theorem 4.5 in \cite{J}),
 the restriction of the random measure $p$ to $(t,T]\times E\times A$
 admits $\tilde{p}^\nu=(\nu\,d_1 + d_2)\,\tilde{p}$
 as compensator   under $\P^{t,x,a}_{\nu}$.
 We set $q^\nu: = p-\tilde{p}^\nu$,  and we denote by $\mathbb{E}_{\nu}^{t,x,a}$ the expectation operator under $\P_{\nu}^{t,x,a}$.
 Previous considerations are formalized in  the following lemma, for a proof see   
 Lemma 3.2 in \cite{BandiniFuhrman}.
 \begin{lemma}\label{Sec:PDP_lemma_P_nu_martingale}
 	Let Hypotheses {\bf (Hb$\lambda$Q)} 
	and  {\bf (H$\lambda_0$)} hold. Then, for every  $(t,x,a)\in [0,\,T] \times E\times A$  and
 	$\nu \in \mathcal{V}$, under the probability
 	$\P^{t,x,a}$, 
 	the process $(L^{t,\nu}_s)_{s \geq t}$ is a   martingale.
 	Moreover,   $L_T^{t,\nu}$ is square integrable,		
 	and, for every  $\mathcal{P}_T \otimes \mathcal E\otimes \mathcal A$-measurable  function $H: \Omega \times [0,T] \times E \times A \rightarrow \R$ such that  $\spertxa{\int_t^T \int_{E\times A} |H_s(y,b)|^2\, \tilde{p}(ds\,dy\,db)}$ $< \infty$,  the process $\int_{t}^{\cdot}\int_{E \times A}H_r(y,b)\, q^\nu(dr\,dy\, db)$ is a $\P^{t,x,a}_{\nu}$-martingale on $[t,T]$.
 \end{lemma}
 Finally, for every $(t,x)\in [0,\,T] \times E$, $a \in A$ and $\nu \in \mathcal V$, we introduce the dual  functional cost
 \begin{equation}\label{Sec:PDP_dual_functional_cost}
 J(t,x,a,\nu) := \spernutxa{\int_{t}^{T} f(X_s,I_s)\,ds + g(X_T)},
 \end{equation}
 and the dual value function
 \begin{equation}\label{Sec:PDP_dual_value_function}
 V^{\ast}(t,x,a) := \inf_{ \nu \in \mathcal{V}} J(t, x,a,\nu).
 \end{equation}

 \section{A constrained BSDEs representation for the value function}
  \label{Sec:PDP_Sec_ConstrainedBSDE}
 In the present section we introduce a BSDE with a sign constraint on its martingale part for which we give  existence and uniqueness of a maximal solution in an appropriate sense.
 This constrained BSDE will provide a probabilistic representation formula for the dual value function introduced in \eqref{Sec:PDP_dual_value_function}.
 
 Throughout the section we still  assume that  \textbf{(H$\textup{b$\lambda$Q}$)},
 {\bf (H$\lambda_0$)}
 and {\bf (Hfg)} hold.
 The random measures $p$, $\tilde p$ and $q$, as well as  the dual control setting $\Omega, \mathbb F^t, (X,I), \P^{t,x,a}$, are the same as in Section \ref{Sec:PDP_Section_control_rand}.
 For any $(t,x,a) \in [0,\,T] \times E \times A$ we introduce the following notation.
 \begin{itemize}
 	\item $\textbf{L}^\textbf{2}_{\textbf{t,x,a}}(\mathcal{F}^t_\tau)$,  the set of $\mathcal{F}^t_\tau$-measurable random variables $\xi$ such that $\spertxa{|\xi|^2} < \infty$; here  $\tau \geqslant 0$ is an $\mathbb F^t$-stopping time.
 	\item $\textbf{S}^\infty$,  the set of real-valued c\`adl\`ag 	adapted processes $Y = (Y_t)_{t \geqslant 0}$ which are uniformly bounded.
 	\item $\textbf{L}_\textbf{t,x,a}^{\textbf{2}}(\textup{q})$,  the set of $\mathcal{P}_T\otimes \mathcal{B}(E) \otimes \mathcal A$-measurable maps $Z: \Omega \times [0,\,T] \times E \times A \rightarrow \R$ such that
 	\begin{align*}
 	||Z||^2_{\textbf{L}_{\textbf{t,x,a}}^{\textbf{2}}(\textup{q})}:
 	&  = \mathbb{E}^{t,x,a}\Big[  \int_{t}^{T}\int_{E} |Z_s(y,I_s)|^2 \, \lambda(X_s,I_s)\,Q(X_s,I_s,dy)\,ds\Big]\\ 
 	&+ \mathbb{E}^{t,x,a}\Big[\int_{t}^{T}\int_{A} |Z_s(X_s,b)|^2 \, \lambda_0(db)\,ds \Big] < \infty.
 	\end{align*}
 	\item $\textbf{K}^\textbf{2}_\textbf{t,x,a}$,  the set of nondecreasing c\`adl\`ag predictable  processes $K = (K_s)_{t \leqslant s \leqslant T}$ such that $K_t = 0$ and $\spertxa{|K_T|^2}< \infty$. 
 \end{itemize}
We consider the following family of BSDEs with partially nonnegative jumps over a finite horizon $T$, parametrized by $(t,x,a)$: 
 ${\mathbb{P}}^{t,x,a}$-a.s.,
 \begin{align}
 Y^{t,x,a}_{s} &= g(X_T) + \int_{s}^{T}f(X_r,I_r)\,dr -( K^{t,x,a}_T - K^{t,x,a}_s) \label{Sec:PDP_BSDE}\\
 &\hspace{-10mm}- \int_{s}^{T}\int_{A}Z^{t,x,a}_r(X_r,\,b)\, \lambda_0(db)\,dr - \int_{s}^{T}\int_{E \times A}Z^{t,x,a}_r(y,\,b)\, q(dr\,dy\, db), \quad  t\leqslant s \leqslant T,\nonumber
 \end{align}
 with
 \begin{equation}\label{Sec:PDP_BSDE_constraint}
 Z_s^{t,x,a}(X_{s-},b)\geqslant 0, \qquad  ds \otimes d\P^{t,x,a} \otimes \lambda_0(db)\text{ -a.e. on } [0,\,T]\times \Omega \times A.
 \end{equation}
We are interested in the \emph{maximal solution} $(Y^{t,x,a},Z^{t,x,a},K^{t,x,a})\in \textbf{S}^{\infty}\times \textbf{L}^{\textbf{2}}_\textbf{t,x,a}(\textup{q})\times \textbf{K}^\textbf{2}_\textbf{t,x,a}$ to \eqref{Sec:PDP_BSDE}-\eqref{Sec:PDP_BSDE_constraint},
 in the sense that for any other solution $(\tilde{Y}, \tilde{Z},\tilde{K})\in \textbf{S}^{\infty}\times \textbf{L}^{\textbf{2}}_\textbf{t,x,a}(\textup{q})\times \textbf{K}^\textbf{2}_\textbf{t,x,a}$ to \eqref{Sec:PDP_BSDE}-\eqref{Sec:PDP_BSDE_constraint}, we  have $Y_s^{t,x,a} \geqslant \tilde{Y}_s$, $\P^{t,x,a}$-a.s., for all $s\geqslant t$.

Let us  introduce the following penalized BSDE, associated to \eqref{Sec:PDP_BSDE}-\eqref{Sec:PDP_BSDE_constraint}, parametrized by the integer $n \geq 1$:
\begin{align}
 Y^{n,t,x,a}_{s} &= g(X_T) + \int_{s}^{T}f(X_r,I_r)\,dr -( K^{n,t,x,a}_T - K^{n,t,x,a}_s) \label{Sec:penalized_BSDE}\\
 &\hspace{-10mm}- \int_{s}^{T}\int_{A}Z^{n,t,x,a}_r(X_r,\,b)\, \lambda_0(db)\,dr - \int_{s}^{T}\int_{E \times A}Z^{n,t,x,a}_r(y,\,b)\, q(dr\,dy\, db), \quad  t\leqslant s \leqslant T,\nonumber
 \end{align}
 where $K_s^{n,t,x,a}:=n \int_0^s \int_A [Z^{n,t,x,a}_r(X_r,\,b)]^{-}\lambda_0(db) \,dr$, $s \in [t,\,T]$.
  \begin{theorem}\label{Sec:PDP_Thm_ex_uniq_maximal_BSDE}
 	Let Hypotheses  \textup{\textbf{(H$\textup{b$\lambda$Q}$)}}, \textup{\textbf{(H$\lambda_0$)}} and  \textup{\textbf{(H$\textup{fg}$)}} hold. 
 	Then, for every  $(t,x,a) \in [0,\,T] \times E \times A$, there exists a unique maximal solution $(Y^{t,x,a},Z^{t,x,a},K^{t,x,a})\in \textup{\textbf{S}}^{\infty}\times \textup{\textbf{L}}^{\textbf{2}}_\textbf{t,x,a}(\textup{q})\times \textup{\textbf{K}}^\textbf{2}_\textbf{t,x,a}$  to the BSDE with partially nonnegative jumps \eqref{Sec:PDP_BSDE}-\eqref{Sec:PDP_BSDE_constraint}, where
	$Y^{t,x,a}$ is the nonincreasing limit of $(Y^{n,t,x,a})_n$, $Z^{t,x,a}$ is the weak limit of $(Z^{n,t,x,a})_n$ in $\textup{\textbf{L}}^{\textbf{2}}_\textbf{t,x,a,\textup{loc}}(\textup{q})$ and  $K_{s}^{t,x,a}$ is the weak limit of $(K_{s}^{n,t,x,a})_n$ in $\textup{\textbf{L}}^{\textbf{2}}_{\textbf{t,x,a}}(\mathcal{F}_{s})$, for any $s \geqslant 0$.
 	Moreover, $Y^{t,x,a}$  has the explicit representation:
 	\begin{equation}\label{Sec:PDP_rep_Y}
 	Y_s^{t,x,a} = \essinf_{\nu \in \mathcal{V}}\spernutxa{\int_s^T f(X_{r},I_{r})\, dr + g(X_T)\Big| \mathcal{F}_s}, \,\, \forall \,\,s \in [t,T].
 	\end{equation}
In particular, setting $s=t$ in \eqref{Sec:PDP_rep_Y}, 
we have the following representation formula for the value function of the randomized control problem:
 	\begin{equation}\label{Sec:PDP_Vstar_Y0}
 	V^{\ast}(t,x,a)
 	= Y_t^{t,x,a}, \quad (t,x,a)\in [0,\,T] \times E \times A.
 	\end{equation}
 \end{theorem}
 \proof
The proof of this result
 is  analogous  
to the one  for the BSDE \eqref{Sec:PDP_BSDE} with underlying   finite-dimensional process $X$,   see Theorem 4.7 in \cite{BandiniPDMPsNoBordo},  and we do not report it for sake of brevity.
 \endproof

 Our main purpose is to show how maximal solutions to BSDEs with nonnegative jumps of the form \eqref{Sec:PDP_BSDE}-\eqref{Sec:PDP_BSDE_constraint} provide actually a Feynman-Kac representation to the value function $V$ associated to our optimal control problem for infinite-dimensional  PDMPs.
 Let us introduce a deterministic function $v: [0\,T] \times E \times A \rightarrow \R$ as
 \begin{equation}
 v(t,x,a):= Y_t^{t,x,a}, \quad (t,x,a) \in [0,\,T] \times E\times A. \label{Sec:PDP_def_v}
 \end{equation}
  \begin{proposition}\label{Sec:Prop_Feynman_Kac_HJB}
 	Assume that Hypotheses \textup{\textbf{(HL)}}, \textup{\textbf{(H$\textup{b$\lambda$Q}$)}}, \textup{\textbf{(H$\lambda_0$)}}, 
 	and  \textup{\textbf{(H$\textup{fg}$)}}  hold. Then the function $v$  
 	in \eqref{Sec:PDP_def_v} does not depend on the variable $a$:
 	\begin{equation}\label{Sec:PDP_vdelta_not_dep_a}
 	v(t,x,a)= v(t,x,a'),\quad t \in [0,\,T],\,x \in E,\,a,a' \in A.
 	\end{equation}
By abuse of notation, we define the function $v$ on $[0,\,T] \times E$ by
 	$v(\cdot,\cdot)= v(\cdot,\cdot,a)$, for any $a \in A$.
Moreover, $v$ admits the representation formula: $\P^{t,x,a}$-a.s.
\begin{equation}\label{Sec:PDP_ident_vdelta}
 	v(s, X_{s})= Y_s^{t,x,a},\quad s \geqslant t.
 	\end{equation}
 	 \end{proposition}
 	 \proof
By  Lemma 5.3 and Remark 5.5 in  \cite{BandiniPDMPsNoBordo}, we have that 
for any $(t,x,a)\in [0,\,T] \times E \times A$, $\P^{t,x,a}$-a.s., 
 	\begin{equation}\label{Sec:PDP_ident_vdelta}
 	v(s, X_{s}, I_s)= Y_s^{t,x,a},\quad s \geqslant 0.
 	\end{equation}
 Now we recall that,
 by \eqref{Sec:PDP_Vstar_Y0} and \eqref{Sec:PDP_def_v}, $v$ coincides with the  value function $V^{\ast}$ of the dual control problem introduced in Section \ref{Sec:PDP_Section_dual_optimal_control}.
 Therefore, identity \eqref{Sec:PDP_vdelta_not_dep_a} corresponds to the fact
that $V^{\ast}(t,x,a)$ does not depend on $a$. 
Proceeding as  in the finite-dimensional case (see the proof of  Proposition 5.6 in \cite{BandiniPDMPsNoBordo}), one can prove that:
\begin{align}
&\textup{for any}\,\, t \in [0,\,T],\,\, x\in E,\,\, a, a' \in A,  \,\,
 	\nu\in \mathcal V, \,\,
\textup{there exists}\,\, (\nu^{\varepsilon})_{\varepsilon}\in \mathcal V: \notag\\
&\qquad \qquad \qquad  \lim_{\varepsilon \rightarrow 0^+} J(t,x,a',\nu^{\varepsilon}) = J(t,x,a,\nu). \label{Sec:PDP_claim_PDMP}
\end{align}
Property  \eqref{Sec:PDP_claim_PDMP}
implies that $V^{\ast}(t,x,a') \leq J(t,x,a,\nu)$ 
for all $t \in [0,\,T]$, $x \in E$, $a,a' \in A$,
 and by the arbitrariness of $\nu$ one can conclude that
 $V^{\ast}(t,x,a') \leq V^{\ast}(t,x,a)$ for all  $t \in [0,\,T]$, $x \in E$, $a,a' \in A$.
 In other words $V^\ast(t,x,a)=v(t,x,a)$ does not depend on $a$, and \eqref{Sec:PDP_vdelta_not_dep_a} holds.  
 	 \endproof
 \begin{theorem}\label{Sec:PDP_THm_Feynman_Kac_HJB}
 	Assume that Hypotheses \textup{\textbf{(HL)}}, \textup{\textbf{(H$\textup{b$\lambda$Q}$)}}, \textup{\textbf{(H$\lambda_0$)}}, 
 	and  \textup{\textbf{(H$\textup{fg}$)}}  hold. 
Then $v$ is  bounded and uniformly continuous in the $|\cdot| \times ||\cdot||_{-1}$ norm. 
 Moreover, $v$ satisfies the so called randomized dynamic programming principle:
\begin{equation}\label{RandDynProgPr}
 	v(t,x) = \inf_{\nu \in \mathcal{V}} \spernutxa{\int_t^{T \wedge T_1} f(X_{r},I_{r})\, dr + v(T \wedge T_1, X_{T \wedge T_1})}.
	 	\end{equation}
 	 \end{theorem}
\proof
See Section \ref{Sec:proofThm5.1}.
\endproof	 
We can now give the  following important result.
 \begin{theorem}\label{Sec:first_main result}
 	Assume that Hypotheses \textup{\textbf{(HL)}}, \textup{\textbf{(H$\textup{b$\lambda$Q}$)}}, \textup{\textbf{(H$\lambda_0$)}} 
 	and  \textup{\textbf{(H$\textup{fg}$)}}  hold. Then the function $v$  
 	in \eqref{Sec:PDP_def_v} is a  viscosity  solution to \eqref{Sec:PDP_HJB}-\eqref{Sec:PDP_HJB_T}.
 	 \end{theorem}
 	 \proof
	 See Section \ref{Sec:proof_viscprop_rand}. 
	 \endproof
Finally, we provide a comparison theorem for viscosity sub and supersolutions to the first order IPDE  of HJB type \eqref{Sec:PDP_HJB}-\eqref{Sec:PDP_HJB_T} on Hilbert spaces.
To this end, we will need the following additional hypothesis on 
the transition measure $Q$:

 \medskip

 \noindent \textbf{(H$\textbf{Q}$')} \,\, For any $x, x_{\varepsilon} \in E$,  $S_{\varepsilon} \subset E$, such that $x_\varepsilon \rightarrow x$ and   $\cap_{\varepsilon} S_\varepsilon  = \emptyset$,
  \begin{align*}
 \sup_{a \in A}Q(x_\varepsilon, a, S_\varepsilon) \underset{\varepsilon \rightarrow 0}{\rightarrow}0.
 \end{align*}
 \begin{theorem}\label{Sec:PDP_Thm_uniqueness}
 	Let \textbf{\textup{(HL)}}, \textup{\textbf{(H$\textup{b$\lambda$Q}$)}}, \textup{\textbf{(H$\textup{fg}$)}} and  
 	  \textup{\textbf{(H$\textbf{Q}$')}}   hold. Let $u: [0,\,T] \times E \rightarrow \R$ (resp. $v: [0,\,T] \times E \rightarrow \R$) be a bounded and uniformly continuous function in the $|\cdot| \times ||\cdot||_{-1}$ norm, providing a viscosity subsolution  (resp. viscosity supersolution) to \eqref{Sec:PDP_HJB}-\eqref{Sec:PDP_HJB_T}. Suppose that $u(T,x)\leq v(T,x)$ for all $x\in E$. Then $u \leq v$.
 \end{theorem} 
 \proof
See Section \ref{Subsec:compthm}. 
 \endproof
By means of Theorems \ref{Sec:PDP_Thm_existence}, \ref{Sec:PDP_THm_Feynman_Kac_HJB}, \ref{Sec:first_main result}, together with  the comparison Theorem \ref{Sec:PDP_Thm_uniqueness}, we can finally obtain the following probabilistic  representation formula for the value function $V$.
 \begin{theorem}\label{Sec:second_main result}
 Let  \textup{\textbf{(HL)}}, \textup{\textbf{(H$\textup{b$\lambda$Q}$)}}, \textup{\textbf{(H$\lambda_0$)}}, 
\textup{\textbf{(H$\textup{fg}$)}}, 
and 
 \textup{\textbf{(H$\textbf{Q}$')}} hold.
 	Then the function $v$  
 	in \eqref{Sec:PDP_def_v} coincides with the value function $V$,  and  the following Feynman-Kac representation formula holds:
 	\begin{equation}\label{Sec:PDP_equality_value_functions}
 	V(t,x)= Y^{t,x,a}_t,\quad (t,x,a)\in [0,\,T] \times  E \times A.
 	\end{equation}
 	 \end{theorem}

\section{Application to a Hodgkin-Huxley model of neuronal dynamics 
}\label{example}

In the present section we apply our theory to an infinite-dimensional stochastic Hodgkin-Huxley model of neuronal dynamics.  The deterministic Hodgkin-Huxley system was first introduced in \cite{HH}, while  stochastic versions as Hilbert space valued PDMP have been studied in e.g. \cite{Austin}, \cite{Riedler2011},\cite{Ge} and \cite{WainribRiedler}, \cite{ReTreThieu}. 


We  focus on the model considered in \cite{ReTreThieu}.
The axon is modeled by the interval $[0,1]$. We consider ion channels of type $Na$ (sodium) or $K$ (potassium), and we assume that they are located along the axon at positions in $I_N=\frac{1}{N}(\mathbb{Z}\cap N(0,1))$ for some fixed $N\in \mathbb{N^*}$, that we will denote $\frac{i}{N}$ or $z_i$.  
The set of possible states of $K$ and $Na$ channels are denoted respectively by $D_1$ and $D_2$, and are given by
\begin{eqnarray*}
D_1:=\{n_0,n_1,n_2,n_3,n_4\},\quad D_2:=\{m_0h_1,m_1h_1,m_2h_1,m_3h_1,m_0h_0,m_1h_0,m_2h_0,m_3h_0\}.
\end{eqnarray*}
In the control problem new (rhodopsin) channels that are sensitive to light are inserted in the neuron. Such a rhodopsin  channel (denoted by $ChR2$) can have the four possible states $O_1,O_2,C_1,C_2$, among which  $O_1$ and $O_2$ are conductive. Experimentally,  the channel is illuminated and the effect of the illumination is to put the channel in one of its conductive states.
We set  ${\overline D}:= D_1\cup D_2 \cup D_{ChR2}$ with $D_{ChR2}:=\{ O_1,O_2,C_1,C_2\}$, and  ${\overline D}_N:={\overline D}^{I_N}$.

We consider  the Hilbert space $E:=L^2(0,1)$ and the operator $L:=-\Delta$.
The controlled PDMP consists in a set of PDEs written as ODEs in the Hilbert space $E$ indexed by $d\in \bar D_N$, 
\begin{equation}\label{PDEHHsto}
\left\{
\begin{aligned}
&{\dot v}(t)= \frac{1}{C_m}\Delta v(t) + b_d(v(t)),\\
&v(0) = v,\\
&v(t)(0)=v(t)(1)=0, \quad \forall t > 0,
\end{aligned}
\right.
\end{equation}
where the membrane capacitance $C_m>0$  is constant and, for each $(v,d)\in E \times {\overline D}_N$,
\begin{align}
b_d(v) &:= \frac{1}{N}\sum_{i\in I_N} \Big\{{\overline g}_K\mathbf{1}_{\{d_i=n_4\}} ({\overline V}_K-\Phi_i(v))  + {\overline g}_{Na} \mathbf{1}_{\{d_i=m_3h_1\}}({\overline V}_{Na}-\Phi_i(v))+ {\overline g}_l({\overline V}_l-\Phi_i(v)))\nonumber\\
& \qquad\qquad \quad + {\overline g}_{ChR2}( \mathbf{1}_{\{d_i=O_1\}}+\rho  \mathbf{1}_{\{d_i=O_2\}})({\overline V}_{ChR2}-\Phi_i(v))\Big\}\, \,  \phi_{z_i}\label{fdHHChR2},
\end{align}
with 
\begin{equation}\label{Phi}
\Phi_i(v):=\langle v,\phi_{z_i}\rangle, \quad  z_i \in I_N,
\end{equation}
where $\phi_{z_i}$ is  a mollifier function supported on a neighborhood of $z_i$. 
For a channel of type $K$,  
${\overline g}_K > 0$ is the normalized conductance  and ${\overline V}_K\in \mathbb{R}$ is the reversal potential; the same notation holds for $Na, l, ChR2$
($\overline V_l$ and ${\overline g}_l$ denote respectively the leaky reversal  potential and conductance).
The conductance  depends on the number of channels in the conductive state: for $K$ (resp. $Na$)  such a state is  unique, and it  is $n_4$ (resp. $m_3 h_1$).   
The leaky conductance ${\overline g}_l$ remains constant. 
Formula \eqref{Phi} models  the mean value of the membrane potential on a neighborhood of $z_i$.

For any $x = (v,d) \in E \times   {\overline D}_N$, we  denote by $v_t=\phi^d_t(v)$ the corresponding unique mild solution to the PDE (\ref{PDEHHsto}), that models the membrane potential evolution between two successive changes in the channels configuration. 
The transitions take place inside the discrete domain  $\bar D_N$, and correspond to a  continuous-time Markov chain $d_t$. 
Denoting by $(T_n, d_{T_n})$ the jump times and   post-jump location,  the controlled PDMP starting from $x=(v,d) \in E \times \bar D_N$ is 
   \begin{equation*}
   X_s =(v_s,d_s)=
   \left\{
   \begin{array}{ll}
   (\phi_{s-t}^{d}(v), d) \quad &\textup{if}\,\,s \in [t,\, T_{1}),\\
   (\phi^{d_{T_n}}_{s- T_n}(v),d_{T_n})\quad &\textup{if}\,\, s \in [T_n,\, T_{n+ 1}),\,\, n \in \N \setminus \{0\}.
   \end{array}	
   \right.
   \end{equation*}
   
   The control process $\alpha_t$ is proportional to the intensity of light (which is necessarily bounded),  so that we take as control space $A:=[0, a_{\max}]$ with $a_{\max}>0$. 
Introducing a family of smooth functions  $\sigma_{\zeta,\xi} :\mathbb{R}\rightarrow \mathbb{R}_+^{*}$ depending on $(\zeta,\xi)\in  {\overline D}\times  {\overline D}$ 
for all $x= (v,d)\in E\times  {\overline D}_N$, $a \in A$, we define the jump rate function $\lambda : E\times  {\overline D}_N \times A\rightarrow \mathbb{R}_+$    by
\begin{equation}\label{lambdaexplecontrolled}
\lambda((v,d),a) := \sum_{i\in I_N} \sum_{\substack{\xi\in {\overline D},\\ \xi\neq d_i}} \sigma_{d_i,\xi}(\Phi_i(v),a).
\end{equation}
The transition measure $Q:E\times {\overline D}_N\times A\rightarrow \mathcal{P}( {\overline D}_N)$ is such that, for any $x= (v,d)\in E\times  {\overline D}_N$, $a \in A$, the measure $Q((v,d),a,\cdot)$ is supported by the set ${\cal S}$ of $y=(\tilde v,\tilde d)$ such that $\tilde v=v$ (the trajectories of $(v_t)$ are continuous) and $\tilde d$ differs from $d$ only by one component. For $y=(\tilde v,\tilde d)\in {\cal S}$ such that $\tilde d$ differs from $d$ only by its component $i$, 
 \begin{equation}
Q((v,d),a, y) := \sum_{\substack{\xi\in {\overline D},\\ \xi\neq d_i}} \, \, \frac{\sigma_{d_i,\xi}(\Phi_i(v),a)}{\lambda((v,d),a)}\, \, \delta_v(\tilde v)\delta_\xi({\tilde d}_i), \label{jumpmeasurecontrolled}
\end{equation}
if $y\notin {\cal S}$, $Q(x,a;dy):=0 $.  
The  transition functions $\sigma$  from $C_1$ to $O_1$ and from $C_2$ to $O_2$ are assumed to be proportional to the control $\alpha$  
while the other ones are uncontrolled functions.   
  More precisely (see \cite{HH}, \cite{ReTreThieu}): 
  \begin{align*}
&\sigma_{c_1,o_1}(v,a) = \varepsilon_1 a, \quad \sigma_{o_1,c_1}(v,a) = K_{d1}, \quad \sigma_{o_1,o_2}(v,a) = e_{12}, \quad \sigma_{o_2,o_1}(v,a) = e_{21}, \\
&\sigma_{o_2,c_2}(v,a) = K_{d2}, \quad \sigma_{c_2,o_2}(v,a) = \varepsilon_2 a,\quad \sigma_{c_2,c_1}(v,a) = K_r, 
\end{align*} 
and
\begin{align*}
\sigma_{n_0,n_1}(z) = 4\alpha_n(z), \quad \sigma_{n_1,n_2}(z) = 3\alpha_n(z), \quad &\sigma_{n_2,n_3}(z) = 2\alpha_n(z), \quad \sigma_{n_3,n_4}(z) = \alpha_n(z),\\
\sigma_{n_4,n_3}(z) = 4\beta_n(z),  \quad \sigma_{n_3,n_2}(z) = 3\beta_n(z), \quad &\sigma_{n_2,n_1}(z) = 2\beta_n(z), \quad \sigma_{n_1,n_0}(z) = \beta_n(z)\\
\sigma_{m_0h_1,m_1h_1}(z)=\sigma_{m_0h_0,m_1h_0}(z) = 3\alpha_m(z), 
\quad &\sigma_{m_1h_1,m_2h_1}(z) =\sigma_{m_1h_0,m_2h_0}(z) = 2\alpha_m(z),\\
\sigma_{m_2h_1,m_3h_1}(z) = \sigma_{m_2h_0,m_3h_0}(z) =  \alpha_m(z),
\quad &\sigma_{m_3h_1,m_2h_1}(z) = \sigma_{m_3h_0,m_2h_0}(z) = 3\beta_m(z),\\
\sigma_{m_2h_1,m_1h_1}(z) =\sigma_{m_2h_0,m_1h_0}(z) = 2\beta_m(z), 
\quad &\sigma_{m_1h_1,m_0h_1}(z) = \sigma_{m_1h_0,m_0h_0}(z) = \beta_m(z), 
\end{align*}
where 
\begin{eqnarray*}
&&\alpha_n(z)=\frac{0.1-0.01z}{e^{1-0.1z}-1},\, \,  \beta_n(z)=0.125e^{-\frac{z}{80}},\quad \alpha_m(z)=\frac{2.5-0.1z}{e^{2.5-0.1z}-1},\, \,  \beta_m(z)=4e^{-\frac{z}{18}},\\
&&\alpha_h(z)=0.07e^{-\frac{z}{20}},\, \,  \beta_h(z)=\frac{1}{e^{3-0.1z}+1}.
\end{eqnarray*}

The optimal control problem consists in mimicking a desired output reference potential $V_{ref}$, that encodes a given biological behavior while minimizing the intensity of the light applied to the neuron. 
This corresponds to setting, for any $x=(v,d) \in E \times \bar D_N$,  
\begin{align}\label{fgexample}
f(x,a)=f((v,d),a) =\kappa||v-V_{ref}||^2+a,  \quad g(x)&=0, 
\end{align}
so that the cost functional and the value function of the control problem are  
 \begin{align*}
 J(t,x,\alpha) =
 \E^{t,x}_{\alpha}\left[\int_{t}^{T}  \left(\kappa||v_s-V_{ref}||^2+\alpha(X_s)\right)ds\right], \quad 
 V(t,x) = \inf_{\alpha \in \mathcal{A}^t_{ad}}J(t,x,\alpha).
 \end{align*}
The reference signal $V_{ref}$ (that we assume not depending on time) may correspond to a healthy behavior that we want the system to recover thanks to the light stimulation. The intensity of the light is modeled by  the control $\alpha_s = \alpha(X_s)$. Getting the intensity minimal is crucial for the feasability of the experiment in relation to the technical characteristics of the devices that are used.

\begin{remark}
The control  of general infinite-dimensional PDMP is considered in \cite{ReTreThieu}, \cite{Ren}. As in the present paper, in \cite{ReTreThieu} the authors deal with piecewise open loop controls (see \cite{Ver}),   and the control may act on the three characteristics of the PDMP; however, the main tools were relaxed controls and the optimal control of Markov Decision Processes, see e.g. \cite{BeSh}. As an application, other types of models can also be considered: 
 the PDEs in (\ref{PDEHHsto}) may depend on the control variable corresponding to the case where $b_d$ depends on the control, $\phi_{z_i}$ may be replaced by $\delta_{z_i}$ or finally the set $D_{ChR2}$ may have three elements, in which case a $ChR2$ channel has a unique conductive state.
\end{remark}

The rest of this section is devoted to check that the  Hodgkin-Huxley stochastic model  described  above  can be put into the framework of the theory developped in the previous sections. 

\begin{proposition}\label{L:propsemigr} 
\begin{itemize}
\item[(i)]
The operator $L:=-\Delta$  is densely defined, maximal monotone and self-adjoint. Moreover,  $B:=(I-\Delta)^{-1}$ satisfies the {\it strong $B$-condition} with $c_0=1$,  namely $-\Delta B+ B\geq I$ which implies the weak B-condition (\ref{LstarB_property}). 
\item [(ii)] The semigroup $(S(r))_{r \geq 0}:=(e^{-rL})_{r \geq 0}=(e^{r\Delta})_{r\geq 0}$  generated by $L:=-\Delta$ is strongly continuous,  and for all $r>0$, and $S(r)$  is a contraction with respect to $|| \cdot ||$ and also with respect to $|| \cdot ||_{-1}$.
\end{itemize}
\end{proposition}
 \proof (i)
 From \cite{gozswi15}, example   3.14 at page 155 (see also \cite{Ren}) $B:=(I-\Delta)^{-1}$ satifies the {\it strong $B$-condition} with $c_0=1$ namely $-\Delta B+ B\geq I$,  which implies in particular the weak B-condition \eqref{LstarB_property}.
  
(ii) 
For  any $k \in \N$, let us define 
\begin{equation}\label{basisfk}
f_k=\sqrt 2\, \sin(k\pi). 
\end{equation}
$(f_k)_{k\geq 1}$ is an orthonormal basis of $E$, $\Delta f_k=-k^2\pi^2 f_k$,   and, for any $v \in E$, 
\begin{align*}
||v||^2&=\sum_{k\geq 1}\, (v,f_k)^2,\quad ||v||_{-1}^2=((I-\Delta)^{-1} v,v)_H =\sum_{k\geq 1}\,\frac{1}{(1+k^2 \pi^2)}\, (v,f_k)^2.
\end{align*}
Moreover $S(r) =e^{r\Delta}$ is such that  $S(r)v\in {\cal D}(\Delta)$ for all $r>0$, $v\in E$, and satisfies
\begin{align*}
&S(r)v=\sum_{k\geq  1} e^{-rk^2\pi^2}(v,f_k)\, f_k\label{Sr}, \quad  r\geq 0,\,\, v\in E,\\
&||S(r)v||_{-1}^2=\sum_{k\geq 1}\,\frac{1}{(1+k^2 \pi^2)}\, (S(r)v,f_k)^2=\sum_{k\geq 1}\, \frac{1}{(1+k^2\pi^2)}\,  e^{-2rk^2\pi^2}(v,f_k)^2, \quad r\geq 0.  
\end{align*}
We have  $||S(r)v||^2\leq e^{-2r\pi^2}||v||^2$. 
Moreover, 
$||S(r)v||_{-1}^2\leq e^{-2r\pi^2}||v||_{-1}^2$.
\endproof
\begin{lemma}\label{L:Phi}
For any $i \in I_N$, let $\Phi_i$ be the function in \eqref{Phi}.  
Then there exists a positive constant $C_i$ such that,  for all $v,v'$ in $E$,
\begin{equation}\label{Phi i Lip_1-1}
|\Phi_i(v')-\Phi_i(v)|\leq C_i\, || v'-v||_{-1}.
\end{equation}
\end{lemma}
\proof  We have  $\Phi_i(v')-\Phi_i(v)=(v'-v,\phi_{z_i})$, so  taking the basis $(f_k)_{k\geq 1}$ in \eqref{basisfk},
\begin{align*}
(v-v',\phi_{z_i})&=\sum_{k\geq 1}\, (v-v',f_k)(\phi_{z_i},f_k)=\sum_{k\geq 1}\, \frac{1}{\sqrt{1+k^2\pi^2}}(v-v',f_k)\, \sqrt{1+k^2\pi^2}(\phi_{z_i},f_k).
\end{align*}
By the Cauchy-Schwarz inequality,
\begin{align*}
|(v-v',\phi_{z_i})|&\leq\left[\sum_{k\geq 1}\, \frac{1}{(1+k^2\pi^2)}(v-v',f_k)^2\right]^{\frac{1}{2}}\, \left[\sum_{k\geq 1}(1+k^2\pi^2)(\phi_{z_i},f_k)^2\right]^{\frac{1}{2}}\nonumber\\
&=||v'-v||_{-1}\, [((I-\Delta)\phi_{z_i},\phi_{z_i})]^{\frac{1}{2}}.
\end{align*}
It remains to prove that $((I-\Delta)\phi_{z_i},\phi_{z_i})<+\infty$,  so that \eqref{Phi i Lip_1-1} holds with $C_i=[((I-\Delta)\phi_{z_i},\phi_{z_i})]^{\frac{1}{2}}$.
To this end, we take $\phi_{z_i}(z):=\frac{1}{\gamma}M(\frac{z-z_i}{\gamma})$ with $M(z)=\one_{(-1,1)}(z)\, e^{-\frac{1}{1-z^2}}$.  We have $\phi'_{z_i}(z):=\frac{1}{\gamma^2}M'(\frac{z-z_i}{\gamma})$ and $\phi''_{z_i}(z):=\frac{1}{\gamma^3}M''(\frac{z-z_i}{\gamma})$. Moreover
\begin{align*}
M'(\zeta)&=-\frac{2\zeta}{(1-\zeta^2)^2}M(\zeta),\quad M''(\zeta)=M(\zeta)\, \left[\frac{4\zeta^2}{(1-\zeta^2)^4}-2\frac{1+3\zeta^2}{(1-\zeta^2)^3}\right]=M(\zeta)\, \frac{2(3\zeta^4-1)}{(1-\zeta^2)^4}.
\end{align*}
Therefore, setting  $\zeta=\frac{z-z_i}{\gamma}$, 
\begin{equation}
(I-\Delta)\phi_{z_i}(z)=\frac{1}{\gamma}M\left(\frac{z-z_i}{\gamma}\right)-\frac{1}{\gamma^3}M''\left(\frac{z-z_i}{\gamma}\right)=\phi_{z_i}(z)\left(1-\frac{2}{\gamma^2}\frac{(3\zeta^4-1)}{(1-\zeta^2)^4}\right)\label{I-Deltaphi i}.
\end{equation}
\qed
\begin{proposition}\label{vectorfields_flows} 
Let  $d \in \bar D_N$, $v, v' \in E$.
 Then 
\begin{align}
 	&||b(v,d)-b(v',d)||_{-1} \leqslant C ||v-v'||_{-1}\label{Lip-1}.
\end{align}
Moreover, for all $R>0$,  there exists a positive constant $C_R$ such that, for all $a \in A$, 
	\begin{equation}\label{lambda_localLip}
|\lambda((v,d),a)-\lambda((v',d),a)| \leqslant C_R \,||v-v'||_{-1},\quad v,\,v' \in E\,\,  {\rm  s.t.}\, \, ||v|| \vee ||v'||\leq R. 
 	\end{equation} 
 \end{proposition}
 \proof
Let $d \in  \bar D_N$ and $v,v'\in E$.
By \eqref{fdHHChR2} we have 
\begin{equation}\label{formfd}
b(v,d)=\sum_{i\in I_N}\, \gamma_i\, \phi_{z_i}-\sum_{i\in I_N}\, c_i\, \Phi_i(v)\,  \phi_{z_i}. 
\end{equation}
Therefore 
\begin{equation}\label{difference}
|| b(v',d)-b(v,d) ||_{-1}\leq\sum_{i\in I_N}\, c_i\,|\Phi_i(v)-\Phi_i(v')|\,  ||\phi_{z_i}||_{-1}.
\end{equation}
Since $||\phi_{z_i}||_{-1}\leq C_i ||\phi_{z_i}||$ and $I_N$ is a finite set, \eqref{Lip-1} follows from \eqref{difference} and   
Lemma \ref{L:Phi}.

 Let us finally prove \eqref{lambda_localLip}.   
We assume that 
 $||v||\vee||v'||\leq R$.
By definition \eqref{lambdaexplecontrolled},  it is sufficient to check that,  for any $i \in I_N$,  
	\begin{equation*}
|\sigma_{d_i,\xi}(\Phi_i(v),a)-\sigma_{d_i,\xi}(\Phi_i(v'),a)| \leqslant C_R \,||v-v'||_{-1}, 
 	\end{equation*}
which in turn corresponds  to  prove the same property   for  the functions $\alpha_q(\Phi_i(v))$, $\beta_q(\Phi_i(v))$,   $q=n,m,h$.
Recalling  \eqref{Phi} and applying  the Cauchy-Scwartz inequality, we see that   $\Phi_i(v)$, $\Phi_i(v')$ belong to a bounded interval $J_R$ depending on $R$. Then, denoting by $K_{q,R}$ the Lipschitz constant of $\alpha_q$ on $J_R$, from Lemma \ref{L:Phi}
\begin{align*}
|\alpha_q(\Phi_i(v))-\alpha_q(\Phi_i(v'))|&\leq K_{q,R} \,  |\Phi_i(v)-\Phi_i(v')|\leq K_{q,R}\,C_i\, ||v-v'||_{-1}, 
\end{align*}
where $C_i$ is the positive constant in  \eqref{Phi i Lip_1-1}. 
The conclusion follows  recalling  that $I_N$ is a finite set.
\endproof

\begin{lemma}\label{bded given x} For any $d \in \bar D_N$, and  $s, s' \in [t,\,T]$, 
\begin{itemize}
	\item[(i)]$||\phi^d_{s-t}(v)||\leq C(1+||v||)$, \quad $v \in E$, 
	\item[(ii)] $||\phi^d_{s-t}(v)- \phi^d_{s'-t}(v)||\leq C\,\sigma_R(|s-s'|)$,\quad $v \in E: ||v||\leq R$, 
	\item[(iii)] $||\phi^d_{s-t}(v)- \phi^{d}_{s-t}(v')||\leq C\,\omega(||v-v'||)$, \quad $v,v' \in E$, 
	\item[(iv)] $||\phi^d_{s-t}(v)- \phi^d_{s'-t}(v)||_{-1}\leq C\,\sigma_R(|s-s'|)$,\quad $v \in E: ||v||\leq R$, 
	\item[(v)] $||\phi^d_{s-t}(v)- \phi^{d}_{s-t}(v')||_{-1}\leq C\,\omega(||v-v'||_{-1})$, \quad $v,v' \in E$.
\end{itemize}
\end{lemma}
\proof
We first prove (i) and (iii). 
   Setting $S(r) = e^{-r L}$,  the equation for the mild solution to \eqref{PDEHHsto} starting from $x=(v,d)\in E \times \bar D_N$ reads  
 \begin{align*}
 	\phi^d_{s-t}(v) = S(s-t)v + \int_t^s S(s-r) b(\phi^d_{r-t}(v))dr.
 \end{align*}
 Concerning (i), 
using the contraction property of $S(u)$ with respect to $||\cdot||$ given in   Proposition \ref{L:propsemigr}-(ii), we obtain
 \begin{align*}
 	||\phi^d_{s-t}(v)|| \leq ||v|| + \int_t^s ||b(\phi^d_{r-t}(v))|| dr.
 \end{align*}
 On the other hand, recalling (\ref{formfd}),
 \begin{align*}
 	||\phi^d_{s-t}(v)|| \leq ||v|| + \sum_{i\in I_N}\, \int_t^s (|\gamma_i|+|c_i|\, \,  |\Phi_i(\phi^d_{r-t}(v))|) \, dr\, \,  || \phi_{z_i}||.
 \end{align*}
Using Lemma \ref{L:Phi} we get 
 \begin{align*}
 	||\phi^d_{s-t}(v)|| \leq (||v||+CT) +\Gamma \,  \int_t^s ||\phi^d_{r-t}(v)) || dr,
 \end{align*}
and item (i) follows by  Gronwall's Lemma.

Let us now turn to (iii). 
 For any $d\in \bar D_N$, $v, v' \in E$, 
    we  have 
  \begin{equation*}
\phi^d_{s-t}(v)-\phi^d_{s-t}(v') = S(s-t)(v'-v) + \int_t^{s} S(s-r) (b_d(\phi^d_{r-t}(v))-b_d(\phi^d_{r-t}(v')))\, dr.
\end{equation*}
Taking the norm $||\cdot  ||$,  and applying    Proposition \ref{L:propsemigr}-(ii) together with \eqref{Lip-1}, we obtain   
\begin{align*}
||\phi^d_{s-t}(v)-\phi^d_{s-t}(v')||_{-1}
\leq ||v'-v||_{-1} +C \int_0^s ||\phi^d_{r-t}(v)-\phi^d_{r-t}(v')||_{-1} \,dr.
\end{align*}
The conclusion follows again from   the Gronwall Lemma.

Properties (iv) and (v) can be proved analogously, using     the contraction property of $S(u)$ with respect to $||\cdot||_{-1}$ given in Proposition \ref{L:propsemigr}-(ii).  \qed

\paragraph{Additional results on   $V=H_0^1(I)$.} The space $V=H_0^1(I)$ is  continuously embedded in the set of continuous functions on $I$.
For any  $k \in \N$, let us set 
\begin{equation}\label{basisek}
e_k=\frac{\sqrt 2}{\sqrt{1+k^2\pi^2}}\, \sin(k\pi).
\end{equation}
 Then  $(e_k)_{k\geq 1}$ is an orthonormal basis of $V=H_0^1(I)$, and $\Delta e_k=-k^2\pi^2 e_k$. 
 For all $v\in V$, we set $(v,e_k)_V:= \int_0^1 v(z)e_k(z)dz+\int_0^1 v'(z)e'_k(z)dz$. We have
\begin{align}
&||v||_V^2=\sum_{k\geq 1}\, (v,e_k)_V^2, \quad \quad ||v||_{-1,V}^2=((I-\Delta)^{-1} v,v)_V =\sum_{k\geq 1}\,\frac{1}{(1+k^2 \pi^2)}\, (v,e_k)_V^2, \notag
\\
&||S(r)v||_{-1,V}^2=\sum_{k\geq 1}\,\frac{1}{(1+k^2 \pi^2)}\, (S(r)v,e_k)_V^2=\sum_{k\geq 1}\, \frac{1}{(1+k^2\pi^2)_V}\,  e^{-2rk^2\pi^2} \, (v,e_k)_V^2 \quad \forall r\geq 0.\notag 
\end{align}

\begin{remark}\label{R_HV}
Lemma \ref{L:Phi} and Propositions \ref{vectorfields_flows}.  
hold true with $V=H_0^1(I)$ in place of $E = L^2(0,1)$.
\end{remark}

The following  result for the 
 PDEs \eqref{PDEHHsto}, given in   Lemma 4.1 in \cite{ReTreThieu}, plays a fundamental role.
\begin{lemma}\label{boundedflow}  
Set 
$V_-:=\min({\overline V}_{Na}, {\overline V}_{K}, {\overline V}_{L}, {\overline V}_{ChR2})$,  $V_+:=\max({\overline V}_{Na}, {\overline V}_{K}, {\overline V}_{L}, {\overline V}_{ChR2})$, and 
let $d \in \bar D_N$. If $v \in H_0^1(I)$ is continuous in $I=[0,1]$, 
  and  $v(z)\in [V_-,V_+]$ for all $z\in I$, then, for every $d \in \bar D_N$, 
\begin{equation}\label{invariance}
\phi_r^d(v)(z)\in [V_-,V_+], \quad r\in[0,T], \,\,\, z\in I.
\end{equation}
\end{lemma}
Physiologically speaking, we are only interested in the domain $[V_{-}, V_{+}]$. Since Lemma \ref{boundedflow} shows that this domain is invariant for the controlled PDMP, we can modify the local characteristics of the PDMP outside the domain $[V_{-}, V_{+}]$ without changing its dynamics inside of  $[V_{-}, V_{+}]$.  
We will do so for the rate functions $\sigma_{d_i, \xi}$. From now on, consider a compact set $K$ containing the closed ball of $E$, centered in $0$ with radius $\max\{V_{-}, V_{+}\}$. We will rewrite $\sigma_{d_i, \xi}$ outside $K$ such that they all become Lipschitz and  bounded functions. We also  take $V_{\rm ref}$  taking values in $K$ and $\tilde f$ bounded and globally Lipschitz such that 
\begin{equation}\label{bounded f}
\tilde f(v)=||v-V_{\rm ref}||^2, \quad \forall v\in K.
\end{equation} 
Since the control set $A=[0,a_{\rm max}]$ is bounded, the corresponding value function and cost 
  \begin{align*}
 \tilde J(t,x,\alpha) =
 \E^{t,x}_{\alpha}\left[\int_{t}^{T}  \left(\kappa\tilde f(v_s)+\alpha(X_s)\right)ds\right], \quad 
 V(t,x) = \inf_{\alpha \in \mathcal{A}^t_{ad}}\tilde J(t,x,\alpha),
 \end{align*}
are bounded as well.

The next two results show that the case of  the stochastic controlled infinite-dimensional Hodgkin-Huxley model can be actually covered by the theory on controlled infinite-dimensional PDMPs developed in the present paper.  

\begin{proposition}\label{locLiplambdaV_fV} 
Let $v, v' \in V$ such that 
$v(z)$ and $v'(z)$ belong to $[V_-,V_+]$ for all $z\in [0,1]$. 
The following hold. 
\begin{itemize}
\item[(i)]	
There exist a  positive constants $C_1$ such that, for all  $d\in \bar D_N$,  and  $a \in A$,
	\begin{equation}\label{lambdaV}
 |\lambda((\phi_s^d(v),d),a)-\lambda((\phi^d_s(v'),d),a)| \leq C_1 \,||v-v'||_{-1,V}, \quad r\in [0,T].
	\end{equation} 
\item [(ii)] If in addition $||(I-\Delta)V_{ref}||<+\infty$,
 there exists a positive constant $C_2$ such that, for all $d\in \bar D_N$, $a \in A$, the function $f$ in \eqref{fgexample} satisfies 
	\begin{equation}\label{f_lipschitz}
 |f(\phi^d_s(v),a)-f(\phi^d_s(v'),a)| \leq C_2 \,||v-v'||_{-1,V},\quad r\in [0,T]. 
 	\end{equation} 
 	\end{itemize}
\end{proposition}

\proof
Let us prove item (i). 
Recalling \eqref{Phi} and using  the Cauchy-Schwarz inequality we have  
\begin{equation}\label{normePhiduflot}
|\Phi_i(\phi_s^d(v))|\leq ||\phi^d_s(v)||\,||\phi_{z_i}||.
\end{equation}
Since
\begin{equation}\label{borne ineq}
||\phi_s^d(v)||\leq ||\phi_s^d(v)||_\infty\leq \max\{|V_- |,|V_+|\},
\end{equation}
and the same inequalities hold for $\phi_s^d(v')$, we have  
\begin{equation}
|\alpha_q(\Phi_i(\phi^d_s(v)))-\alpha_q(\Phi_i(\phi^d_s(v')))|\leq K_{q,R} \,  |\Phi_i(\phi^d_s(v))-\Phi_i(\phi^d_s(v'))|,
\end{equation}
with $R=\max\{|V_ -|,|V_+|\}$ and $K_{q,R}$ the  Lipschitz constant of $\alpha_q$ depending on $R$. Taking into account Remark \ref{R_HV}, we conclude by applying the $V$-versions of Lemmas \ref{L:Phi} 
and \ref{vectorfields_flows}.

Let us now consider item (ii). Using  the basis $(e_k)$ introduced in \eqref{basisek}, and applying the Cauchy-Schwarz inequality, 
	\begin{align}\label{est f}
 &|f(\phi^d_s(v),a)-f(\phi^d_s(v'),a)|=\kappa \, \Big|\sum_{k\geq 1} ((\phi^d_s(v)-V_{ref},e_k)^2-(\phi^d_s(v')-V_{ref},e_k)^2)\Big|\\
&\leq \kappa \, \sum_{k\geq 1}  |(\phi^d_s(v_0)-\phi^d_s(v'),e_k)|\, |(\phi^d_s(v)+\phi^d_s(v')-2V_{ref},e_k)|\leq \kappa \, ||\phi^d_s(v)-\phi^d_s(v')||_{-1,V}\, \cal T,\notag 
\end{align} 
where 
\begin{equation*}
{\cal T}:=\sum_{k\geq 1} (1+k^2\pi^2)(\phi^d_s(v)+\phi^d_s(v')-2V_{ref},e_k)^2 =||(I-\Delta)(\phi^d_s(v)+\phi^d_s(v')-2V_{ref})||_V^2.
\end{equation*} 
By Proposition  \ref{vectorfields_flows} and Remark \ref{R_HV},   it remains to study the boundedness properties of $\cal T$. 
Since by assumption $||(I-\Delta)V_{ref}||<+\infty$, 
\begin{equation}\label{Tau_est}
{\cal T}\leq C(||(I-\Delta)\phi^d_s(v)||_V^2+||(I-\Delta)\phi^d_s(v')||_V^2+||(I-\Delta)V_{ref}||_V^2).
\end{equation}
Let us thus consider the term 
$||(I-\Delta)\phi_d(s,v)||_V$.  
Being $(I-\Delta)$  linear,
we can write
  \begin{align*}
(I-\Delta)\phi^d_s(v)&= (I-\Delta)S(s)v_0 + \int_0^{s} (I-\Delta)S(s-r) b_d(\phi^d_r(v),a)\, dr.
\end{align*}
Moreover, since $(I-\Delta)$ and $S(r)$ commute, 
  \begin{align}\label{I-Delta step}
||(I-\Delta)\phi^d_s(v)||_V&
\leq ||(I-\Delta)v||_V + \int_0^{s} \sum_{i\in I_N}\, (|\gamma_i|+|c_i|\, |\Phi_i(\phi^d_r(v))|)\, ||(I-\Delta) \phi_{z_i}||_V\, dr,
\end{align}
where  we have used that (recall formula \eqref{formfd})
\begin{equation*}
(I-\Delta)b_d(\phi^d_r(v_0),a)=\sum_{i\in I_N}\, \gamma_i\, (I-\Delta)\phi_{z_i}-\sum_{i\in I_N}\, c_i\, \Phi_i(\phi^d_r(v_0))\, (I-\Delta) \phi_{z_i}.
\end{equation*}
Recalling (\ref{normePhiduflot}) and (\ref{borne ineq}), \eqref{I-Delta step} yields 
\begin{align*}
||(I-\Delta)\phi^d_s(v)||_V\leq ||(I-\Delta)v||_V + \int_0^s \sum_{i\in I_N}\, (|\gamma_i|+|c_i|\,  \max\{|V_- |,|V_+|\}||\phi_{z_i}||)\, ||(I-\Delta) \phi_{z_i}||_V \, dr.
\end{align*}
Recalling \eqref{I-Deltaphi i} we see that, for any $i \in I_N$,  $||(I-\Delta) \phi_{z_i}||_V\leq C_i$. Since $I_N$ is finite, we conclude from the above inequality 
that there exists some   constant $\Gamma$ such that
  \begin{equation}\label{I-Delta step2}
||(I-\Delta)\phi_d(s,v)||_V\leq ||(I-\Delta)v||_V + \Gamma \, T;
\end{equation}
analogous  inequalities holds true for $\phi^d_s(v')$ and $v'$. Then \eqref{Tau_est}, together with \eqref{I-Delta step2}, yields
\begin{equation}\label{00}
{\cal T}\leq 2||(I-\Delta)v||_V + 2\, \Gamma \, T+4\,  ||(I-\Delta)V_{ref})||_V^2.
\end{equation}
and the conclusion follows.
\endproof

\begin{proposition}\label{boundednessVflow_lambda} Let $v_0\in V$ such that 
$v_0(z)\in [V_-,V_+]$ for all $z\in [0,1]$.
Then there exist
two 
positive constants $C_1,C_2$, 
 only depending on $T, N, \max\{|V_- |,|V_+|\})$, such that, for all $d \in D$, $a \in A$, 
	\begin{align}
& ||\phi^d_s(v_0)||_V \leq C_1, \quad s\in [0,T],\label{bound_phi}\\
 &|f((\phi^d_s(v_0),d),a)|+|\lambda((\phi^d_s(v_0),d),a)|+||b(\phi^d_s(v_0),d)||_V \leq C_2, \quad s\in [0,T].\label{Bound}
 	\end{align}  
\end{proposition}
\proof 
Estimate \eqref{bound_phi}
is obtained arguing as in  Lemma \ref{bded given x}-(i).
The boundedness of $f(\phi^d_s(v_0),a)$  follows  from 
\eqref{bound_phi},
 recalling that
$$
 |f(\phi^d_s(v_0),a)|=\kappa \, \Big|\sum_{k\geq 1} ((\phi^d_s(v_0)-V_{ref},e_k)^2\Big|\leq\kappa(||\phi^d_s(v_0)||^2_V + V_{ref}^2 -2 V_{ref}||\phi^d_s(v_0)||_V).
$$
On the other hand, recalling \eqref{formfd} and \eqref{normePhiduflot}, 
\begin{align*}
||b_d(\phi^d_s(v_0))||_V 
&\leq \sum_{i\in I_N}\, |\gamma_i|\, ||\phi_{z_i}||_V+\sum_{i\in I_N}\, |c_i|\, ||\phi^d_s(v_0)||\,||\phi_{z_i}||\,  ||\phi_{z_i}||_V,
\end{align*}
and we obtain the bound   from 
Lemma \ref{bded given x}-(i) and the fact that $||\cdot||\leq ||\cdot||_V$.

The boundedness of $\lambda_d(\phi^d_s(v_0), a)$ follows from the form of the functions $\alpha_q$, $\beta_q$, together with (\ref{normePhiduflot}) and the fact that $||\cdot||\leq ||\cdot||_V$.
\qed

\section{Proofs of the  results in Section \ref{Sec:PDP_Sec_control_problem}}\label{Sec_mainproofs_Sec2}

\subsection{Proof of Proposition \ref{P:dynprogpr}}\label{Sec_proof_contV}
The boundedness of $V$ directly comes from \eqref{Sec:PDP_value_function} and the boundedness of $f$ and $g$.  
Let $B([0,\,T] \times E)$ be the set of all bounded functions on $[0,\,T] \times E$,   and define the map $\mathcal T: B([0,\,T] \times E) \rightarrow B([0,\,T] \times E)$ as 
\begin{align*}
\mathcal T \psi(t,x) &:=  \inf_{\alpha \in \mathcal A_{ad}^t}  \E^{t,x}_{\alpha}\left[\int_{t}^{T_1 \wedge T} f(X_s,\alpha_s)\, ds  +   g(X_T)\one_{T \leq T_1} +\psi(T_1, X_{T_1})\one_{T >T_1}\right].
\end{align*}
Set $\mathcal U =\{u : [0,\,+ \infty) \rightarrow A\,\,\textup{measurable}\}$.
One can  show that 
\begin{align}\label{new_G}
\mathcal T \psi(t,x) & =\inf_{u\in \mathcal U} \bigg\{\int_0^{T-t}\chi^{u}(s, x) (f^{u}(s,  x)+L_{\psi}^{u}(s,x))\, ds  + \chi^{u}(T- t, x) g(\phi^{u}(T-t,x)\bigg\},
\end{align}
where $\chi^{u}(s, x)= e^{-\int_0^{s}  \lambda(\phi^{u}(r, x), u_r)\,dr}$, $f^{u}(s,  x) = f(\phi^{u}(s, x), u_s)$,  and 
\begin{align*}
L_{\psi}^{u}(s,  x)&=  \int_E \psi(s,y)\,\lambda(\phi^{u}(s, x), u_s)\,Q(\phi^{u}(s, x), u_s,dy).
\end{align*} 
It can be proved that $\mathcal T$ is a contracting map in $B([0,\,T] \times E)$ and $V$ is its unique fixed point, see e.g.  Theorem 4.6 in \cite{CalviaInf} or Theorem 3.3 and Lemma 3.4 in \cite{ReTreThieu}. 
In particular, $V$ satisfies the DPP \eqref{DynProgPr}.

 Denote by $C_b([0, T] \times E)$ the set of bounded functions, continuous on $[0,\,T] \times E$  with the $|\cdot| \times ||\cdot||_{-1}$ norm.  In order to prove that $V \in C_b([0, T] \times E)$, it suffices to show that for any function  $\psi  \in C_b([0, T] \times E)$ one has $\mathcal T \psi \in C_b([0, T] \times E)$. As a matter of fact, we know that  
 $\mathcal T$ is a contracting map  in  $B([0,\,T] \times E)$ and that $V$ is its unique fixed point, namely  $V = \mathcal T V$. Assume now that,  for any function  $\psi  \in C_b([0, T] \times E)$, one has $\mathcal T \psi \in C_b([0, T] \times E)$. Then $\mathcal T$ is a contracting map in $C_b([0, T]$
 and has a unique fixed point in $C_b([0, T]$, that we denote   $w^\ast$. We have 
 $$
 ||w^\ast -V||_{\infty}=||\mathcal Tw^\ast -\mathcal TV||_{\infty}\leq \rho||w^\ast -V||_{\infty}, \quad \rho \in (0,1), 
 $$
 so that $ ||w^\ast -V||_{\infty}=0$.

In the following $C$ will denote a generic constant, that may vary from line to line, and that may depend  on $T$.
We start by noticing that,  by \eqref{new_G}, 
 	$
 	\mathcal T \psi(t,x) = \inf_{u \in \mathcal U}\bar J(t,x, u),
 	$
 	with
 	\begin{align*}
\bar J(t,x, u)  
& = \int_0^{T-t}\chi^{u}(s, x) (f^u(s,x) + L_{\psi}^{u}(s,x))\, ds + \chi^{u}(T-t, x) g(\phi^{u}(T-t, x)).
\end{align*}
Let $t,t',s \in [0,\,T]$, $t'\leq t \leq s$, $x,x' \in E$, $u \in \mathcal U$. Recalling hypotheses \textbf{(H$\textup{b$\lambda$Q}$)}-(i),  \textup{\textbf{(H$\textup{fg}$)}} and \eqref{contrflowestimate}, we have $|\chi^{u}(s, x)| \leq 1$, $|f^{u}(s, x)| \leq C$, and, for any $s' \leq s$,  
\begin{align}
&|\chi^u(s',x)-\chi^u(s,x') |\leq 
(1 - e^{-C ||x-x'||_{-1}})+(1 - e^{-C (s-s')}),\label{chiest}\\
&|f^{u}(s, x)-f^{u}(s, x')| \leq C ||x-x'||_{-1}.\label{fest},\\
&|g(\phi(T-t, x))- g(\phi(T-t', x'))|\leq \omega(||x-x'||_{-1})\label{gest1}.
\end{align}
On the other hand,  by  \textup{\textbf{(H$\textup{b$\lambda$Q}$)}}-(i)-(ii), together with the boundedness and continuity of $\psi$,  we have  $|L_{\psi}^{u}(s, x)\leq C$ and, for $s < T-t$, 
\begin{align}
&|L_{\psi}^{u}(s, x)-L_{\psi}^{u}(s, x')|\leq |\lambda(\phi^{u}(s, x),u_s)-\lambda(\phi^{u}(s, x'),u_s)|\,||\psi||_{\infty}\notag\\
&+||\lambda||_{\infty}\bigg|\int_E \psi(s,y)\,[Q(\phi^{u}(s, x),u_s,dy)-Q(\phi^{u}(s, x'),u_s,dy)]\bigg|\notag\\
&\leq C \,\sigma(||\phi^{u}(s,x)-\phi^{u}(s,x')||_{-1})\leq C \omega(||x-x'||_{-1}),\label{Lest1}
\end{align}
where the latter inequality follows from \eqref{contrflowestimate-BIS}.
Then, for any $t, t' \in [0,\,T]$, $x, x' \in E$,  $u \in \mathcal U$, 
\begin{align*}
&|\bar J(t,x, u)-\bar J(t',x', u)|\notag\\
&\leq\left|\int_0^{T-t}\chi^u(s,x) f^u(s,x)\,ds-\int_{0}^{T-t'}\chi^u(s,x') f^u(s,x')\,ds\right|\notag\\
 	&+ \left|\int_0^{T-t} \chi^u(s,x) L^u_\psi(s,x)\,ds-\int_{0}^{T-t'} \chi^u(s,x') L^u_\psi(s,x')\,dr\right|\\
 	& + |\chi^{u}(T-t, x) g(\phi^{u}(T-t, x)) - \chi^{u}(T-t', x') g(\phi^{u}(T-t', x')) |\notag\\
&\leq \int_{0}^{T-t}|\chi^u(s,x) f^u(s,x)- \chi^u(s,x') f^u(s,x')|\,ds\notag\\	
&+ \int_{0}^{T-t} |\chi^u(s,x) L^u_\psi(s,x)- \chi^u(s,x') L^u_\psi(s,x')|\,ds +  C |t-t'|\notag\\
 	& + C|g(\phi^{u}(T-t, x)) -  g(\phi^{u}(T-t', x')) |+ C|\chi^{u}(T-t, x)  - \chi^{u}(T-t', x') |\notag\\
&\leq C \Bigg(\int_{0}^{T-t}|\chi^u(s,x))-\chi^u(s,x')|\,ds+\int_{0}^{T-t} |f^u(s,x)-f^u(s,x')|\,ds\notag\\	
&+ \int_{0}^{T-t} | L^u_\psi(s,x)  -  L^u_\psi(s,x')| \,ds  + |g(\phi^{u}(T-t, x)) -  g(\phi^{u}(T-t', x')) |\\
&+ |\chi^{u}(T-t, x)  - \chi^{u}(T-t', x') |+   |t-t'|\Bigg)\\
&\leq C (\omega(t-t')+ \omega'(||x-x'||_{-1}))
\end{align*}
for some modulus of continuity $\omega$, $\omega'$, where the latter inequality follows from  \eqref{chiest}, \eqref{fest}, \eqref{gest1}, \eqref{Lest1}.
This  shows in particular   that $V$ is uniformly continuous  in the $|\cdot| \times ||\cdot||_{-1}$ norm.
 \endproof

\subsection{Proof of Proposition \ref{P:DPPNEW}} \label{proofP:DPPNEW} 
We first show that the left-hand side of \eqref{DPP2} is smaller than the right-hand side.   
To this end, let us fix $\alpha \in \mathcal A_{ad}^t$. By \eqref{Sec:PDP_open_loop_controls}, we have that, under $\P^{t,x}_\alpha$,  
$
\theta = \tau_d \wedge T \wedge T_1, 
$
with 
$$
\tau_d :=\inf \left\{s \geqslant t: 
(s,\phi^{\alpha_0}(s-t,x)) \notin B\left((t,x); 
\rho\right)
\right\}. 
$$
From the selection theorem (see e.g.  Proposition 7.50  in \cite{BeSh}), 
for any $\varepsilon >0$,  there exists  a Borel-measurable map  $\gamma^\varepsilon \in \mathcal A_{ad}^\theta$ of the form 
 \begin{align*}
 \gamma_s^\varepsilon(\theta)
&:=\beta_0^\varepsilon(s- \tau_d, \phi^{\alpha_0}(\tau_d-t,x))\,\one_{(\theta,\,T_1]}(s)+ \sum_{n=1}^{\infty}\beta_n^\varepsilon(s-T_n,E_n)\,\one_{(T_n,\,T_{n+1}]}(s),\quad s \in [\theta,\,T], 
 \end{align*}
such that   $\gamma^{\varepsilon}(\theta)$  is an $\varepsilon$-optimal control for $V(\theta, X_\theta)$, namely 
 \begin{equation}\label{eps_opt}
 	V(\theta, X_{\theta}) \geq J(\theta, X_{\theta}, \gamma^{\varepsilon}(\theta)) - \varepsilon, \quad \P^{t,x}_{\alpha}\textup{-a.s.}
 \end{equation}
Set 
   \begin{equation*}
   \bar \alpha^\varepsilon_s=
   \left\{
   \begin{array}{ll}
   \alpha^0_s \quad &\textup{if}\,\,s \in (t,\, \theta],\\
   \gamma_s^\varepsilon(\theta)\quad &\textup{if}\,\, s \in (\theta,\, T],
   \end{array}	
   \right.
   \end{equation*}
and define (by a slight abuse of notation),   for any real  $\eta \geq \tau_d$ and  $\chi \in E$, 
\begin{align*}
\gamma_s^\varepsilon(\eta, \chi)&:= \beta_0^\varepsilon(s- \eta, \chi)\,\one_{(\eta,\,T_1]}(s)+ \sum_{n=1}^{\infty}\beta_n^\varepsilon(s-T_n,E_n)\,\one_{(T_n,\,T_{n+1}]}(s),\quad s \in [\eta,\,T].
\end{align*} 
We have 
 $\bar \alpha^\varepsilon_s = \alpha^0_s \one_{s \in (t,\, \tau_d]} + \gamma_s^\varepsilon(\tau_d, X_{\tau_d}) \one_{s \in (\tau_d,\, T]}$,  if $\theta = \tau_d$, and  $\bar \alpha^\varepsilon_s = \alpha^0_s \one_{s \in (t,\, T_1]}+ \gamma_s^\varepsilon(T_1, E_1) \one_{s \in (T_1,\, T]}$ if $\theta =  T \wedge T_1$. Therefore, 
  $\bar \alpha^\varepsilon \in \mathcal A^t_{ad}$. 
In particular, 
   \begin{align}\label{first}
	V(t,x)\leq J(t,x,\bar \alpha^\varepsilon)&= \E^{t,x}_{\bar \alpha^\varepsilon}\left[\int_{t}^{\theta}  f(X_{s},\bar \alpha_s^\varepsilon)\,ds + \int_{\theta}^{T}  f(X_{s},\bar \alpha_s^\varepsilon)\,ds  + g(X_T)\right]\notag\\
	&= \E^{t,x}_{\bar \alpha^\varepsilon}\left[\int_{t}^{\theta}  f(X_{s},\bar \alpha_s^\varepsilon)\,ds\right] + \E^{t,x}_{\bar \alpha^\varepsilon}\left[\left(\int_{\tau_d}^{T}  f(X_{s},\gamma^\varepsilon_s(\tau_d, X_{\tau_d}))\,ds  + g(X_T)\right) \one_{\tau_d < T_1 \wedge T}\right]\notag\\
	& + \E^{t,x}_{\bar \alpha^\varepsilon}\left[\left(\int_{T_1}^{T}  f(X_{s},\gamma^\varepsilon_s(T_1, E_1))\,ds  + g(X_T)\right)\one_{\tau_d \geq T_1 \wedge T}  \right]\notag\\
	&= \E^{t,x}_{\bar \alpha^\varepsilon}\left[\int_{t}^{\theta}  f(X_{s},\bar \alpha_s^\varepsilon)\,ds\right] + I + II. 
	\end{align}
At this point,  we aim at proving   that  
      \begin{align}
I  
&=J(\tau_d,\phi^{\alpha^0}(\tau_d-t,x), \gamma^\varepsilon(\tau_d, X_{\tau_d}) \, \P^{t,x}_{\bar \alpha^\varepsilon}\left[\tau_d < T_1 \wedge T\right],  \label{Jid_I}\\
II &=  \E^{t,x}_{\bar \alpha^\varepsilon}\left[\one_{\tau_d \geq T_1 \wedge T} \,J(\eta, \chi, \gamma^\varepsilon(\eta, \chi))|_{\eta = T_1, \chi= E_1}\right] \label{Jid_II}. 
	\end{align}
We notice that 
 \begin{align}
		&\E^{t,x}_{\bar \alpha^\varepsilon}\left[\int_{T_1}^{T}  f(X_{s},\gamma^\varepsilon_s(T_1, E_1))\,ds  + g(X_T) \Big| \mathcal F_{T_1}\right] \notag\\
		&= \lim_{n \rightarrow \infty} \E^{t,x}_{\bar \alpha^\varepsilon}\left[\int_{T_1}^{T_n \wedge T}  f(X_{s},\gamma^\varepsilon_s(T_1, E_1))\,ds  + g(X_T)\one_{T_n >T}\Big| \mathcal F_{T_1}\right].\label{lim2}
	\end{align}
	Moreover, setting for any $(\eta, \chi) \in [0,\,\infty) \times E$
	\begin{align*}
	&J^n(\eta, \chi, \gamma^\varepsilon(\eta, \chi)) := 
	\E^{\eta,\chi}_{\gamma^\varepsilon(\eta, \chi)}\left[\int_{\eta}^{T_{n} \wedge T}  f(X_{s},\gamma^\varepsilon_s(\eta, \chi))\,ds  + g(X_T)\one_{T_{n} >T} \right], 
	\end{align*}
	we see that  
	\begin{align}
	&J(\eta, \chi, \gamma^\varepsilon(\eta, \chi))= \lim_{n \rightarrow \infty} J^n(\eta, \chi, \gamma^\varepsilon(\eta, \chi)). \label{lim1}
	\end{align}
Since
   \begin{align}\label{second_J}
&II = \E^{t,x}_{\bar \alpha^\varepsilon}\left[\one_{\tau_d \geq T_1 \wedge T}  \left(\int_{T_1}^{T}  f(X_{s},\gamma^\varepsilon_s(T_1, E_1))\,ds  + g(X_T)\right)\right] \notag\\
&= \E^{t,x}_{\bar \alpha^\varepsilon}\left[\one_{\tau_d \geq T_1 \wedge T} \,\E^{t,x}_{\bar \alpha^\varepsilon}\left[\int_{T_1}^{T}  f(X_{s},\gamma^\varepsilon_s(T_1, E_1))\,ds  + g(X_T) \Big| \mathcal F_{T_1}\right]\right], 
	\end{align}
it follows  from 
\eqref{lim2} and \eqref{lim1} that \eqref{Jid_I} and \eqref{Jid_II} hold  true if and only if, for any $n > 1$, 
\begin{align}
&\E^{t,x}_{\bar \alpha^\varepsilon}\left[\left(\int_{\tau_d}^{T_{n-1} \wedge T}  f(X_{s},\gamma^\varepsilon_s(\tau_d, X_{\tau_d})\,ds  + g(X_T)\one_{T_{n-1} >T}\right) \one_{\tau_d < T_1 \wedge T}\right]\notag\\
&=J^{n-1}(\tau_d,X_{\tau_d}, \gamma^\varepsilon(\tau_d, X_{\tau_d})) \, \P^{t,x}_{\bar \alpha^\varepsilon}\left[\tau_d < T_1 \wedge T\right],  \label{ToProve1}\\
&\E^{t,x}_{\bar \alpha^\varepsilon}\left[\int_{T_1}^{T_n \wedge T}  f(X_{s},\gamma^\varepsilon_s(T_1, E_1))\,ds  + g(X_T)\one_{T_n >T}\Big| \mathcal F_{T_1}\right]= J^{n-1}(\eta, \chi, \gamma^\varepsilon(\eta, \chi))|_{\eta = T_1, \chi= E_1}. \label{ToProve2}
	\end{align}
We show for simplicity  the case  $n=2$, the other cases are obtained analogously.
Concerning   \eqref{ToProve1},  we have 
\begin{align*}
&\E^{t,x}_{\bar \alpha^\varepsilon}\left[\left(\int_{\tau_d}^{T_{1} \wedge T}  f(X_{s},\gamma^\varepsilon_s(\tau_d, X_{\tau_d}))\,ds  + g(X_T)\one_{T_1 >T}\right) \one_{\tau_d < T_1 \wedge T}\right]\notag\\
&= e^{- \int_t^{\tau_d} \lambda(\phi^{\alpha_0}(r- t, x),\alpha_0(r-t,x))  \,dr } \int_{\tau_d}^{\infty} du \Big\{\left(\int_{\tau_d}^{u\wedge T}  f(\phi^{\beta_0^\varepsilon}(s- \tau_d, X_{\tau_d}),\beta_0^\varepsilon(s- \tau_d, X_{\tau_d}))\,ds\right)  \cdot \\
&\cdot\lambda(\phi^{\beta_0^\varepsilon}(u- \tau_d, X_{\tau_d}),\beta_0^\varepsilon(u- \tau_d, X_{\tau_d})   \, e^{- \int_{\tau_d}^u  \lambda(\phi^{\beta_0^\varepsilon}(r- \tau_d, X_{\tau_d}),\beta_0^\varepsilon(r- \tau_d, X_{\tau_d}))  \,dr }\\
&+ g(\phi^{\beta_0^\varepsilon}(T- \tau_d, X_{\tau_d}))\, e^{- \int_{\tau_d}^T  \lambda(\phi^{\beta_0^\varepsilon}(r- \tau_d, X_{\tau_d})),\beta_0^\varepsilon(r-\tau_d, X_{\tau_d})))  \,dr }\Big\} \\
&= \P^{t,x}_{\bar \alpha^\varepsilon}\left[\tau_d < T_1 \wedge T\right] \, J^{1}(\tau_d,X_{\tau_d}, \gamma^\varepsilon(\tau_d, X_{\tau_d})). 
	\end{align*}

On the other hand,  \eqref{ToProve2} with $n = 2$ reads 
	\begin{align*}
		&\E^{t,x}_{\bar \alpha^\varepsilon}\left[\int_{T_1}^{T_2 \wedge T}  f(X_{s},\gamma^\varepsilon_s(T_1, E_1))\,ds  + g(X_T)\one_{T_2 >T}\Big| \mathcal F_{T_1}\right]\\
		&= 
		  \int_{\eta}^\infty  \Big\{ \int_E\left(\int_{\eta}^{u \wedge T}  f(\phi^{\beta_1^\varepsilon}(s- \eta, \chi),\beta_1^\varepsilon(s-\eta,\chi))\,ds\right) e^{- \int_\eta^u  \lambda(\phi^{\beta_1^\varepsilon}(r- \eta, \chi),\beta_1^\varepsilon(r-\eta,\chi))  \,dr }\cdot  \\
		&\cdot\lambda(\phi^{\beta_1^\varepsilon}(u- \eta, \chi),\beta_1^\varepsilon(u-\eta,\chi)) \,Q(\phi^{\beta_1^\varepsilon}(u- \eta, \chi),\beta_1^\varepsilon(u-\eta,\chi), dy)  \, \\
		&+ g(\phi^{\beta_1^\varepsilon}(T- \eta, \chi))\, e^{- \int_\eta^T  \lambda(\phi^{\beta_1^\varepsilon}(r- \eta, \chi),\beta_1^\varepsilon(r-\eta,\chi))  \,dr }\Big\}\,du\,\Big|_{\eta = T_1, \chi = E_1} \\
		&= J^1(\eta, \chi, \gamma^\varepsilon(\eta, \chi))|_{\eta = T_1, \chi= E_1}.
	\end{align*} 
  Thus \eqref{Jid_I}   and \eqref{Jid_II} hold true, 
 and  \eqref{first} gives  
   \begin{align*}
	V(t,x)&\leq 
		\E^{t,x}_{\bar \alpha^\varepsilon}\left[\int_{t}^{\theta}  f(\phi^{\alpha^0}(s-t,x),\alpha^0_s)\,ds + J(\eta,\chi, \gamma^\varepsilon)|_{\eta= \theta, \chi= X_{\theta}}\right]\\
	&=	\E^{t,x}_{\alpha}\left[\int_{t}^{\theta}  f(\phi^{\alpha^0}(s-t,x),\alpha^0_s)\,ds + J(\eta,\chi, \gamma^\varepsilon)|_{\eta= \theta, \chi= X_{\theta}}\right]\\
	&   \leq \E^{t,x}_{\alpha}\left[\int_{t}^{\theta}  f(X_s,\alpha_s)\,ds + V(\eta,\chi)|_{\eta= \theta, \chi= X_{\theta}}\right] + \varepsilon
	\end{align*}
	where we have also used  the fact that $\E^{t,x}_{\bar \alpha^\varepsilon} \left[\varphi(\theta, X_\theta)\right]=\E^{t,x}_{\alpha} \left[\varphi(\theta, X_\theta)\right]$ for any measurable function $\varphi$, since $\bar \alpha^\varepsilon_s$ and  $\alpha_s$ coincide on $(t, \theta] \subset (t, T_1]$.
	The result  follows from the arbitrariness of $\varepsilon >0$ and $\alpha \in \mathcal A_{ad}^t$.

It remains to prove that the left-hand side of \eqref{DPP2} is greater than the left-hand side. 
	To this end, let $\alpha \in \mathcal A_{ad}^t$. 
By  
	using 
	the same argument as in  the previous step, with $\gamma_s^\varepsilon(\eta, \chi)$ replaced by 
	 \begin{align*}
\alpha_0(s- \eta, \chi)\,\one_{(\eta,\,T_1]}(s)+ \sum_{n=1}^{\infty}\alpha_n(s-T_n,E_n)\,\one_{(T_n,\,T_{n+1}]}(s),\quad s \in [\eta,\,T], 
 \end{align*}
	we have 
	\begin{align*}
	J(t,x, \alpha)&= \E^{t,x}_{\alpha}\left[\int_{t}^{\theta}  f(X_{s},\alpha_s)\,ds + J(\eta,\chi, \alpha)|_{\eta = \theta, \chi = X_{\theta}}\right] \\
	&\geq \E^{t,x}_{\alpha}\left[\int_{t}^{\theta}  f(X_{s},\alpha_s)\,ds + V(\eta, \chi)|_{\eta= \theta, \chi = X_{\theta}}\right]\\
	&\geq \inf_{\alpha \in \mathcal{A}^t_{ad}}\E^{t,x}_{\alpha}\left[\int_{t}^{\theta}  f(X_{s},\alpha_s)\,ds + V(\eta, \chi)|_{\eta= \theta, \chi = X_{\theta}}\right], 
	\end{align*}
and the result follows by taking the infimum over $\alpha \in \mathcal A_{ad}^t$ in the left-hand side term.	
\qed

\subsection{Proof of Lemma \ref{Gpsi}}\label{Proof Lemma Gpsi}
Let us fix $t \in (\varepsilon,\,T-\varepsilon)$, $\varepsilon >0$. We first prove that  the map
\begin{align*}
r \mapsto& 
-\frac{\partial \psi}{\partial s} (r,\phi^{\alpha_0}(r-t,x))+\langle \phi^{\alpha_0}(r-t,x),\,L^\ast \, D \varphi(r,\phi^{\alpha_0}(r-t,x)) \rangle \\
&+\langle \phi^{\alpha_0}(r-t,x), h(||\phi^{\alpha_0}(r-t,x))||)\,L^\ast D\delta(r,\phi^{\alpha_0}(r-t,x))\rangle
\end{align*}
is continuous on $[t,\,T-\varepsilon)$, uniformly in $\alpha_0$, and  on  $B_R(x):=\{x \in E: ||x||\leq R\}$, $R >0$. To this end, let $r, r' \in [t,\,T-\varepsilon)$. Since $\psi$ satisfies Definition \ref{D:testfunc}, in particular $\frac{\partial \psi}{\partial s}$, $L^\ast \, D \varphi$, $L^\ast \, D \delta$ are bounded on bounded sets of $E$. In the following $C$ will denote a generic constant that may depend on $T$ and $x$, and that may vary from line to line.  For any $x \in E: ||x|| \leq R$, $R>0$, we have 
\begin{align*}
&\left|\frac{\partial \psi}{\partial s} (r,\phi^{\alpha_0}(r-t,x))-\frac{\partial \psi}{\partial s} (r',\phi^{\alpha_0}(r'-t,x))\right|\\
&\leq C \omega(|r-r'| + ||\phi^{\alpha_0}(r-t,x))-\phi^{\alpha_0}(r'-t,x))||)\leq C\sigma_R( |r-r'|), 
\end{align*}
where in the latter inequality we have used  \eqref{contrflowestimate_s}. Using again the properties of the test functions, together with  \eqref{contrflowestimate_s}-\eqref{contrflowestimate_bound}, we get
\begin{align*}
	&|\langle \phi^{\alpha_0}(r-t,x),\,L^\ast \, D \varphi(r,\phi^{\alpha_0}(r-t,x)) \rangle-\langle \phi^{\alpha_0}(r'-t,x),\,L^\ast \, D \varphi(r',\phi^{\alpha_0}(r'-t,x)) \rangle|\\
	&\leq |\langle \phi^{\alpha_0}(r-t,x)- \phi^{\alpha_0}(r'-t,x),\,L^\ast \, D \varphi(r',\phi^{\alpha_0}(r'-t,x)) \rangle|\\
	&+ |\langle \phi^{\alpha_0}(r-t,x),\,L^\ast \, D \varphi(r,\phi^{\alpha_0}(r-t,x)) -L^\ast \, D \varphi(r',\phi^{\alpha_0}(r'-t,x)) \rangle|.\\
	&\leq 
	 C \sigma_R(|r-r'|)+ C(1 + ||x||) \, \omega(|r-r'|+ ||\phi^{\alpha_0}(r-t,x)-\phi^{\alpha_0}(r'-t,x)||)\\
	&\leq C\sigma_R(|r-r'|).
\end{align*}
Analogously, 
\begin{align*}
	&|\langle \phi^{\alpha_0}(r-t,x),h(||\phi^{\alpha_0}(r-t,x))||)\,L^\ast D\delta(r,\phi^{\alpha_0}(r-t,x))\rangle\\
	&- \langle \phi^{\alpha_0}(r'-t,x), h(||\phi^{\alpha_0}(r'-t,x))||)\,L^\ast D\delta(r',\phi^{\alpha_0}(r'-t,x))\rangle|\\
		&\leq C\sigma_R(|r-r'|).
\end{align*}

 Moreover, for any  $x\in E$,  $a \in A$, and any measurable function $\alpha_0 : \R_+ \times E \rightarrow A$, the map
 \begin{equation*}
r\mapsto f(\phi^{\alpha_0}(r-t,x),a)
 \end{equation*}
 is continuous on $[t,\,T-\varepsilon)$, uniformly in $\alpha_0$ and in $a$.
 Indeed, 
 from \textbf{(H$\textup{fg}$)} and (\ref{contrflowestimate-BIS}), 
for any $r \in [t,\,T-\varepsilon)$,  $x \in E: ||x|| \leq R$, $R>0$, 
\begin{align*}
&|f(\phi^{\alpha_0}(r'-t,x),a)-f(\phi^{\alpha_0}(r-t,x),a)| \\
&\leq   C \omega(||\phi^{\alpha_0}(r-t,x),a)-\phi^{\alpha_0}(r'-t,x),a)||_{-1})\leq C   \sigma_R(|r-r'|).
 \end{align*}

Let us finally study the continuity of the maps 
\begin{align*}
&r \mapsto \mathcal L^a \psi(r,\phi^{\alpha_0}(r-t,x),a))\\
&= \langle b(\phi^{\alpha_0}(r-t,x),a), D \psi(r,\phi^{\alpha_0}(r-t,x))\rangle  \\
&+\lambda(\phi^{\alpha_0}(r-t,x), a)\int_{E} (\psi(r,y)-\psi(r,\phi^{\alpha_0}(r-t,x)))\,   Q(\phi^{\alpha_0}(r-t,x),a, dy).
\end{align*}
Since  by Definition \ref{D:testfunc} $D \psi$ is bounded on bounded sets of $E$, and using assumption \textbf{(H$\textup{b$\lambda$Q}$)} for $b$,  for any $x \in E: ||x|| \leq R$, $R>0$ we get 
\begin{align*}
&|\langle b(\phi^{\alpha_0}(r-t,x),a), D \psi(t,\phi^{\alpha_0}(r-t,x))\rangle-\langle b(\phi^{\alpha_0}(r'-t,x),a), D \psi(r,\phi^{\alpha_0}(r'-t,x))\rangle|	\\
&\leq |\langle b(\phi^{\alpha_0}(r-t,x),a)-\langle b(\phi^{\alpha_0}(r'-t,x),a), D \psi(r,\phi^{\alpha_0}(r'-t,x))\rangle|\\
&+
|\langle b(\phi^{\alpha_0}(r-t,x),a), D \psi(t,\phi^{\alpha_0}(r-t,x))\rangle-D \psi(r,\phi^{\alpha_0}(r'-t,x))\rangle|\\
&\leq C\sigma_R
(|r-r'|)+ C\omega (|r-r'| + ||\phi^{\alpha_0}(r-t,x),a)-\phi^{\alpha_0}(r'-t,x)||)\\
&\leq C\,\sigma_R 
(|r-r'|).
\end{align*}
 On the other hand, by assumption \textbf{(H$\textup{b$\lambda$Q}$)} for $\lambda$ and $Q$, recalling that $\psi$ is uniformly continuous on $(\varepsilon, T-\varepsilon) \times E$,  for any $x \in E: ||x|| \leq R$, $R>0$ we have  
 \begin{align*}
 	&\Big|\lambda(\phi^{\alpha_0}(r-t,x), a)\int_{E} (\psi(r,y)-\psi(r,\phi^{\alpha_0}(r-t,x)))\,   Q(\phi^{\alpha_0}(r-t,x),a, dy)\\
 	&-\lambda(\phi^{\alpha_0}(r'-t,x), a)\int_{E} (\psi(r',y)-\psi(r',\phi^{\alpha_0}(r'-t,x)))\,   Q(\phi^{\alpha_0}(r'-t,x),a, dy)\Big|\\
 	&\leq |\lambda(\phi^{\alpha_0}(r-t,x), a)
 	-\lambda(\phi^{\alpha_0}(r'-t,x), a)| \Big|\int_{E} (\psi(r',y)-\psi(r',\phi^{\alpha_0}(r'-t,x)))\,   Q(\phi^{\alpha_0}(r'-t,x),a, dy)\Big|\\
 	&+|\lambda(\phi^{\alpha_0}(r-t,x), a)| \,
 	|\psi(r',\phi^{\alpha_0}(r'-t,x)) -\psi(r,\phi^{\alpha_0}(r-t,x))|\\
 	& +|\lambda(\phi^{\alpha_0}(r-t,x), a)| \, \Big|\int_{E} \psi(r,y)\,   Q(\phi^{\alpha_0}(r-t,x),a, dy)-\int_{E} \psi(r',y)\,   Q(\phi^{\alpha_0}(r'-t,x),a, dy)\Big|\\
 	&\leq C ||\phi^{\alpha_0}(r-t,x)-\phi^{\alpha_0}(r'-t,x)||_{-1}	+ C \omega(|r-r'| + ||\phi^{\alpha_0}(r'-t,x)-\phi^{\alpha_0}(r-t,x)||)\\
 	& + C \sigma_R
 	(|r-r'|) + C \omega(||\phi^{\alpha_0}(r-t,x), a)-\phi^{\alpha_0}(r'-t,x), a)||_{-1})\\
 	&\leq C \sigma_R
 	(|r-r'|),
\end{align*}
where the latter inequality follows from 
\eqref{contrflowestimate_s}-\eqref{contrflowestimate_bound}-\eqref{contrflowestimate-1}.
\qed
 
 \subsection{Proof of Theorem \ref{Sec:PDP_Thm_existence}}\label{Sec_visc_prop_V}
 
We start by giving the following preliminary result. 

\begin{lemma}\label{L:ineqIto1}
 Assume that Hypotheses \textup{\textbf{(HL)}} and \textup{\textbf{(H$\textup{b$\lambda$Q}$)}} hold.
Let $0 < t < \bar T < T$, $\hat \tau$ be a stopping time such that $\hat \tau  \in [t,\,\bar T]$, $x \in E$, $\alpha  \in \mathcal A_{ad}^t$, and $X$  be the process  in \eqref{Sec:PDP_controlledX} under $\P^{t,x}_{\alpha}$. For $R>0$, let $\tau_R$ be the exit time of $X$ from $\{y:\,||y||\leq R\}$, and set $\tau= \hat \tau \wedge \tau_R$. Let $\psi = \varphi + h(||\cdot||)\,\delta$ be a test function. Then, 
\begin{align}\label{ineq_ito}
	\E^{t,x}_{\alpha}\left[\psi(\tau,X_\tau)\right] &\leq \psi(t,x)+ \E^{t,x}_{\alpha}\left[\int_t^\tau \left(\frac{\partial \psi}{\partial t}(r,X_r) + \langle b(X_r, \alpha_r),\,D\psi(r,X_r)\rangle\right)dr\right]\notag\\
	&-\E^{t,x}_{\alpha}\left[\int_t^\tau \langle X_r,\,L^\ast D\varphi(r,X_r) + h(||X_r||)\,L^\ast D\delta(r,X_r)\rangle \,dr\right]\notag\\
	&+\E^{t,x}_{\alpha}\left[\int_t^\tau\int_E (\psi(r,y)-\psi(r,X_r))\,\lambda(X_r, \alpha_r)\,Q(X_r, \alpha_r,\,dy) \,dr\right].
\end{align}
\end{lemma}
\noindent \emph{Proof of Lemma \ref{L:ineqIto1}.}
The result follows from the Dynkin formula \eqref{itoformula} and the properties of the test functions $\psi$ in Definition \ref{D:testfunc}. In particular, 
$D\psi(r,X_r) = D \varphi(r, X_r) + h(||X_r||)\,D\delta(r, X_r) + \delta(r, X_r) \, \frac{h'(||X_r||)}{||X_r||} X_r$, 
and $\langle L X_r, \delta(r, X_r) \, \frac{h'(||X_r||)}{||X_r||} X_r\rangle \geq 0$, being $L$  monotone. 
 \endproof

  \noindent \paragraph{Proof of Theorem \ref{Sec:PDP_Thm_existence}.} 
  \textbf{Viscosity subsolution property.}
 Let $\psi(s,y) = \varphi(s,y) + \delta(s,y) \,h(||y||)$ be a test function of the type introduced in Definition \ref{D:testfunc}, such that  $V- \psi$ has a global maximum at $(t,x) \in ]0,\,T[ \times E$. We also assume that 
\begin{equation}\label{v+psi=0BIS}
V(t,x)= \psi(t,x), 
\end{equation}
and consequently
\begin{equation}\label{ineqvpsi1}
V(s,y) \leq \psi(s,y), \quad \forall \,(s,y).
\end{equation}
Remember that $T_1$ denotes the first jump time of $X$. 
From the dynamic programming principle \eqref{DPP2}  applied to $\theta :=(t+\eta) \wedge T_1$ where $\eta>0$ is such that $(t+\eta)<T$ otherwise arbitrary for the moment, 
and since  by (\ref{ineqvpsi1}),
$V(s,y) \leq \psi(s,y)$ for all $(s,y)$, we have
 \begin{eqnarray}\label{byDPP}
 	\psi(t, x)  \leqslant \E^{t, x}_\alpha\left[\psi(\theta, X_{\theta}) + \int_{t}^{\theta} f(X_{r}, \alpha_r)\,dr\right],\quad \forall \alpha \in \mathcal  A_{ad}^t.
 \end{eqnarray}
All elements of $\mathcal  A_{ad}^t$ have the form (\ref{Sec:PDP_open_loop_controls}). Let us fix $a\in A$, and let us take  $\alpha \in \mathcal  A_{ad}^t$ such that $\alpha_0\equiv a$. Notice that, $\P^{t, x}_{\alpha}$-a.s., $X_r = \phi^{a}(r-t,x)$ for $r \in [t,\,\theta)$. In particular, by \eqref{contrflowestimate_bound}, 
$$
||X_s|| \leq C (1 + ||x||)=:R_x. 
$$
Denoting by $\tau_R$  the exit time of $X$ from $\{y:\,||y||\leq R\}$, it follows that $\theta=(t+h) \wedge T_1 \wedge \tau_{R_x}$. As a matter of fact, if $(t+h) \leq T_1$, then $(t+h) \wedge T_1 \wedge \tau_{R_x}=t+h$. On the other hand, if $(t+h) > T_1$, we have two cases: if $X_{T_1} \notin B_{R_x}$, then $(t+h) \wedge T_1 \wedge \tau_{R_x}= \tau_{R_x}= T_1$,  if $X_{T_1} \in B_{R_x}$, then $(t+h) \wedge T_1 \wedge \tau_{R_x}= T_1$.
 Then (\ref{byDPP}) for such an $\alpha$ and Lemma \ref{L:ineqIto1} imply that 
\begin{align}\label{Sec:PDP_sub_sol_int_ineq_a}
 &\E^{t, x}_{\alpha}\left[\int_{t}^{\theta} \left[\frac{\partial \psi}{\partial t}(r,X_r)  + \mathcal L^{a}\psi(r,X_r)+f(X_{r}, a)
\right]dr\right]\notag\\
	&-\E^{t, x}_{\alpha}\left[\int_{t}^{\theta} \langle X_r,\,L^\ast D\psi(r,X_r) + h(||X_r||)\,L^\ast D\delta(r,X_r)\rangle \,dr\right] \geq 0,
\end{align}
with $X_r = \phi^{a}(r-t,x)$. 
Moreover, by Lemma \ref{Gpsi}, the (deterministic) map  
\begin{align*}
 r\mapsto &\frac{\partial \psi}{\partial t}(r,\phi^a(r-t,x))+\mathcal L^{a}\psi(r,\phi^a(r-t,x))+f(\phi^a(r-t,x),a)\\
 &-\langle \phi^a(r-t,x),\,L^\ast D\psi(r,\phi^a(r-t,x)) + h(||\phi^a(r-t,x)||)\,L^\ast D\delta(r,\phi^a(r-t,x))\rangle
\end{align*}
 is continuous  at $t$,
 uniformly in $a$. 
Therefore, for any $\epsilon>0$ there exists $\eta>0$ such that (\ref{Sec:PDP_sub_sol_int_ineq_a}) with $\theta$ associated to $\eta$ becomes
 \begin{equation}\label{sign1}
 \left(\varepsilon + \frac{\partial \psi}{\partial t}(t,x) +\mathcal L^{a} \psi(t, x)+f(x,a)-\langle x,\,L^\ast D\psi(t,x) + h(||x||)\,L^\ast D\delta(t,x)\rangle\right) \E^{t, x}_{\alpha}\left[\theta-t\right]\geq 0,
 \end{equation}
 valid for any $a\in A$.
 Now we observe that $\E^{t, x}_{\alpha}\left[\theta-t\right]\geq 0$ by definition of $\theta$. Then (\ref{sign1}) implies 
  \begin{equation}\label{sign2}
 \left(\varepsilon + \frac{\partial \psi}{\partial t}(t,x) +\mathcal L^{a} \psi(t, x)+f(x,a)-\langle x,\,L^\ast D\psi(t,x) + h(||x||)\,L^\ast D\delta(t,x)\rangle\right)\geq 0,
 \end{equation}
 for any $\varepsilon>0$ and $a\in A$. The conclusion follows by the arbitrariness of $\varepsilon$ and $a$.

\medskip
 
\medskip
 
 \noindent \textbf{Viscosity supersolution property.}
Let $\psi(s,y) = \varphi(s,y) + \delta(s,y) \,h(||x||)$ be a test function of the type introduced in  Definition \ref{D:testfunc}, such that  $V+ \psi$  has a global minimum at $(t,x) \in ]0,\,T[\times E$.  We also assume that 
\begin{equation}\label{v+psi=02}
V(t,x)+ \psi(t,x) = 0, 
\end{equation}
so  
we have 
\begin{equation}\label{ineqvpsi2}
V(s,y) \geq -\psi(s,y), \quad \forall \,(s,y).
\end{equation}

We will show that $V$ is a viscosity supersolution by contradiction. Let us use the notations of Lemma \ref{Gpsi}.  Assume that 
\begin{equation}\label{Gmu}
\inf_{a \in A} G^\psi_a (t,x)= \mu >0.
\end{equation}
By Lemma \ref{Gpsi},
there exists $\eta>0$, independent from $\alpha_0$,  such that 
\begin{equation}\label{cont G alpha}
{\cal G}^{\alpha_0}(r)\geq \frac{\mu}{2}>0,  \quad \forall r\in [t,t+\eta), 
\end{equation}
Let us now set $\theta:=(t+\eta)\wedge T_1$ where $\eta$ satisfies $(t+\eta)<T$.  
By the dynamic programming principle  \eqref{DPP2},  for all $\gamma>0$ there exists $\alpha \in \mathcal A_{ad}^t$ such that 
\begin{equation*}
V(t,x)+\gamma \geq \E^{t,x}_{\alpha}\left[\int_{t}^{\theta}  f(X_{r},\alpha_r)\,dr + V(\theta,X_{\theta})\right],
\end{equation*}
and therefore, recalling \eqref{v+psi=02} and \eqref{ineqvpsi2}, 
\begin{equation} \label{Sec:PDP_ineq_num1-gamma}
-\psi(t,x) +\gamma \geq \E^{t,x}_{\alpha}\left[\int_{t}^{\theta}  f(X_{r},\alpha_r)\,dr -\psi(\theta,X_{\theta})\right].
\end{equation}
As in the proof of the viscosity subsolution property, we set $R_x$ to the the bound in \eqref{contrflowestimate_bound}, and we notice  that $\theta=(t+h) \wedge T_1 \wedge \tau_{R_x}$, where $\tau_R$  denotes  the exit time of $X$ from $\{y:\,||y||\leq R\}$. 
Applying Lemma \ref{L:ineqIto1} to $\psi$ between $t$ and $\theta$, we get  
\begin{align*}
\gamma &\geqslant \E^{t, x}_{\alpha}\left[\int_{t}^{\theta} f(X_{r}, \alpha_r)\,dr \right] 
 -\E^{t, x}_{\alpha}\left[\int_{t}^{\theta} \left[\frac{\partial \psi}{\partial t}(r,X_r)  +\langle b(X_r, \alpha_r),\,D\psi(X_r, \alpha_r)\rangle\right]dr\right]\notag\\
	&-\E^{t, x}_{\alpha}\left[\int_{t}^{\theta}\int_E (\psi(r,y)-\psi(r,X_r))\,\lambda(X_r, \alpha_r)\,Q(X_r, \alpha_r,\,dy) \,dr\right]\notag\\
	&+\E^{t, x}_{\alpha}\left[\int_{t}^{\theta} \langle X_r,\,L^\ast D\varphi(r,X_r) + h(||X_r||)\,L^\ast D\delta(r,X_r)\rangle \,dr\right].
\end{align*}
Then 
\begin{align}\label{Sec:PDP_ineq_num_BIS1-2}
\gamma
& \geqslant  \E^{t, x}_{\alpha}\left[\int_{t}^{\theta} \left(-\frac{\partial \psi}{\partial t}(r,X_r)+\langle X_r,\,L^\ast D\varphi(r,X_r) + h(||X_r||)\,L^\ast D\delta(r,X_r)\rangle \right)dr\right]\notag\\
	&+\E^{t, x}_{\alpha}\left[\int_{t}^{\theta}  \inf_{a \in A} (-{\cal L}^a\psi(r,X_r)+f(X_r,a))\,dr \right]\notag\\
	&=\E^{t, x}_{\alpha}\left[\int_{t}^{\theta}  \, \inf_{a \in A} G_a^\psi(r,X_r)\, dr\right].
\end{align}
By the definition of $\theta$, together with \eqref{Sec:PDP_open_loop_controls}  and \eqref{Sec:PDP_controlledX},  for all $r\in [t,\theta)$, $\alpha_r=\alpha_0(r-t,x)$, with $\alpha_0$ as in (\ref{Sec:PDP_open_loop_controls}) 
and $X_r=\phi^{\alpha_0}(r-t,x)$. Thus (\ref{cont G alpha}) yields 
\begin{equation}
\gamma\geq \E^{t, x}_{\alpha}\left[\int_{t}^{\theta}  \, {\cal G}^{\alpha_0}(r)\, dr\right]\geq \frac{\mu}{2}\, \E^{t, x}_{\alpha}\left[(\theta-t)\right].
\end{equation}
Now we notice that
\begin{eqnarray*}
 \E^{t, x}_{\alpha}(\theta-t)&=& \eta \E^{t, x}_{\alpha}(\one_{T_1>t+\eta})+\E^{t, x}_{\alpha}((T_1-t)\one_{T_1\leq t+\eta})\notag\\
&\geq& \eta \, \P^{t, x}_{\alpha}(\one_{T_1>t+\eta})\\
&=&\eta \, \, e^{-\int_{t}^{t+\eta}\lambda(\phi^{\alpha_0}(s,x), \alpha_0(s,x))\,dr}\\
&\geq&  \eta\,\, e^{-\eta\, M},
\end{eqnarray*} 
 where in the latter inequality we have used that by assumption  $\lambda$ is bounded by some constant $M$. 
By letting $\gamma$ go to zero we obtain the contradiction.

\section{Proofs of the  results in Section \ref{Sec:PDP_Sec_ConstrainedBSDE}}\label{Sec_mainproofs_Sec4}

\subsection{Proof of Theorem \ref{Sec:PDP_THm_Feynman_Kac_HJB}}\label{Sec:proofThm5.1}

The boundedness of $v$ follows from \eqref{Sec:PDP_ident_vdelta},  \eqref{Sec:PDP_Vstar_Y0},  together with the  definition of $V^\ast$ in \eqref{Sec:PDP_dual_value_function} and the assumption \textup{\textbf{(H$\textup{fg}$)}}.

Let us now turn to the continuity properties.
We argue as in the proof of Proposition \ref{P:dynprogpr}.   
We denote by  $B([0,\,T] \times E)$  the set of all bounded functions on $[0,\,T] \times E$,   and, for any $a \in A$,  we define the map $\mathcal T_a: B([0,\,T] \times E) \rightarrow B([0,\,T] \times E)$ as 
\begin{align*}
\mathcal T_a \psi(t,x) &:=  \inf_{\nu \in \mathcal V}  \E^{t,x,a}_{\nu}\left[\int_{t}^{T_1 \wedge T} f(X_s,I_s)\, ds  +   g(X_T)\one_{T \leq T_1} +\psi(T_1, X_{T_1})\one_{T >T_1}\right].
\end{align*}
We have 
\begin{align*}
\mathcal T_a \psi(t,x) & =\inf_{\nu \in \mathcal V} \bigg\{\int_0^{T-t}\chi^{\nu}(s, x,a) (f(s,  x,a)+L_{\psi}(s,x,a))\, ds  + \chi^{\nu}(T- t, x,a) g(\phi(T-t,x,a)\bigg\},
\end{align*}
where $\chi^{\nu}(s, x,a)= e^{-\int_0^{s}  \left(\lambda(\phi(r, x,a), a)+ \int_A \nu_s(b)\lambda_0(db)\right)\,dr}$, $f(s,  x,a) = f(\phi(s, x,a), a)$,  and 
\begin{align*}
L_{\psi}(s,  x,a)&=  \int_E \psi(s,y)\,\lambda(\phi(s, x,a), a)\,Q(\phi(s, x,a), a,dy).
\end{align*} 
As in the case of the map in \eqref{new_G}, for any $a \in A$,  $\mathcal T_a$ is a contracting map in $B([0,\,T] \times E)$ and $v$ is its unique fixed point. 
In particular,   $v$ 
satisfies the randomized DPP \eqref{RandDynProgPr}.

Then, we denote by $C_b([0, T] \times E)$ the set of bounded functions, continuous on $[0,\,T] \times E$  with the $|\cdot| \times ||\cdot||_{-1}$ norm. 
We aim at showing that, for any function  $\psi  \in C_b([0, T] \times E)$, for any $a \in A$ one has $\mathcal T_a \psi \in C_b([0, T] \times E)$.  This would prove that  $v \in C_b([0, T] \times E)$.  

In what follows $C$ will denote a generic constant, that may vary from line to line, and that may depend  on $T$.
Let $t,t',s \in [0,\,T]$, $t'\leq t \leq s$, $x,x' \in E$, $a \in A$, $\nu \in \mathcal V$. Recalling hypotheses \textbf{(H$\textup{b$\lambda$Q}$)}-(i),  \textup{\textbf{(H$\textup{fg}$)}} and \eqref{flowestimateRand}, we have $|\chi^{\nu}(s, x,a)| \leq 1$, $|f(s, x,a)| \leq C$, and, for any $s' \leq s$,  
\begin{align}
&|\chi^\nu(s',x,a)-\chi^\nu(s,x',a) |\leq 
(1 - e^{-C ||x-x'||_{-1}})+(1 - e^{-C (s-s')}),\label{chiest2}\\
&|f(s, x,a)-f(s, x',a)| \leq C ||x-x'||_{-1}\label{fest2},\\
&|g(\phi(T-t, x,a))- g(\phi(T-t', x',a))|\leq \omega(||x-x'||_{-1})\label{gest12}.
\end{align}
On the other hand,  by  \textup{\textbf{(H$\textup{b$\lambda$Q}$)}}-(i)-(ii), together with the boundedness and continuity of $\psi$,  we have  $|L_{\psi}(s, x,a)\leq C$ and, for $s < T-t$, 
\begin{align}
&|L_{\psi}(s, x,a)-L_{\psi}(s, x',a)|\leq |\lambda(\phi(s, x,a),a)-\lambda(\phi(s, x',a),a)|\,||\psi||_{\infty}\notag\\
&+||\lambda||_{\infty}\bigg|\int_E \psi(s,y)\,[Q(\phi(s, x,a),a,dy)-Q(\phi(s, x',a),a,dy)]\bigg|\notag\\
&\leq C \,\sigma(||\phi(s,x,a)-\phi(s,x',a)||_{-1})\leq C \omega(||x-x'||_{-1}),\label{Lest12}
\end{align}
where the latter inequality follows from \eqref{flowestimateRand}.
Then, for any $t, t' \in [0,\,T]$, $x, x' \in E$, $a \in A$, $\nu \in \mathcal V$, 
\begin{align*}
&|J(t,x, a, \nu)-J(t',x', a,\nu)|\notag\\
&\leq\left|\int_0^{T-t}\chi^\nu(s,x,a) f(s,x,a)\,ds-\int_{0}^{T-t'}\chi^\nu(s,x',a) f(s,x',a)\,ds\right|\notag\\
 	&+ \left|\int_0^{T-t} \chi^\nu(s,x,a) L_\psi(s,x,a)\,ds-\int_{0}^{T-t'} \chi^u(s,x') L_\psi(s,x',a)\,dr\right|\\
 	& + |\chi^{\nu}(T-t, x) g(\phi(T-t, x,a)) - \chi^{\nu}(T-t', x') g(\phi(T-t', x',a)) |\notag\\
&\leq \int_{0}^{T-t}|\chi^\nu(s,x,a) f(s,x,a)- \chi^\nu(s,x',a) f(s,x',a)|\,ds\notag\\	
&+ \int_{0}^{T-t} |\chi^\nu(s,x,a) L_\psi(s,x)- \chi^\nu(s,x',a) L_\psi(s,x',a)|\,ds +  C |t-t'|\notag\\
 	& + C|g(\phi(T-t, x,a)) -  g(\phi(T-t', x',a)) |+ C|\chi^{\nu}(T-t, x,a)  - \chi^{u}(T-t', x') |\notag\\
&\leq C \Bigg(\int_{0}^{T-t}|\chi^\nu(s,x,a))-\chi^\nu(s,x',a)|\,ds+\int_{0}^{T-t} |f(s,x,a)-f(s,x',a)|\,ds\notag\\	
&+ \int_{0}^{T-t} | L_\psi(s,x,a)  -  L_\psi(s,x',a)| \,ds  + |g(\phi(T-t, x,a)) -  g(\phi(T-t', x',a)) |\\
&+ |\chi^{\nu}(T-t, x,a)  - \chi^{\nu}(T-t', x',a) |+   |t-t'|\Bigg)\\
&\leq C (\omega(t-t')+ \omega'(||x-x'||_{-1}))
\end{align*}
for some modulus of continuity $\omega$, $\omega'$, where the latter inequality follows from  \eqref{chiest2}, \eqref{fest2}, \eqref{gest12}, \eqref{Lest12}.
This  shows in particular   that $v$ is uniformly continuous  in the $|\cdot| \times ||\cdot||_{-1}$ norm.

 \subsection{Proof of Theorem \ref{Sec:first_main result}}\label{Sec:proof_viscprop_rand}

We first give the following preliminary result. 
\begin{lemma}\label{L:ineqIto}
Let $0 < t < \bar T < T$, $\hat \tau$ be a stopping time such that $\hat \tau  \in [t,\,\bar T]$, $x \in E$, $a \in A$, $\nu \in \mathcal V$, and $(X,I)$ be the PDMP constructed in Section \ref{Sec:PDP_Section_control_rand} under the probability $\P^{t,x,a}_\nu$. For $R>0$, let $\tau_R$ be the exit time of $X$ from $\{y:\,||y||\leq R\}$, and set $\tau= \hat \tau \wedge \tau_R$. Let $\psi = \varphi + h(||\cdot||)\,\delta$ be a test function. Then, 
\begin{align}\label{ineq_ito}
	\spernutxa{\psi(\tau,X_\tau)} &\leq \psi(t,x)+ \spernutxa{\int_t^\tau \left(\frac{\partial \psi}{\partial t}(r,X_r) + \langle b(X_r, I_r),\,D\psi(r,X_r)\rangle\right)dr}\notag\\
	&-\spernutxa{\int_t^\tau \langle X_r,\,L^\ast D\varphi(r,X_r) + h(||X_r||)\,L^\ast D\delta(r,X_r)\rangle \,dr}\notag\\
	&+\spernutxa{\int_t^\tau\int_E (\psi(r,y)-\psi(r,X_r))\,\lambda(X_r, I_r)\,Q(X_r, I_r,\,dy) \,dr}.
\end{align}
\end{lemma}
\noindent \emph{Proof of Lemma \ref{L:ineqIto}.}
By Proposition \ref{T_ItoFormula}, applying the Dynkin formula to $\psi(s, X_s)$ between $t$ and $\tau$ and taking the expectation under $\P^{t,x,a}_\nu$, we get 
\begin{align*}
&\spernutxa{\psi(\tau,X_\tau)} = \psi(t,x)+ \spernutxa{\int_t^\tau \left(\frac{\partial \psi}{\partial t}(r,X_r) + \langle b(X_r, I_r),\,D\psi(r,X_r)\rangle\right)dr}\notag\\
	&+ \spernutxa{\int_t^\tau  \bigg(-\langle L\,X_r, \,D\psi(r,X_r)\rangle 
	+\int_E (\psi(r,y)-\psi(r,X_r))\,\lambda(X_r, I_r)\,Q(X_r, I_r,\,dy)\bigg) \,dr}.
\end{align*}
We conclude noticing that 
$
D\psi(r,X_r) = D \varphi(r, X_r) + h(||X_r||)\,D\delta(r, X_r) + \delta(r, X_r) \, \frac{h'(||X_r||)}{||X_r||} X_r$, and 
that $\langle L X_r, \delta(r, X_r) \, \frac{h'(||X_r||)}{||X_r||} X_r\rangle \geq 0$, being $L$ is monotone.

\noindent \paragraph{Proof of Theorem \ref{Sec:first_main result}.}
\noindent \textbf{Viscosity subsolution property.}
 Let $\psi(s,y) = \varphi(s,y) + \delta(s,y) \,h(||x||)$ be a test function of the type introduced in  Definition \ref{D:testfunc}, such that  $v- \psi$  has a global maximum at $(t,x) \in [0,\,T]\times E$. We also assume that 
\begin{equation}\label{v+psi=0BIS}
v(t,x)= \psi(t,x), 
\end{equation}
so  
we have 
\begin{equation}\label{ineqvpsi2BIS}
v(s,y) \leq \psi(s,y), \quad \forall \,(s,y).
\end{equation}
Fix $(t,x,a)$ and $\nu \in \mathcal  V$. Let $\eta>0$ and define $\theta =(t+\eta) \wedge T_1$,  where $T_1$ denotes the first jump time of $(X,I)$. 
Using the identification property \eqref{Sec:PDP_ident_vdelta}, from the randomized dynamic programming principle \eqref{RandDynProgPr}, together with \eqref{ineqvpsi2BIS},  we get 
 \begin{eqnarray*}
 	\psi(t, x)  \leqslant \E^{t, x,a}_\nu\left[\psi(\theta, X_{\theta}) + \int_{t}^{\theta} f(X_{r}, I_r)\,dr\right].
 \end{eqnarray*}
 Applying Lemma \ref{L:ineqIto}, we obtain 
\begin{align}\label{Sec:PDP_sub_sol_int_ineq}
 &\E^{t, x,a}_{\nu}\left[\int_{t}^{\theta} \left[\frac{\partial \psi}{\partial t}(r,X_r)  + \mathcal L^{I_r}\psi(r,X_r)+f(X_{r}, I_r)
\right]dr\right]\notag\\
	&-\E^{t, x,a}_{\nu}\left[\int_{t}^{\theta} \langle X_r,\,L^\ast D\psi(r,X_r) + h(||X_r||)\,L^\ast D\delta(r,X_r)\rangle \,dr\right] \geq 0,
\end{align} 
 where 
\begin{align}\label{LI}
\mathcal L^{I_r}\psi(r,X_r)= \langle b(X_r, I_r),\,D\psi(r, X_r)\rangle + \int_E (\psi(r,y)-\psi(r,X_r))\,\lambda(X_r, I_r)\,Q(X_r, I_r,\,dy).
\end{align}
 Now we notice that $\P^{t, x,a}$-a.s., for all $r\in(t,\theta)$, $(X_r,I_r) = (\phi(r-t,x,a),a)$.  
 Moreover, by Lemma \ref{Gpsi}, the map 
  \begin{align*}
r\mapsto &
\frac{\partial \psi}{\partial t}(r,\phi(r-t,x,a))+\mathcal L^{a}\psi(r,\phi(r-t,x,a))+f(\phi(r-t,x,a),a)\\
&-\langle \phi(r-t,x,a)\,L^\ast D\psi(r,\phi(r-t,x,a)) + h(||\phi(r-t,x,a)||)\,L^\ast D\delta(r,\phi(r-t,x,a))\rangle
 \end{align*}
 is continuous,  uniformly with respect to  $a\in A$. We can proceed as  in the proof of Theorem \ref{Sec:PDP_Thm_existence}. By the latter continuity property, for any $\epsilon>0$, we can find $\eta>0$ independent of $a$ such that  (\ref{Sec:PDP_sub_sol_int_ineq}) holds true for $\theta$ corresponding to $\eta$. Since $\E^{t, x}_{\alpha}\left[\theta-t\right]\geq 0$ by definition of $\theta$, then identity (\ref{Sec:PDP_sub_sol_int_ineq}) implies 
  \begin{equation}\label{sign2}
 \left(\varepsilon + \frac{\partial \psi}{\partial t}(t,x) +\mathcal L^{a} \psi(t, x)+f(x,a)-\langle x,\,L^\ast D\psi(t,x) + h(||x||)\,L^\ast D\delta(t,x)\rangle\right)\geq 0,
 \end{equation}
 for any $\epsilon>0$ and $a\in A$.
  As in the proof of Theorem \ref{Sec:PDP_Thm_existence}, we conclude by the arbitrariness of $\varepsilon$ and $a$.

\medskip

 \noindent \textbf{Viscosity supersolution property.}
Let $\psi(s,y) = \varphi(s,y) + \delta(s,y) \,h(||x||)$ be a test function of the type introduced in  Definition \ref{D:testfunc}, such that  $v+ \psi$  has a global minimum at $(t,x) \in [0,\,T]\times E$. We also assume that 
\begin{equation}\label{v+psi=0}
v(t,x)+ \psi(t,x) = 0, 
\end{equation}
so  
we have 
\begin{equation}\label{ineqvpsi}
v(s,y) \geq -\psi(s,y), \quad \forall \,(s,y).
\end{equation}
 We will show that $v$ is a viscosity supersolution by contradiction.  Let us use the notations of Lemma \ref{Gpsi}.  Assume that
\begin{equation}\label{Gmu}
G^\psi\left(t,x,\psi, D\varphi, D\delta\right)= \mu >0.
\end{equation}
By Lemma \ref{Gpsi} that we apply for $\alpha_0\equiv a$, $a\in A$ arbitrary, 
there exists $\eta>0$, independent from $a$,  such that 
\begin{equation}\label{cont G a}
{\cal G}^{a}(r)\geq \frac{\mu}{2}>0 \quad \forall r\in [t,t+\eta). 
\end{equation}
 Let us set $\theta=(t+\eta)\wedge T_1$ and  fix $a \in A$.
By the dynamic programming principle  \eqref{RandDynProgPr} toghether with  the identification property \eqref{Sec:PDP_ident_vdelta}, 
we see that, for all $\gamma >0$,  it exists a strictly  positive, predictable and bounded function $\nu$ such that
\begin{eqnarray*}
v(t,x) + \gamma \geqslant  \E^{t, x,a}_{\nu}\left[\int_{t}^{\theta} f(X_{r}, I_r)\,dr + v(\theta, X_{\theta}) \right]. 
\end{eqnarray*}
Recalling \eqref{v+psi=0} and \eqref{ineqvpsi}, we get
\begin{eqnarray}\label{Sec:PDP_ineq_num}
-\psi(t,x) + \gamma \geqslant  \E^{t, x,a}_{\nu}\left[\int_{t}^{\theta} f(X_{r}, I_r)\,dr - \psi(\theta, X_{\theta}) + \beta(\eta)\,\one_{\tau  \wedge T_1\leq T}\right]. 
\end{eqnarray}
Applying Lemma \ref{L:ineqIto}, inequality \eqref{Sec:PDP_ineq_num} yields  
\begin{align}\label{Sec:PDP_ineq_num_BIS}
\gamma
& \geqslant \E^{t, x,a}_{\nu}\left[\int_{t}^{\theta} f(X_{r}, I_r)\,dr 
-\int_{t}^{\theta} \left(\frac{\partial \psi}{\partial t}(r,X_r)  +\langle b(X_r, I_r),\,D\psi(X_r, I_r)\rangle\right)dr\right]\notag\\
	&-\E^{t, x,a}_{\nu}\left[\int_{t}^{\theta}\int_E (\psi(r,y)-\psi(r,X_r))\,\lambda(X_r, I_r)\,Q(X_r, I_r,\,dy) \,dr\right]\notag\\
	&+\E^{t, x,a}_{\nu}\left[\int_{t}^{\theta} \langle X_r,\,L^\ast D\psi(r,X_r) + h(||X_r||)\,L^\ast D\delta(r,X_r)\rangle \,dr\right].
\end{align}
Noticing that
\begin{align*}
	-\mathcal L^{I_r} \psi(r,X_{r})+f(X_{r}, I_r) = U^{\psi}(r,X_r, I_r, D \psi)\geqslant 
	\inf_{a \in A}U^\psi(r,X_r,a, D\psi),
\end{align*}
with $\mathcal L^I$ is the operator in \eqref{LI},  previous inequality
 gives
\begin{align}\label{Sec:PDP_ineq_num_BIS1-2}
\gamma
& \geqslant  \E^{t, x,a}_{\nu}\left[\int_{t}^{\theta} \left(-\frac{\partial \psi}{\partial t}(r,X_r)+\langle X_r,\,L^\ast D\varphi(r,X_r) + h(||X_r||)\,L^\ast D\delta(r,X_r)\rangle \right)dr\right]\notag\\
	&+\E^{t, x,a}_{\nu}\left[\int_{t}^{\theta}  \inf_{a \in A} U^\psi(r,X_r,a, D\psi)\,dr \right]\notag\\
	&=\E^{t, x,a}_{\nu}\left[\int_{t}^{\theta}  \, G^\psi\left(r,X_r,\psi, D\psi, D\varphi, D\delta\right)\, dr\right].
\end{align}
By the definition of $\theta$, together with \eqref{Sec:PDP_XI_def},
  for all $r\in [t,\theta)$, 
$X_r=\phi(r-t,x,a)$.
Thus,  \eqref{Sec:PDP_ineq_num_BIS1-2} together with \eqref{cont G a} yields 
\begin{equation}
\gamma\geq \E^{t, x,a}_{\nu}\left[\int_{t}^{\theta}  \, G^\psi\left(r,\phi^{\alpha_0}(r-t,x),\psi, D\psi, D\varphi, D\delta\right)\, dr\right]\geq \frac{\mu}{2}\, \E^{t, x,a}_{\nu}\left[(\theta-t)\right].
\end{equation}
We conclude as in the proof of Theorem  \ref{Sec:PDP_Thm_existence} using that 
\begin{eqnarray*}\label{ineqfinal3}
 \E^{t, x}_{\alpha}(\theta-t)&=& \eta \E^{t, x}_{\alpha}(\one_{T_1>t+\eta})+\E^{t, x}_{\alpha}((T_1-t)\one_{T_1\leq t+\eta})\notag\\
&\geq& \eta \, \P^{t, x}_{\alpha}(\one_{T_1>t+\eta})\\
&=&\eta \, \, e^{-\int_{t}^{t+\eta}\lambda(\phi^{\alpha}(s,x), \alpha_0(s,x))\,dr}\\
&\geq&  \eta\,\, e^{-\eta\, M},
\end{eqnarray*} 
where $M$ is an upper bound of $\lambda$. We obtain the contradiction by letting $\gamma$ go to zero.
\qed

\subsection{Proof of the comparison Theorem \ref{Sec:PDP_Thm_uniqueness}
 }\label{A:Comparison}

We begin recalling the following result concerning an equivalent definition of viscosity super and subsolution to \eqref{Sec:PDP_HJB}-\eqref{Sec:PDP_HJB_T}.
\begin{definition}\label{D:5.4}
	Let assumptions \textbf{\textup{(HL)}}, \textup{\textbf{(H$\textup{b$\lambda$Q}$)}} and  \textup{\textbf{(H$\textup{fg}$)}} be satisfied. We will say that a function $\psi$ is a test function in the sense of Definition \ref{D:5.4} if $\psi(s,y) = \varphi(s,y) + h(||y||)$, where $\varphi$, $h$ are as in Definition \ref{D:testfunc} without being bounded, however $\varphi$ is bounded on every set $(\varepsilon, T-\varepsilon)\times \{x\in E: ||x|| \leq R\}$, $\varepsilon \in (0,T)$, $R >0$.
	\begin{itemize}
 		\item[(i)] A bounded $B$-upper-semicontinuous  function $u:(0,\,T] \times E\rightarrow \R$  is  a \emph{viscosity subsolution in the sense of Definition \ref{D:5.4}} of \eqref{Sec:PDP_HJB}-\eqref{Sec:PDP_HJB_T} if $u(T,x) \leq g(x)$ on $E$, and, whenever $u- \psi$ has a global maximum at a point $(t,x)$ for a test function $\psi(s,y) = \varphi(s,y) + h(||y||)$, then 
 		\begin{align*}
 		&
 		\frac{\partial \psi}{\partial t} (t,x)
 		-\langle x,\,L^\ast \, D \varphi(t,x)\rangle \\
 		& + \inf_{a \in A}\left\{\langle b(x,a),  D \psi(t,x) \rangle  +\int_{E} (u(t,y)-u(t,x))\, \lambda(x, a)\,  Q(x,a, dy) +f(x,a) \right\}  \geq  \,\,0.
 		\end{align*}
 	 	\item[(ii)] A bounded $B$-lower-semicontinuous  function $w:(0,\,T] \times E\rightarrow \R$  is  a \emph{viscosity supersolution in the sense of Definition \ref{D:5.4}} of \eqref{Sec:PDP_HJB}-\eqref{Sec:PDP_HJB_T} if
 		 $w(T,x) \geq g(x)$ on $E$, and, whenever $w+ \psi$ has a global minimum at a point $(t,x)$ for a test function $\psi(s,y) = \varphi(s,y) + h(||y||)$, then 
 		\begin{align*}
 		& 
 		-\frac{\partial \psi}{\partial t}(t,x)
 		 +\langle x,\,L^\ast \, D \varphi(t,x)\rangle \\
 		& + \inf_{a \in A}\left\{\langle b(x,a),  -D \psi(t,x) \rangle  +\int_{E} (w(t,y)-w(t,x))\, \lambda(x, a)\,  Q(x,a, dy) +f(x,a) \right\}  \leq  \,\,0.
 		\end{align*}
 		\item[(iii)] A \emph{viscosity solution} of \eqref{Sec:PDP_HJB}-\eqref{Sec:PDP_HJB_T} \emph{in the sense of Definition \ref{D:5.4}} is a  function which is both a viscosity subsolution and a viscosity supersolution.
 	\end{itemize}	
\end{definition}
\begin{lemma}\label{L:equiv_def_viscsol}
	Let assumptions \textbf{\textup{(HL)}}, \textup{\textbf{(H$\textup{b$\lambda$Q}$)}} and  \textup{\textbf{(H$\textup{fg}$)}} be satisfied. If a function $u: (0,\,T) \times E \rightarrow \R$ (resp. $w: (0,\,T) \times E \rightarrow \R$) is  bounded and uniformly continuous in the $|\cdot| \times ||\cdot||_{-1}$ norm, and is a viscosity subsolution (resp. supersolution) of equation \eqref{Sec:PDP_HJB}-\eqref{Sec:PDP_HJB_T}, then it is a viscosity subsolution (resp. supersolution) of equation \eqref{Sec:PDP_HJB}-\eqref{Sec:PDP_HJB_T} in the sense of Definition \ref{D:5.4}.
\end{lemma}

\noindent \emph{Proof of Lemma \ref{L:equiv_def_viscsol}.}
This lemma extends to the infinite-dimensional framework a well known result in the finite-dimensional case, see e.g. Lemma 2.1 in \cite{SoII}.

We consider the subsolution case, the supersolution case can be proved analogously.
Let thus $u: (0,\,T) \times E \rightarrow \R$ be  bounded and uniformly continuous function in the $|\cdot| \times ||\cdot||_{-1}$ norm, providing a viscosity subsolution to \eqref{Sec:PDP_HJB}-\eqref{Sec:PDP_HJB_T}. Let $u-\psi$ has a global maximum at $(t,x)$ for a test function $\psi(s,y)= \varphi(s,y) + h(||y||)$, where without loss of generality we can assume that $\varphi$ and $h(||\cdot||)$ are bounded and that $u(t,x)= \psi(t,x)$.
By assumption, it exists a modulus $\sigma_u$  such that 
\begin{equation}\label{uniformcontinuity}
|u(s,y) - u(s,z)|\leq \sigma_u(||y-z||_{-1})\quad s \in (0,T), \,\,y,z \in E.
\end{equation}
For  any $\varepsilon>0$, let  $\bar u^{\varepsilon}$ be the sup-inf convolution of  $u$ (see e.g. Definition D.24 in  \cite{gozswi15}), namely 
$$
\bar u^\varepsilon(s,x) = \inf_{z \in E} \sup_{w \in E} \left (u(w) - \frac{||z-w||_{-1}^2}{2 \varepsilon}+ \frac{||z-x||_{-1}^2}{ \varepsilon} \right). 
$$
Then, according to Proposition D.26 in  \cite{gozswi15},  $\bar u^\varepsilon$, $\frac{\partial \bar u^\varepsilon}{\partial_t}$, $D \bar u^\varepsilon$
are uniformly continuous in the $|\cdot| \times ||\cdot||_{-1}$ norm and bounded, and for any $s \in [0,\,T]$, $ y \in E$, 
\begin{align}
&u(s,y) \leq \bar u_\varepsilon(s,y),\label{prop_supinf_1}\\
&|u(s,y) - \bar u^\varepsilon(s,y)| \leq \sigma_u(t_\varepsilon), \label{prop_supinf_2}
\end{align}
where $\frac{t_\varepsilon}{\sqrt{\varepsilon}} \rightarrow 0$ as $\varepsilon \rightarrow 0$. This implies in particular that $\bar u^\varepsilon$, $\frac{\partial \bar u^\varepsilon}{\partial_t}$, $D \bar u^\varepsilon$ and  $A^\ast D \bar u^{\varepsilon}$ are uniformly continuous with respect in the $|\cdot| \times ||\cdot||$ norm.

Let $\eta$ be a smooth function,  such that $\eta(\tau)=1$ for $\tau <1$, $\eta(\tau) = 0$ for $\tau >2$, and which is strictly decreasing on $[1,2]$. We define
$$
\psi^\varepsilon(s,y) := \psi(s,y)\,\eta\left(\frac{||y-x||_{-1}}{\varepsilon}\right)+\bar u^{\varepsilon}(s,y)\left[1-\eta\left(\frac{||y-x||_{-1}}{\varepsilon}\right)\right].
$$
 By definition $u(t,x)-\psi^\varepsilon(t,x)=0$. Moreover
\begin{align*}
u(s,y)-\psi^\varepsilon(s,y)&=u(s,y)\,\eta\left(\frac{||y-x||_{-1}}{\varepsilon}\right)+u(s,y)\left[1-\eta\left(\frac{||y-x||_{-1}}{\varepsilon}\right)\right]\\
&-\psi(s,y)\,\eta\left(\frac{||y-x||_{-1}}{\varepsilon}\right)-\bar u^{\varepsilon}(s,y)\left[1-\eta\left(\frac{||y-x||_{-1}}{\varepsilon}\right)\right]\\
&=(u(s,y)-\psi(s,y))\,\eta\left(\frac{||y-x||_{-1}}{\varepsilon}\right)+(u(s,y)-\bar u^{\varepsilon}(s,y))\left[1-\eta\left(\frac{||y-x||_{-1}}{\varepsilon}\right)\right].
\end{align*}
For all $(s,y) \in [0,\,T] \times E$,   $u(s,y) \leq \bar u^\varepsilon (s,y)$ by \eqref{prop_supinf_1},   and $u(s,y)-\psi(s,y)\leq 0$ by assumption. 

It follows that $u-\psi^\varepsilon$ has a global maximum at $(t,x)$.
Therefore, we apply Definition \ref{Sec:PDP_Def_viscosity_sol_HJB} with $\psi^\varepsilon(s,y) = \varphi^\varepsilon(s,y) + h(||y||)\,\delta^\varepsilon(s,y)$, where 
\begin{align*}
	\varphi^\varepsilon(s,y) &= \varphi(s,y) \,\eta\left(\frac{||y-x||_{-1}}{\varepsilon}\right) + \bar u^{\varepsilon}(s,y)\left[1-\eta\left(\frac{||y-x||_{-1}}{\varepsilon}\right)\right],\\
	\delta^\varepsilon(s,y)&= \eta\left(\frac{||y-x||_{-1}}{\varepsilon}\right).
\end{align*}
Notice  that
$\psi^{\varepsilon}(t,x) = \psi(t,x) = u(t,x)$, $\frac{\partial\psi^\varepsilon}{\partial t}(t,x)= \frac{\partial\psi}{\partial t}(t,x)$, $D \psi^\varepsilon(t,x)=  D\psi(t,x)$.
We get 
\begin{align*}
 		0&\leq  
 		\frac{\partial \psi}{\partial t} (t,x)
 		-\langle x,\,L^\ast \, D \varphi(t,x)\rangle \\
 		& + \inf_{a \in A}\bigg\{f(x,a)  +\langle b(x,a),  D \psi(t,x) \rangle  +\int_{E} (\psi^\varepsilon(t,y)-u(t,x))\, \lambda(x, a)\,  Q(x,a, dy)\bigg\} \\
 		& =  
 		\frac{\partial \psi}{\partial t} (t,x)
 		-\langle x,\,L^\ast \, D \varphi(t,x)\rangle \\
 		& + \inf_{a \in A}\bigg\{f(x,a)  +\langle b(x,a),  D \psi(t,x) \rangle +\int_{E} (u(t,y)-u(t,x))\, \lambda(x, a)\,  Q(x,a, dy) \\
 		&\qquad \quad+\int_{E} (\psi^\varepsilon(t,y)-u(t,y))\, \lambda(x, a)\,  Q(x,a, dy)\bigg\}.
 		\end{align*}
At this point we notice that
 \begin{align*}
&|\psi^\varepsilon(t,y)-u(t,y)|= \left|(\psi(t,y)-\bar u^{\varepsilon}(t,y))\,\eta\left(\frac{||y-x||_{-1}}{\varepsilon}\right)+\bar u^{\varepsilon}(t,y)-u(t,y)\right|\\
&\leq |\psi(t,y)-\bar u^{\varepsilon}(t,y))|\,\eta\left(\frac{||y-x||_{-1}}{\varepsilon}\right)+|\bar  u^{\varepsilon}(t,y)-u(t,y)|\\
&\leq |\psi(t,y)-u(t,y))|\,\eta\left(\frac{||y-x||_{-1}}{\varepsilon}\right) + |\bar u^{\varepsilon}(t,y)-u(t,y)|\left[1+ \,\eta\left(\frac{||y-x||_{-1}}{\varepsilon}\right)\right]\\
& \leq |\psi(t,y)-u(t,y))|\,\eta\left(\frac{||y-x||_{-1}}{\varepsilon}\right) + \sigma_u(t_\varepsilon)\left[1+ \,\eta\left(\frac{||y-x||_{-1}}{\varepsilon}\right)\right], 
\end{align*}
where in the latter inequality we have used \eqref{prop_supinf_2}. 
The conclusion follows by the Lebesgue dominated convergence theorem.

\paragraph{Proof of Theorem \ref{Sec:PDP_Thm_uniqueness}.}\label{Subsec:compthm}

We will show the result by contradiction. Assume therefore that $u \nleq v $.

\medskip 

\noindent \emph{Step 1.} Set 
$u^{\eta}(t,x) =u(t,x) - \frac{\eta}{t}$, $v^{\eta}(s,y) =v(s,y) + \frac{\eta}{s}$, $\eta >0$, 
and, for $\varepsilon, \delta, \beta >0$,  define the function
$$
\Phi^{\varepsilon, \delta, \beta}(t,s,x,y):= u^{\eta}(t,x)-v^{\eta}(s,y)- \frac{||x-y||^2_{-1}}{2 \varepsilon} - \delta (||x||^2 + ||y||^2) - \frac{(t-s)^2}{2 \beta}.
$$
By perturbed optimization (see, e.g. \cite{CrandallLionsV}, page 430) there exist sequences $a_n, b_n \in \R$, $p_n, q_n \in E$ such that 
\begin{equation}\label{parameters}
|a_n|+ |b_n| + |q_n| + |p_n| \leq \frac{1}{n},\quad n \delta \rightarrow \infty, 
\end{equation}
and 
$$
\Phi^{\varepsilon, \delta, \beta}(t,s,x,y) + a_n t + b_n s + \langle B p_n, x\rangle + \langle B q_n, y\rangle
$$
attains a strict maximum at some point $(\bar t, \bar s, \bar x, \bar y) \in (0,\,T] \times (0,\,T] \times E \times E$. 
Standard considerations yield (see e.g.  \cite{gozswi15},  page 209)
\begin{align}
&\lim_{\beta \rightarrow 0} \limsup_{ n \rightarrow \infty} \,\frac{|\bar t -\bar s|^2}{2 \beta}=0, \quad \forall \delta, \varepsilon >0,\label{1}\\
&\lim_{\delta \rightarrow 0}\limsup_{\beta \rightarrow 0} \limsup_{ n \rightarrow \infty} \,\delta(||\bar x||^2 + ||\bar y||^2)=0, \quad \forall \varepsilon >0,\label{2}\\
&\lim_{\varepsilon \rightarrow 0}\limsup_{\delta \rightarrow 0}\limsup_{\beta \rightarrow 0} \limsup_{ n \rightarrow \infty} \,\frac{1}{2 \varepsilon}||\bar x- \bar y||_{-1}^2=0.\label{3}
\end{align}
Then, recalling that  by assumption $u \nleq v $, it follows from \eqref{1}-\eqref{2}-\eqref{3} and the uniform continuity of $u$, $v$, that for sufficiently small  $\varepsilon, \eta, \delta, \beta >0$ and $n$ large enough, $\bar t, \bar s <T$.

\medskip

\noindent \emph{Step 2. 
} 
From  Step 1 we deduce that  
\begin{align*}
&u(t,x)- (\varphi(t,x) + h(||x||)) \,\,\textup{has a global maximum over}\,\,(0,T)\times E\,\,\textup{at}\,\,(\bar t, \bar x),\\
&v(s,y)+ (\psi(s,y) + h(||y||))\,\,\textup{has a global minimum over}\,\,(0,T)\times E\,\,\textup{at}\,\,(\bar s,\bar y),
\end{align*}
where $h(||z||) := \delta ||z||^2$, and 
\begin{align*}
\varphi(t,x) &:= 	\frac{\eta}{t} - a_n t - \langle B p_n, x \rangle  + \frac{||x - \bar y||^2_{-1}}{2\varepsilon} + \frac{(t-\bar s)^2}{2\beta},\\
\psi(s,y) &:= 	\frac{\eta}{s} - b_n s - \langle B q_n, y \rangle  + \frac{||\bar x - y||^2_{-1}}{2\varepsilon} + \frac{(\bar t- s)^2}{2\beta}.
\end{align*}

In particular, $\nabla h(||z||) = 2 \delta z$, and 
\begin{align*}
&\frac{\partial \varphi}{\partial t}(\bar t,\bar x) = - \frac{\eta}{\bar t^2}- a_n  + \frac{\bar t - \bar s}{\beta}, \quad  \frac{\partial \psi}{\partial t}(\bar s,\bar y)  =-\frac{\eta}{s^2}- b_n -\frac{\bar t - \bar s}{\beta},
\\
	&B^{-1} D \varphi(\bar t,\bar x)=- p_n +\frac{\bar x- \bar y}{\varepsilon},\quad 	B^{-1} D \psi(\bar s,\bar y)= - q_n -\frac{\bar x- \bar y}{\varepsilon}.
	\end{align*}

\noindent \emph{Step 3. Viscosity inequalities.} 
 By Lemma \ref{L:equiv_def_viscsol}, 
	$u$ is a  viscosity subsolution of equation \eqref{Sec:PDP_HJB}-\eqref{Sec:PDP_HJB_T} in the sense of Definition \ref{D:5.4}.
Therefore,  
using Step 2, 
we have 
\begin{align}\label{Firstineq_limit}
 &\frac{\bar t - \bar s}{\beta}- \frac{\eta}{T^2}- a_n -\left\langle \bar x,\,L^\ast \, \left(
 \frac{ B (\bar x -\bar y)}{\varepsilon}
 - B p_n\right)\right\rangle + \inf_{a \in A}\Bigg\{\left\langle b(\bar x,a),  \frac{ B (\bar x -\bar y)}{\varepsilon} 
 - B p_n + 2 \delta \bar x \right\rangle \notag\\
 & +\int_{E} (u(\bar t,y)-u(\bar t,\bar x))\, \lambda(\bar x, a)\,  Q(\bar x,a, dy) +f(\bar x,a) \Bigg\}  \geq  \,\,0.
 \end{align}
Similarly, 
being 
	$v$ is a  viscosity supersolution of equation \eqref{Sec:PDP_HJB}-\eqref{Sec:PDP_HJB_T} in the sense of Definition \ref{D:5.4}  by  Lemma \ref{L:equiv_def_viscsol},   proceeding as before one gets
\begin{align}\label{Secondineq_limit}
 &\frac{\bar t - \bar s}{\beta}+ \frac{\eta}{T^2}+ b_n -\left\langle \bar y,\,L^\ast \, \left(\frac{B (\bar x -\bar y)}{\varepsilon}+ B q_n\right)\right\rangle + \inf_{a \in A}\Bigg\{\left\langle b(\bar y,a),  \frac{ B (\bar x -\bar y)}{\varepsilon}+ B q_n - 2 \delta \bar y \right\rangle \notag\\
 & +\int_{E} (v(\bar s,y)-v(\bar s,\bar y))\, \lambda(\bar y, a)\,  Q(\bar y,a, dy) +f(\bar y,a) \Bigg\}  \leq  \,\,0.
 \end{align}
 Subtracting \eqref{Secondineq_limit} from \eqref{Firstineq_limit} we obtain 
 \begin{align}\label{Ineq_limit}
&\frac{2\eta}{T^2}  \leq   -(a_n+b_n) -\frac{1}{\varepsilon}\left\langle (\bar x- \bar y),\,L^\ast \, \left( B (\bar x -\bar y)\right)\right\rangle+\left\langle \bar x,\,L^\ast \, B p_n\right\rangle +\left\langle \bar y,\,L^\ast \, B q_n\right\rangle\\
 & + \sup_{a \in A}\Bigg\{\left\langle b(\bar x,a),  \frac{B(\bar x -\bar y)}{\varepsilon} - B p_n + 2 \delta \bar x \right\rangle -\left\langle b(\bar y,a),  \frac{ B (\bar x -\bar y)}{\varepsilon}+ B q_n - 2 \delta \bar y \right\rangle    \notag\\
 & \qquad \quad +\int_{E} (u(\bar t,y)-u(\bar t,\bar x))\, \lambda(\bar x, a)\,  Q(\bar x,a, dy)\notag\\
 & \qquad \quad -\int_{E} (v(\bar s,y)-v(\bar s,\bar y))\, \lambda(\bar y, a)\,  Q(\bar y,a, dy)  +f(\bar x,a)-f(\bar y,a) \Bigg\},\notag
 \end{align}
where we have used that   $\inf A_1 - \inf A_2 \leq \sup(A_1-A_2)$). Using condition  \eqref{LstarB_property},
 together with the   assumptions on the functions  $b$ and $f$, 
\eqref{Ineq_limit} yields 
 \begin{align}\label{Ineq_limit2}
 & \frac{2\eta}{T^2} +a_n+b_n \leq \langle \bar x,\,L^\ast \, B p_n\rangle +\langle \bar y,\,L^\ast \, B q_n\rangle\notag \\
 &+ \sup_{a \in A}\Bigg\{\frac{1}{\varepsilon}\langle b(\bar x,a)-b(\bar y,a),  B(\bar x -\bar y)\rangle  - \langle b(\bar x,a), B p_n \rangle + \langle b(\bar x,a), 2 \delta \bar x \rangle \notag\\
&\qquad  \quad -\langle b(\bar y,a),  B q_n \rangle +\langle b(\bar y,a),2 \delta \bar y \rangle     +\int_{E} (u(\bar t,y)-u(\bar t,\bar x))\, \lambda(\bar x, a)\,  Q(\bar x,a, dy)\notag\\
 &\qquad \quad -\int_{E} (v(\bar s,y)-v(\bar s,\bar y))\, \lambda(\bar y, a)\,  Q(\bar y,a, dy)  +f(\bar x,a)-f(\bar y,a) \Bigg\}\notag\\
 &\leq \langle \bar x,\,L^\ast \, B p_n\rangle +\langle \bar y,\,L^\ast \, B q_n\rangle \notag\\
 &+ C \left(\frac{||\bar x- \bar y||^2_{-1}}{2\,\varepsilon} + \,(|B p_n| + |B q_n|) +\omega(||\bar x-\bar y||_{-1}) + \delta (1 + ||x||^2 + ||y||^2)\right)\notag \\
  &+ \sup_{a \in A}\Bigg\{\int_{E} u(\bar t,y)\, \lambda(\bar x, a)\,  Q(\bar x,a, dy)-\int_{E} u(\bar t,y)\, \lambda(\bar y, a)\,  Q(\bar y,a, dy)   \Bigg\}\notag\\
   &+ \sup_{a \in A}\Bigg\{\int_{E} (u(\bar t,y)-u(\bar t,\bar x)- v(\bar s,y)+v(\bar s,\bar y))\, \lambda(\bar y, a)\,  Q(\bar y,a, dy)   \Bigg\}.
 \end{align}
At this point, by Hypothesis 
\textbf{(H$\textup{b$\lambda$Q}$)}-(i)-(ii), we get 
\begin{equation*}
\sup_{a \in A}\Bigg\{\int_{E} u(\bar t,y)\, \lambda(\bar x, a)\,  Q(\bar x,a, dy)-\int_{E} u(\bar t,y)\, \lambda(\bar y, a)\,  Q(\bar y,a, dy)\Bigg\} \leq C \omega(||\bar x- \bar y||_{-1}).
\end{equation*}
Therefore it remains to prove that 
\begin{equation}\label{sup_lim}
\sup_{a \in A}\Bigg\{\int_{E} (u(\bar t,y)-u(\bar t,\bar x)-v(\bar s,y)+v(\bar s,\bar y))\, \lambda(\bar y, a)\,  Q(\bar y,a, dy)\Bigg\}
\end{equation}
converges to $0$ when the parameters go to their respective limits. 

\medskip

\noindent \emph{Step 4. Proof of the convergence of \eqref{sup_lim} to $0$.}   Set $m:= 2(||u||_{\infty}\vee ||v||_{\infty})$ and
\begin{equation*}
M:=\Phi^{\varepsilon, \delta, \beta}(\bar t,\bar s,\bar x,\bar y) + a_n \bar t + b_n \bar s + \langle B p_n, \bar x\rangle + \langle B q_n, \bar y\rangle.
\end{equation*}
\noindent By Step 1,  we know that $M$ is a strict maximum on $(0,\,T] \times (0,\,T] \times E \times E$ of the function
$$
\Phi^{\varepsilon, \delta, \beta}(t,s,x,y) + a_n t + b_n s + \langle B p_n, x\rangle + \langle B q_n, y\rangle.
$$ 
\noindent The definition of $\Phi^{\varepsilon, \delta, \beta}$ implies that
\begin{align}\label{eq_missed}
M&=u(\bar t,\bar x)-\frac{\eta}{\bar t}-v(\bar s,\bar y)-\frac{\eta}{\bar s}- \frac{||\bar x-\bar y||^2_{-1}}{2 \varepsilon} - \delta (||\bar x||^2 + ||\bar y||^2) - \frac{(\bar t-\bar s)^2}{2 \beta} + a_n \bar t + b_n \bar s + \langle B p_n, \bar x\rangle + \langle B q_n, \bar y\rangle\notag\\
&=u(\bar t,\bar x)-\frac{\eta}{\bar t}-v(\bar s,\bar y)-\frac{\eta}{\bar s}- \frac{(\bar t-\bar s)^2}{2 \beta} - \frac{||\bar x-\bar y||^2_{-1}}{2 \varepsilon} - \delta ||\bar x-\frac{Bp_n}{2\delta}||^2 - \delta ||\bar y-\frac{Bq_n}{2\delta}||^2\notag\\
&+\frac{||Bp_n||^2}{4\delta}+\frac{||Bq_n||^2}{4\delta}+ a_n \bar t + b_n \bar s,
\end{align}
which in  turn implies that 
\begin{align*}
&\frac{\eta}{\bar t}+\frac{\eta}{\bar s}+\frac{(\bar t-\bar s)^2}{2 \beta}+ \frac{||\bar x-\bar y||^2_{-1}}{2 \varepsilon}+\delta ||\bar x-\frac{Bp_n}{2\delta}||^2+\delta ||\bar y-\frac{Bq_n}{2\delta}||^2\\
&=u(\bar t,\bar x)-v(\bar s,\bar y)-M+\frac{||Bp_n||^2}{4\delta}+\frac{||Bq_n||^2}{4\delta} + a_n \bar t + b_n \bar s.
\end{align*}
Moreover $u(\bar t,\bar x)-v(\bar s,\bar y)\leq m$ and $a_n \bar t + b_n \bar s\leq T$ for all $n\geq 2$, since  $|a_n|+|b_n|\leq \frac{1}{2}$ for $n\geq 2$. Therefore
\begin{align}
\delta ||\bar x-\frac{Bp_n}{2\delta}||^2+\delta ||\bar y-\frac{Bq_n}{2\delta}||^2&\leq m-M+\frac{||Bp_n||^2}{4\delta}+\frac{||Bq_n||^2}{4\delta}+T, \label{forxbarybar}\\
M&\leq m+\frac{||Bp_n||^2}{4\delta}+\frac{||Bq_n||^2}{4\delta}+T.
\end{align}

\noindent Let us take $K \in \N$ satisfying 
\begin{equation}\label{choiceK}
2KM > m + T +\frac{||Bp_n||^2}{4\delta}+\frac{||Bq_n||^2}{4\delta}-M,
\end{equation}and define the set
\begin{equation}\label{gamma1d}
\Gamma_{1,d}:=\left\{ (x,y)\in E\times E; \, \, \, ||x-\frac{Bp_n}{2\delta}||^2+||y-\frac{Bq_n}{2\delta}||^2\leq\frac{2KM}{\delta}\right\}.
\end{equation} 
Notice  that from (\ref{forxbarybar}) we have $(\bar x,\bar y)\in \Gamma_{1,d}$. 
Let also $\alpha>0$ be such that 
\begin{equation}\label{choicealpha}
m+T+\frac{||Bp_n||^2}{4\delta}+\frac{||Bq_n||^2}{4\delta}-2KM+\alpha<M
\end{equation}
\noindent and $D$ be a smooth function on $E\times E$ satisfying 
\begin{align}
D(x,y)=
-\delta \Bigg(||x-\frac{Bp_n}{2\delta}||^2+||y-\frac{Bq_n}{2\delta}||^2\Bigg), \quad \forall (x,y)\in\Gamma_{1,d},\\
-2 KM\leq D(x,y)\leq -2KM+\alpha, \quad \forall (x,y)\in\Gamma_{1,d}^c.\notag
\end{align}
Then the function
\begin{equation}\label{newfuncD}
u^{\eta}(t,x)-v^{\eta}(s,y)- \frac{||x-y||^2_{-1}}{2 \varepsilon} - \frac{(t-s)^2}{2 \beta} + a_n t + b_n s + D(x,y)+\frac{||Bp_n||^2}{4\delta}+\frac{||Bq_n||^2}{4\delta}
\end{equation}
admits a strict maximum at $(\bar t,\bar s, \bar x,\bar y)$. Indeed, if $(x,y)\in \Gamma_{1,d}$ the expression (\ref{newfuncD}) coincides with $
\Phi^{\varepsilon, \delta, \beta}(t,s,x,y) + a_n t + b_n s + \langle B p_n, x\rangle + \langle B q_n, y\rangle$,  and if $(x,y)\notin \Gamma_{1,d}$, by the definition of $D(x,y)$  the expression (\ref{newfuncD}) is smaller or equal to
\begin{align*}
&u^{\eta}(t,x)-v^{\eta}(s,y)- \frac{||x-y||^2_{-1}}{2 \varepsilon} - \frac{(t-s)^2}{2 \beta} + a_n t + b_n s -2KM+\alpha+\frac{||Bp_n||^2}{4\delta}+\frac{||Bq_n||^2}{4\delta}\\
&\leq u^{\eta}(t,x)-v^{\eta}(s,y) + a_n t + b_n s-2KM+\alpha+\frac{||Bp_n||^2}{4\delta}+\frac{||Bq_n||^2}{4\delta}\\
&\leq m+T-2KM+\alpha+\frac{||Bp_n||^2}{4\delta}+\frac{||Bq_n||^2}{4\delta}, 
\end{align*}
the latter being strictly smaller than $M$ by the choice of $\alpha$ (cf. (\ref{choicealpha})).
Using Step 1 with $x=y$ we obtain that,  for all $y\in E$,
\begin{eqnarray*}
\Phi^{\varepsilon, \delta, \beta}(\bar t,\bar s,y,y) + a_n \bar t + b_n \bar s + \langle B p_n, y\rangle + \langle B q_n, y\rangle 
\leq\Phi^{\varepsilon, \delta, \beta}(\bar t,\bar s,\bar x,\bar y) + a_n \bar t + b_n \bar s + \langle B p_n, \bar x\rangle + \langle B q_n, \bar y\rangle,
\end{eqnarray*}
which implies
\begin{align}
&u(\bar t,y)-u(\bar t,\bar x)-v(\bar s,y)+v(\bar s,\bar y)\notag
\\
&\leq- \frac{||\bar x-\bar y||^2_{-1}}{2 \varepsilon}+\delta ||y-\frac{Bp_n}{2\delta}||^2+\delta ||y-\frac{Bq_n}{2\delta}||^2 - \delta ||\bar x-\frac{Bp_n}{2\delta}||^2 - \delta ||\bar y-\frac{Bq_n}{2\delta}||^2\nonumber\\
&\leq\delta \Bigg( ||y-\frac{Bp_n}{2\delta}||^2+||y-\frac{Bq_n}{2\delta}||^2\Bigg)\label{majbis}.
\end{align}

\noindent Let us set 
$$
\Sigma_1:=\left\{ y\in E: \, \, \, ||y-\frac{Bp_n}{2\delta}||^2+||y-\frac{Bq_n}{2\delta}||^2\leq \frac{2KM}{\sqrt\delta}\right\}.
$$ 
For any $y\in {\Sigma_1}$ we obtain by (\ref{majbis}) that
\begin{equation}\label{est_1}
u(\bar t,y)-u(\bar t,\bar x)-v(\bar s,y)+v(\bar s,\bar y)\leq 2KM \, \sqrt{\delta}.
\end{equation}
Let us now set  (since we are interested in $\delta\in (0,1)$ we have $\frac{2KM}{\sqrt\delta}<\frac{2KM}{\delta}$)
$$
\Sigma_2:=\left\{ y\in E:  \, \, \, \frac{2KM}{\sqrt\delta}<||y-\frac{Bp_n}{2\delta}||^2+||y-\frac{Bq_n}{2\delta}||^2\leq\frac{2KM}{\delta}\right\}.
$$
Inequality (\ref{majbis}) yields 
\begin{equation}\label{Gammatilde}
u(\bar t,y)-u(\bar t,\bar x)-v(\bar s,y)+v(\bar s,\bar y)\leq 2KM, \quad \forall y\in \Sigma_2.
\end{equation}
Finally set 
\begin{equation*}
\Sigma_3:=\left\{y\in E: \, \, \, ||y-\frac{Bp_n}{2\delta}||^2+||y-\frac{Bq_n}{2\delta}||^2> \frac{2KM}{\delta}\right\}.
\end{equation*}
Let us now take $y\in \Sigma_3$. Then $(y,y)$ belongs to $\Gamma_{1,d}^c$. From the previous arguments 
\begin{equation*}
u(\bar t,y)-v(\bar s,y)- \frac{\eta}{\bar t} -\frac{\eta}{\bar s}- \frac{(\bar t-\bar s)^2}{2 \beta} + a_n \bar t + b_n \bar s +D(y,y)+\frac{||Bp_n||^2}{4\delta}+\frac{||Bq_n||^2}{4\delta} \leq M,
\end{equation*}
which together with (\ref{eq_missed}) implies 
\begin{eqnarray*}
u(\bar t,y)-u(\bar t,\bar x)-v(\bar s,y)+v(\bar s,\bar y)
\leq  -\frac{||\bar x-\bar y||^2_{-1}}{2 \varepsilon}- \delta ||\bar x-\frac{Bp_n}{2\delta}||^2 - \delta ||\bar y-\frac{Bq_n}{2\delta}||^2-D(y,y).
\end{eqnarray*}
\noindent Thus we obtain
\begin{equation}\label{Gamma2}
 u(\bar t,y)-u(\bar t,\bar x)-v(\bar s,y)+v(\bar s,\bar y)\leq 2KM, \quad \forall y\in\Sigma_3.
\end{equation}
At this point, let us go back to \eqref{sup_lim}. Using the partitioning $E= \Sigma_1  \cup \Sigma_2 \cup  \Sigma_3$,  in \eqref{sup_lim} we split the integral on $E$ in the integrals over the sets $\Sigma_i$.    
From \eqref{est_1} together with  \textbf{(H$\textup{b$\lambda$Q}$)}, we get 
\begin{equation*}
\sup_{a\in A}\int_{\Sigma_1}(u(\bar t,y)-u(\bar t,\bar x)-v(\bar s,y)+v(\bar s,\bar y))\, \lambda(\bar y, a)\,  Q(\bar y,a, dy)\leq 
||\lambda||_\infty \, 2KM\, \sqrt\delta, 
\end{equation*} 
which obviously converges to zero. 
On the other hand, by \eqref{Gammatilde} and \eqref{Gamma2}, we obtain 
\begin{align*}
&\sup_{a \in A}\Bigg\{\int_{\Sigma_2} (u(\bar t,y)-u(\bar t,\bar x)-v(\bar s,y)+v(\bar s,\bar y))\, \lambda(\bar y, a)\,  Q(\bar y,a, dy)\\
&\qquad +\int_{{\Sigma_3}} (u(\bar t,y)-u(\bar t,\bar x)-v(\bar s,y)+v(\bar s,\bar y))\, \lambda(\bar y, a)\,  Q(\bar y,a, dy)\Bigg\}\\
&\qquad \leq 
||\lambda||_\infty \, 2KM   \sup_{a \in A} \,(Q(\bar x,a,\Sigma_2)+Q(\bar x,a, \Sigma_3)).
\end{align*}
We have chosen the parameters according to (\ref{parameters}). Then in particular $\frac{||Bp_n||}{\delta}\leq \frac{1}{n\delta}$ converges to $0$. This completes the proof recalling assumption \textbf{(H$\textbf{Q}$')} (see Section \ref{Sec:PDP_Sec_ConstrainedBSDE}) and the respective definitions of $\Sigma_2, \Sigma_3$.
\qed

\small

\paragraph{Acknowledgements.}
The first author would like to thank Prof. Fausto Gozzi for his helpful
discussions and valuable suggestions to improve this paper.
The  first author   has been financed by "Progetto di Ricerca GNAMPA - INdAM 2018",  and partially benefited
 from the support of the  Italian MIUR-PRIN 2015-16 "Deterministic and
stochastic evolution equations". The financial support of the Laboratoire de Probabilit\'es, Statistique et Mod\'elisation (LPSM, UMR 8001) of Sorbonne Universit\'e is also greatly acknowledged.

\end{document}